\magnification=\magstep 1

\parskip 5pt
\vsize 8.9truein
\global \topskip .3truein
\input epsf.sty
\def \e {\eqno}

\def\bye{\endpage\end}
\def\endpage{\par \vfill \eject}
\font\twelveb=cmbx12

\font\double=msbm10
\def\CC{\hbox{\double\char'103}}

\def\RR{\hbox{\double\char'122}}

\centerline{\twelveb Cauchy Type Integrals of Algebraic Functions}
\centerline{\twelveb and a tangential Center-Focus 
problem for Abel equations} 
\noindent\footnote{}{This research of the last two authors was 
supported by the ISF, Grant No. 264/02, by the BSF, Grant No. 2002243, 
and by the Minerva Foundation.}

\centerline{F. Pakovich, N. Roytvarf and Y. Yomdin}

\centerline{The Weizmann Institute of Science}
\centerline{Rehovot 76100, Israel}
\vskip 1truecm

\noindent {\bf Abstract}
\noindent We consider Cauchy type integrals 
$I(t)={1\over 2\pi i}\int_{\gamma} {g(z)dz\over z-t}$ with $g(z)$ an algebraic
function. The main goal is to give constructive (at least, in principle)
conditions for $I(t)$ to be an algebraic function, a rational function, 
and ultimately an identical zero
near infinity. This is done by relating the Monodromy group of the algebraic 
function $g$, the geometry of the integration curve $\gamma$, and the analytic
properties of the Cauchy type integrals. The motivation for the study of
these conditions is provided by the fact that certain Cauchy type integrals of
algebraic functions appear in the infinitesimal versions of two classical open
questions in Analytic Theory of Differential Equations: the Poincar\'e
Center-Focus problem and the second part of the Hilbert 16-th problem.
 
\vskip 1truecm

\noindent {\twelveb 1. Introduction}

\noindent In this paper we study integrals
$$
I(t)=I(\gamma,g,t)={1\over 2\pi i}\int_{\gamma} {g(z)dz\over z-t}\ , \e(1.1)
$$
where $\gamma$ is a curve in the complex plane $\CC$ and $g(z)$ is an 
algebraic function.
More accurately, we assume that after removing from $\gamma$ a finite
set of points $\Sigma$ (which includes all the double points of $\gamma$ and its
end points) on each segment of $\gamma\setminus \Sigma$ the function $g(z)$
is given by an analytic continuation of a germ of an algebraic function. 
Accordingly, it is always assumed below that the branches of
$g(z)$ on each segment of $\gamma\setminus \Sigma$ are
chosen in advance and in this sense $g(z)$ is univalued on $\gamma$.
Of course, analytic continuation of $g(z)$ outside $\gamma\setminus
\Sigma$ may ramify. Furthermore, we assume that at the
points of $\Sigma$ the function $g(z)$ can ramify but does not have poles.

The main problem considered in this paper is to give conditions 
for the identical vanishing of $I(t)$ near infinity. Since for 
$\vert t\vert\gg 1, $
$$
I(t)= -{1\over 2\pi i}\sum^{\infty}_{k=0}m_kt^{-k-1}
$$
with $m_k=\int_{\gamma}z^kg(z)dz$ this is equivalent to finding conditions
for all the moments $m_k$ to vanish. By the reasons explained below,
we call this problem a ``Moment problem" for integral (1.1). Notice that our
Moment problem overlaps with the classical Moment problems (described for 
example in [1]) only to a rather limited extent. 

The motivation for the vanishing problem for the Cauchy type integrals of
algebraic functions is provided by two classical open   
questions in Analytic Theory of Differential Equations: the Poincar\'e
Center-Focus problem and the second part of the Hilbert 16-th problem
(see Sections 1.1.1 and 1.1.3 below, where the appearance of the
Moment problem in Differential Equations is explained).

As we show below, it is natural to consider the vanishing problem for $I(t)$
together with the conditions for $I(t)$ to be an algebraic or a rational
function.

\medskip 

Recall that for a closed curve $\gamma$ without self-intersection the
vanishing of all the moments $m_k$ is a necessary and sufficient condition
for a function $g$ to be the boundary value of a certain holomorphic function 
$G$ in the compact domain bounded by $\gamma.$ On the other hand, if $\gamma$
is a non-closed curve without self-intersection then the identical vanishing
of $I(t)$ near infinity can happen only for $g(z)\equiv 0$ on $\gamma$
(see, for example, [38]). Nevertheless, the condition of {\it algebraicity} of 
$I(t)$ even in these 
simplest situations turns out to be a rather delicate requirement on the global
monodromy group of $g$ and its local ramifications at the points of $\Sigma$. 
We give this condition and some examples in Section 4 below. 

\medskip

For $\gamma$ a curve with self-intersections the situation
becomes much more complicated already for the vanishing problem. 
Of course, a classical ``homological" condition
for the identical vanishing of $I(t)$ near infinity (i.e. that $g$ 
on $\gamma$ bounds a holomorphic chain) remains valid.
However, this condition does not provide in general a ``constructive"
answer as we would like to have in the Moment problem. Indeed, the
classical homological condition is not easy to translate into an
explicit condition on our finite-dimensional input data (the algebraic
function $g$, the set $\Sigma$, and the homotopy class of $\gamma$).
Even in the simplest case of the ``Polynomial Moment problem'' (see
Section 1.1.1 and Section 6 below) the explicit answer 
is not known in spite of a very
classical setting of the question. Moreover, the results of [39-41] as
well as some of the results of the present paper imply that this answer can 
not be too simple. As we shall see below, the problem of algebraicity of $I(t)$
brings certain additional difficulties related to the ``global'' behavior
of the monodromy of $g$ with respect to $\gamma$.

\medskip

The main goal of the present paper is to give a constructive (at least,
in principle) answer to the Moment and the Algebraicity problems. This is
done by relating the Monodromy group of the algebraic function $g$, the
geometry of the curve $\gamma$, and the analytic properties of the Cauchy
type integrals. Let us describe shortly the main results.  

Let $\CC \setminus \gamma$ be the union of the domains $D_i$ (with $D_0$
being the infinite domain). The expression (1.1) defines $I(t)$ as a collection
of regular analytic functions $I_i$ on the domains $D_i$. A simple
classical description exists for the behavior of $I(t)$ in the process of
crossing the curve $\gamma$: for the adjacent domains $D_i$ and $D_j$ the
function $I_j$ is obtained from $I_i$ by the analytic continuation into $D_j$
combined with the addition of the local branch of $g$ at the crossing point
(also analytically continued into $D_j$). This last operation (as extended to
several crossings of $\gamma$) is ``combinatorial'' in its
nature. It is captured by the
notion of the ``combinatorial monodromy of $I$'' introduced in Section 4.
The combinatorial monodromy depends only on the monodromy of $g$ and on the
geometry of $\gamma$ and in principle it can be explicitly computed.

We analyze the analytic continuation of $I_i(t)$ from each of the domains $D_i$
and show that it is essentially described by the
combinatorial monodromy. On this base we get a necessary and sufficient condition 
for $I(t)$ to be algebraic: its combinatorial monodromy must be finite.
A necessary and sufficient condition for rationality is that the combinatorial
monodromy is trivial. Finally, the vanishing of $I(t)$ is provided by the
additional condition of the absence of the poles in certain sums of the branches
of $g$.

We translate the above conditions into certain local (and local-global) 
branching conditions of $g$ with respect to $\gamma$. Most of these conditions
have a form of the vanishing of a certain sum of the branches of $g$. Besides, we  
give an accurate analytic description of $I(t)$ at its singular points.

We give a number of examples illustrating the above results. Some of them we
consider as rather unexpected. This includes
a non-closed $\gamma$ with nonzero $g$ and $I(t)\equiv 0$ near $\infty$, a non-algebraic
$I(t)$ with a finite ramification of all its branches at each singular point,
and $g$ with non-trivial ``jumps'' on a closed $\gamma$ and with $I(t)\equiv 0$
near $\infty$.

As the first specific application we give an essentially complete solution of 
the ``Rational Double Moment problem'' on the non-closed curve. Remind that
in the case of the closed integration curve the answer is given
by the classical result of Wermer and Harwey-Lawson: double moments vanish
if and only if the path bounds a complex 1-chain (see [2,23,30,59,60] and
Section 1.1.3 below). We show that on a non-closed curve the vanishing of the double moments
(and in fact just an algebraicity of the appropriate generating functions)
is equivalent to a certain composition factorization of the integrand functions
which ``closes up'' the integration path, combined together with the Wermer and 
Harwey-Lawson condition for their ``left factors''.

Another application provided is a significant extension of the class of
``definite'' polynomials (those for which the answer to the Polynomial Moment
problem is given by the Composition condition--see Section 1.1.2 and
Section 5 below). As it is shown in [9,17,18,61], definite polynomials play an important
role in the explicit analysis of the Center-Focus problem for the Abel 
equation. We characterize some classes of definite polynomials $P$ through the
geometry of the images $P(\Gamma)$ of the curves $\Gamma$ joining $a$ and $b$.
This leads also to an interesting geometric invariant of complex polynomials.

In some aspects the present paper provides just the approach to (or the first examples of) 
the phenomena which we expect to be of a major importance in the circle of the problems
considered. This concerns first of all the fact that the Cauchy type integrals of algebraic 
functions satisfy Fuchsian linear differential equations. We do not prove this fact in the
present paper (providing just the idea of the proof in Remark 2 after Theorem 4.4 of
Section 4) but all the necessary tools are prepared here. Another case is the following:
in the present paper apparently new examples appear of a specific
type of functions arising as Cauchy Integrals of algebraic functions: those with
an infinite global ramification but with a finite branching of each of its leaves
at each of the finite number of singular point (see Example 5 of Section 4).
It turns out that such functions are closely related to certain Kleinian groups
and automorphic functions. Once more, in the present paper we restrict ourselves to
a short discussion of this phenomenon in a remark after Example 5 of Section 4, not
providing the proofs. We plan to present separately the rigorous results in these
directions.  

We hope also that the tools introduced in this paper around the notion of the 
combinatorial monodromy can be further developed to provide a really strong approach
for the investigation of the analytic and algebraic properties of the
Cauchy type integrals of algebraic functions. In the remarks at the end of
Section 4 we outline some of the natural directions of such development. We
believe that ultimately it may provide a much deeper understanding of the structure 
of these integrals and of their role in the open questions of the Analytic
Theory of Differential Equations. The present paper is the first step in this direction.

\medskip

The authors would like to thank M. Briskin, A. Eremenko, J.-P. Francoise,
L. Gavrilov, G. Henkin, S. Natanzon, and M. Sodin for inspiring discussions,
and the Max-Planck Institute fur Mathematik, Bonn, where the final 
version of this paper has been prepared, for its kind hospitality.
We would like to thank the referee for a constructive critics which
led to a serious improvement of the paper.
\medskip
\medskip

\noindent {\twelveb 1.1. Motivations}

In this section we discuss some questions in Analytic Theory
of Differential Equation where the Moment problem naturally appears.

\noindent {\twelveb 1.1.1. Classical Center-Focus problem and Moments}

\noindent Our study of the Moment Problem for the Cauchy type
integrals is motivated by the classical Poincar\'e Center-Focus Problem 
for plane polynomial vector fields.

Let $F(x,y)$, $G(x,y)$ be analytic functions of $x$, $y$ in a neighborhood
of the origin in ${\RR}^2$ vanishing at $0$ together with their first derivatives.
Consider the system of differential equations

$$
\cases{\dot x = -y+F(x,y) \cr \dot y=x+G(x,y)} \e(A)
$$

The system (A) has a center at the origin if all its solutions around zero are
closed. The (part of the) classical Center-Focus problem is to find conditions on
$F$ and $G$ which are necessary and sufficient for the system (A) to
have a center at the origin. See [6,7,32,42,46,50,51,53,55,62] for a detailed discussion 
of this problem and of a closely related second part of Hilbert's 16-th problem
(which asks for the maximal possible number of isolated closed trajectories
(limit cycles) of (A)).

It was shown in [20] that one can reduce the system (A) with
homogeneous polynomials $F$, $G$ of degree $d$ to the {\bf trigonometric Abel equation}
$$
r'=p(t)r^2+q(t)r^3, \ \ t \in [0, 2 \pi], \e(B)
$$
where $p(t)$, $q(t)$ are polynomials in $\sin t$, $\cos t$ of the degrees
$d+1$, $2d+2$ respectively. Then (A) has a center if and only if (B) has
all the solutions $r=r(t)$ periodic on $[0,2 \pi]$, i.e. satisfying $r(0)=r(2 \pi)$.
So the classical Center-Focus problem is to find for the trigonometric Abel equations
(B) (obtained via the Cherkas transformation [20]) the necessary and sufficient
condition on $p$ and $q$ for all its solutions $r=r(t)$ to be periodic on $[0,2 \pi]$.
A natural generalization of the
classical Center-Focus problem is to find for any trigonometric Abel equation
(not only for those obtained via the Cherkas transformation [20] from the plane
vector fields) the necessary and sufficient conditions on $p$ and $q$ for all the solutions
of (B) to be periodic.

In turn, the trigonometric Abel equation (B) can be transformed by an exponential
substitution into the  equation
$$
y'=p(x)y^2 + q(x)y^3 
$$
with $p$ and $q$ Laurent polynomials on the
unit circle $S^1$. The Center problem becomes in this setting a problem of 
non-ramifying of all the solutions on $S^1$.

Finally, deviating slightly from the original setting we can consider the Abel
equation
$$
y'=p(x)y^2 + q(x)y^3 \e(C)     
$$
with meromorphic $p,q$ on
any curve, not necessarily closed.

Let $\Gamma$ be a curve in $\CC$ avoiding poles of $p$ and $q$ and joining
two points $a, b \in {\CC}$. The points $a$ and $b$
are called {\it conjugated with respect to (C) along a curve $\Gamma$}
if $y(a)=y(b)$ for any solution $y(x)$ of (C) analytically
continued from $a$ to $b$ along $\Gamma$, with
the initial value $y(a)$ sufficiently small. Equivalently, we shall
say that {\it (C) has a center at $(a, b)$ along $\Gamma$}.
The condition on $p$ and
$q$ under which (C) has a center is called the Center
condition. For $a = b$ this means that the solutions of (C) do not
ramify on the closed curve $\Gamma$. In this case we say that
{\it (C) has a center along $\Gamma$}.

The Hilbert 16-th problem can be reformulated in this setting as follows: we say
that a solution $y(x)$ of the equation (C) (analytically continued
from $a$ to $b$ along $\Gamma$) is {\it periodic} if $y(a)=y(b)$.
The problem is to bound the possible number of periodic solutions of (C).

In both the Center and the Hilbert problems the key tool is the Poincar\'e
first return mapping $G(y_a)=y_b$ which associates to each $y_a$ the value 
$y_b$ of the solution $y(x)$ of the equation (C) satisfying $y(a)=y_a$
and analytically continued from $a$ to $b$ along $\Gamma$. The periodic solutions
of (C) along $\Gamma$ correspond exactly to the fixed points $y_a$ of $G$
(so that $G(y_a)=y_a$) and the Center condition is that $G(y_a) \equiv y_a$. 

\medskip

Although the Center and the Hilbert problems for equation (C) on an interval 
(or in general on a non-closed curve)
do not correspond directly to the classical setting,
they presents an interest by their own and they have been
intensively investigated in [3-5,19-22,25,36,37,54,58] and in many
other publications.
It is a general belief that the Center and Hilbert problems for a polynomial
Abel equation (C) on the interval present all the main difficulties of the
classical ones while possibly simplifying essential technical details.

\medskip
In [3,4] the Center problem for the trigonometric Abel equation has been related
to the composition factorization of the coefficients. Recently in [9-18,21,47,61]
the Center problem for both the trigonometric and the polynomial Abel equations has 
been related to the problem of the vanishing of certain generalized moments on one hand
and to the Composition algebra of univariate analytic functions on the other. This
approach has been further developed in [19,31,58].

For $p$ and $q$ as above let $P= \int p$, $Q= \int q$. Consider the ``one-sided"
moments $m_k(P,Q)$ defined by
$$
m_k(P,Q)=\int_{\Gamma} P^k(x)Q(x)p(x)dx\ , \e(1.2)
$$
and the Moments generating function
$$
H(y)=\int_{\Gamma}{QdP\over y-P} = \sum^{\infty}_{k=0}m_k(P,Q)y^{-k-1}. \e(1.3)
$$
It is shown in [14-16] (see also [9,17,18,61]) that if we consider a parametric
version of the equation (C)
$$
y'=p(x) y^2+{\epsilon}q(x) y^3 \e(D)
$$
then the infinitesimal center conditions with respect to $\epsilon$ for (D)   
at $\epsilon=0$ are given by the vanishing of the one-sided moments. 
Essentially, the Moment generating function $H(y)$ of (1.3) is the
derivative with respect to $\epsilon$ (at $\epsilon = 0$) of the Poincar\'e 
first return mapping $G(y,\epsilon)$ of the equation (D) (see [14-16,17,18]). 
Hence, $H(y)\equiv 0$ is the infinitesimal or the tangential Center condition
(corresponding to the fact that the center of the equation (D) for 
$\epsilon = 0$ survives in first approximation also for nonzero $\epsilon$).

The Moment generating function $H(y)$ of (1.3) defines also the behavior
of the periodic solutions of the equation (D) for small $\epsilon$. 
One can show (see [18]) that this periodic solutions correspond to the zeroes of $H(y)$.
This agrees with the standard fact in the analysis of perturbed plane
Hamiltonian systems, where the derivative with respect to $\epsilon$ 
(at $\epsilon = 0$) of the Poincar\'e first return mapping is given by
the Abelian integrals along the level curves of the Hamiltonian. The 
periodic trajectories of the perturbed system for small $\epsilon$
correspond to the zeroes of these Abelian integrals (see [6,7,27,32,46]).

Consequently, the Moment generating function $H(y)$ plays in the analysis
of the Abel equation the same role as the Abelian integrals play in the
investigation of the perturbed plane Hamiltonian systems, i.e. the most
central role. The study of the conditions for the identical vanishing
of $H(y)$ (which is one of the main problems of the present paper) corresponds
to the study of the identical vanishing of the Abelian integrals in [26]. Essentially,
this is the study of the infinitesimal (or ``tangential'') Center problem for (D).
In turn, the study of the distribution of zeroes of $H(y)$ (which has been recently
started in [18]) corresponds to the study of zeroes of the Abelian integrals.
This last problem is one of the most active research areas in the Analytic Theory
of Differential Equations in the last two decades (see [6,7,26-29,32-35,43,46,56,57]).

Assume now that $p$, $q$ are polynomials. Performing a change of
variables $P(x)=z$, $p(x)dx=dz$ we obtain
$$
m_k(P,Q)
= \int_{\gamma}z^kg(z)dz, \ \ \ \
\ \ \ \ \
H(y)= \int_{\gamma}{g(z)dz\over y-z} = -{2\pi i} I(y)\ ,
$$
with $\gamma=P(\Gamma)$ and $g(z)=Q(P^{-1}(z))$.
Thus the Moment generating function $H(y)$ is a special case of the
Cauchy integral (1.1). As it will be clear below this special case is
not much simpler than the general one. To our opinion this justifies 
a detailed investigation of the general Cauchy-type integrals of
algebraic functions. The relation of $H(y)$ with the Poincar\'e first
return mapping motivates also the study of general analytic properties of 
$I(y)$ (singularities, analytic continuation, etc.) that we start in the
present paper.

\bigskip

\noindent {\twelveb 1.1.2. Center-Focus problem and Compositions}

Let us explain now the role of the Composition Algebra of univalent 
functions in the study of the Center-Focus Problem for the Abel equation
(C) (and in particular in the study of the Moment generating function
$H(y)$). It turns out that a basic sufficient condition for the equation
(C) to have a center (as well as for the vanishing of the one-sided moments
(1.2) and thus for the identical vanishing of $H(y)$) is provided by what
we call a ``Composition condition". Let as above $p=P', q=Q'$. 

\medskip

{\bf Composition condition.} {\it $P(x)=\tilde P(W(x))$, $Q(x)=\tilde Q(W(x))$,
where $W$ maps $\CC$ into a Riemann
Surface $X$ in such a way that $W(\Gamma)$ is a closed curve $\delta$ in $X$
(in particular, if $a\ne b$ then $W(a)$ = $W(b)$), $\delta$ is contained
in a simply-connected domain $D$ in $X,$ and $\tilde P$ and
$\tilde Q$ are regular in $D$.}

A special form of this condition with $X=\CC$ is the Polynomial Composition condition
(PCC) below. In [10,11] (and in Section 1.1.3 below) 
the Composition condition appears with $X$ a rational
curve in ${\CC}^2$. The case of $X$ an elliptic curve is considered in [12].
 
The sufficiency of the Composition condition for the Center problem 
follows from the fact that after performing a change of variables
in (C) and taking $W$ an independent variable, we get a regular
equation in a simply-connected domain $D$ in $X$ whose solutions cannot 
ramify along a closed path $\delta$. The same consideration provides
sufficiency of the Composition condition also for the Moment problem i.e. for
the vanishing of the moments (1.2).

In the case where $p=P', q=Q'$ are polynomials the Composition condition
takes the following simple form which we call the Polynomial Composition
condition (PCC):
$$
P(x)=\tilde P(W(x))\ ,\ Q(x)=\tilde Q(W(x))\ , \eqno(PCC)
$$
with $\tilde P,\tilde Q,W$ - polynomials and with an additional requirement
that $W(a)=W(b)$.

A specific case of the Composition condition has been introduced to the study
of the Center-Focus problem for the Abel equation in [3,4]. The condition (PCC) has
been introduced and intensively studied in [9-18,21,47,61] (see also [19,22,31,58]). 
There is a growing evidence supporting the major role played by the
Polynomial Composition condition (and in general, by the polynomial
Composition Algebra) in the structure of the Center conditions for the
polynomial Abel equation (C). In particular, we have no counterexamples
to the following ``Composition conjecture":

{\bf Composition conjecture.} {\it The Abel equation (C) on the interval with $p,q$
polynomials has a center if and only if (PCC) holds.}

This conjecture has been verified for small degrees of $p$ and $q$ and in many
special cases in [9-19,21,58,61].

\medskip

For some time it was conjectured that (PCC) is a necessary and
sufficient condition also for the Polynomial Moment problem i.e. for the vanishing
of the polynomial moments (1.2) (the ``Moment Composition conjecture''). However,
it was recently shown in [39] that the Polynomial Composition condition (PCC) is
only sufficient but {\it not necessary} for the vanishing of one-sided moments 
(1.2). The counterexample given in [39] exploits rather subtle composition properties
of univariate polynomials (in particular, some classical results of Ritt ([44,45,52]).
This fact stresses the role of the Composition Algebra in the structure of the moments
and of the Center equations and illuminates some important features of the Cauchy-type
integrals of algebraic functions (see Sections 4,5,6 below). 

The appearance of counterexamples to the ``Moment Composition conjecture'' 
(together with the recent results of [9,17,18]) underlines also the role 
of those polynomials $P$ for which the vanishing the one-sided moments (1.2)
with any given $q$ {\it does imply} the Composition condition. Following
[9,61] we call such $P$ {\it definite}. We consider characterization of
definite polynomials as an important problem. Indeed, the role of
definite polynomials in the local Center-Focus problem has been demonstrated
in [9,18,61]. Shortly it can be explained as follows: consider the Abel equation
$$
y'=p(x)y^2 + q(x)y^3, \e(C)
$$
where $p$ is fixed, and $q$ is a variable polynomial in the space $V_d$ of
the univariate polynomials of a given degree $d$. Then for $P=\int p$ definite the
local geometry of the Center set $CS\subset V_d$ of the equation (C) near the
origin is completely described by the Moment vanishing equations. In particular,
this implies that $CS$ near $0\in V_d$ coincides with the Composition set ([9]).
Moreover, it turns out that in this case the so-called Bautin ideal $I$ (generated
by all the Taylor coefficients of the Poincar\'e mapping in the local ring of the
polynomials on $V_d$ near the origin) is in fact generated by the moments. This in turn
provides an important information on the fixed points of the Poincar\'e mapping (i.e.
on the periodic solutions of the Abel equation). See [9,18].
 
Even more important role definite polynomials  
play in the global study of the Center equations  
near infinity (as presented in [17]). One of important conclusions of [17] is the  
following: consider a projectivization $PV_d$ of $V_d$. Then the structure of the 
Center set $CS\subset PV_d$ at and near the infinite hyperplane in $PV_d$ is very
accurately described by the zero set of the ``dual'' moments 
$\tilde m_k = \int Q^k p$. In many cases this description can be extended from
the neighborhood of the infinite hyperplane in $PV_d$ to the affine part $V_d$
and thus it includes the original Center set $CS$.  

One of the applications of the
methods developed in the present paper is a description in Section 5 
below of some apparently new (with respect to the results of [9-16,21,40-41,47,61]) 
and natural classes of definite polynomials.

\bigskip
\noindent {\twelveb 1.1.3. Double moments and centers on a closed curve} 

Here we describe shortly an additional question arising naturally on
the way to the Center-Focus problem. 
This question concerns the center conditions for the Abel equation on a
{\it closed} curve $\Gamma$ and their relation to the vanishing of the
one-sided and double moments. As it was explained above, the Center problem
for the Abel equation on a closed curve provides a rather accurate approximation
of the classical Center-Focus problem. Also in the case of a closed $\Gamma$
the Composition algebra plays a central role connecting the vanishing of
the moments with the center conditions. 

As we pass from an open interval to a closed curve $\Gamma$ the vanishing of the
one-sided moments (1.2) becomes a much less restrictive condition. Although the
Composition condition is still sufficient for this vanishing it is now far
from being necessary. A natural stronger analytic condition on $P$ and $Q$ to be
considered (and compared with the Center one for a closed $\Gamma$) is the
vanishing of the double moments.
 
Generalizing situation a little bit let us consider the double moments
$$
m_{i,j}=\int_{\Gamma} P^i(x) Q^j(x) p(x) dx,  i,j=0, 1, \ldots , \e(1.4)
$$
where $\Gamma$ is a path - closed or non-closed - and $P$ and $Q$ are only
assumed to be holomorphic in a neighborhood of $\Gamma$.

\medskip

The double moments (1.4) appear in several fields of Analysis,
Several Complex variables and Banach Algebras.
The classical result of Wermer and Harwey-Lawson (See [2,23,30,59,60])
implies that {\it if the image of $\Gamma$ under the map $z\rightarrow
(P(z),Q(z))\in {\CC}^2$ is a closed curve $\sigma$ then the vanishing of all the
moments $m_{i,j}$ is equivalent to the fact that the curve $\sigma$ bounds
a compact analytic one-chain in ${\CC}^2$}.

As it was mentioned above the Composition algebra still plays an important 
role in relating the moments vanishing, the Wermer and Harwey-Lawson condition,
and the topology of the rational curve $Y\in {\CC}^2$ parametrized by $P(z),Q(z)$.
In [10] we give a simple and constructive necessary and sufficient
condition for vanishing of double moments (1.4) for $P$ and $Q$
rational functions. It is obtained as a combination
of the composition approach with the general Wermer-Harwey-Lawson theorem.
Assume that the curve $\Gamma$ is closed.
Let $W$ be the right composition greatest common factor (CGCF) 
of $P$ and $Q$, i.e. $P=\tilde P(W)$, $Q=\tilde Q(W)$ with
$\tilde P$, $\tilde Q$ having no composition right factor of a positive
degree. (We shall call such rational functions $\tilde P$, $\tilde Q$ relatively 
prime in composition sense). 

\noindent {\bf Theorem 1.1.3.} {\sl
For $P$, $Q$ rational functions $m_{ij}=0$ for all $i, j \geq 0$ if
and only if all the poles of $\tilde P$ and $\tilde Q$ lie on one
side of $W(\Gamma)$. In particular, if $P$, $Q$ are relatively prime
in composition sense rational functions then  $m_{ij}=0$ for all
$i, j \geq 0$ if and only if all the poles of $P$ and $Q$ lie on one
side of $\Gamma$.}

For $\Gamma$ with self-intersections the ``sides'' of $\Gamma$ are
accurately defined in the next section.

\medskip

The examples of the centers in homogeneous planar systems (A) of degrees 2
and 3 (which lead via the Cherkas transformation to the Abel equations (C)
with $P$ and $Q$ - Laurent polynomials) show that {\it in general on a closed
curve $\Gamma$ the Center condition for the Abel equation (C) with rational
coefficients $p$ and $q$ does not imply the vanishing of the double 
moments} (see [12]).

On the other hand, the vanishing of the double moments implies Center in
many special cases. In particular, in [10,11] it is shown that this is true 
for For $P$ and $Q$ -- Laurent polynomials and $\Gamma=S^1$.
For $P$ and $Q$ - general rational functions the validity of the above 
implication depends on the geometry of the curve $\Gamma$ and its
image $W(\Gamma)$ under the Composition Greatest Common Factor $W$ of $P$ and $Q$.
In particular, if $\Gamma$ and $W(\Gamma)$ are simple closed curves then the
vanishing of the double moments implies Center.

The proof of these facts in [10-11] relays on the general Composition condition given in 
Section 1.1.2 above. The Riemann surface $X$ in our situation is the rational
algebraic curve $Y= (P,Q)(\CC) \subset {\CC}^2$. $P$ and $Q$ provide a rational
parametrization of $Y$ and map $\Gamma$ into the curve $\delta$ in $Y$. The
vanishing of the double moments implies (via the general Wermer-Harwey-Lawson theorem)
that $\delta$ is homological to zero in $Y$. But in many cases a loop $\delta$ on an
affine {\it  rational} complex curve $Y$ which is homological to zero in $Y$ must be
contained in a certain simply-connected domain in $Y$.
By the the general Composition condition this implies Center.
 
In contrast, consider the example of the Abel equation (C) with $P$ and $Q$ the Weierstrass
function and its derivative, respectively, and $\Gamma$ a small circle around zero in $\CC$.
It is shown in [12] that in this case all the double moments vanish while the Abel equation (C)
does not have Center. The difference with the rational case is that for $P$ and $Q$ the 
Weierstrass function and its derivative the curve $Y= (P,Q)(\CC) \subset {\CC}^2$
is an affine {\it elliptic} curve. The loop $\delta = (P,Q)(\Gamma)$ is still homological to  
zero in $Y$. But it bounds in $Y$ a {\it not simply-connected domain}.

All these results stress the role of the vanishing of double moments and of the
Composition factorizations in the case of the closed curve $\Gamma$. Having this in mind we
can now step back and ask for the vanishing condition for the {\it double} moments on the
{\it non-closed} curve $\Gamma$. This question is considered in Section 6 below.

\medskip
\medskip

\noindent {\twelveb 1.2. Organization of the paper}

\noindent {\bf In Section 2} a convenient classical description of the partition of 
$\CC \setminus \gamma$ is given. Then as the initial example we give (by
a direct computation) a description of $I(t)$ for $g$ a rational function.

\noindent {\bf In Section 3} we first remind some classical results on the behavior of
the Cauchy type integrals near the integration curve. Then we obtain an
analytic description of $I(t)$ at the ramification points of $g$ (up to addition
of a regular analytic germ) in terms of the Puiseux expansion of $g$.
As a result we obtain our first necessary condition for algebraicity of
$I(t)$ (which naturally overlaps with the results of Section 4).

\noindent {\bf Section 4} is devoted to the global analytic properties of $I(t)$. 
In particular, here we obtain our main conditions for algebraicity, rationality,
and vanishing of $I(t)$.

We start with the introduction of the notion of the ``sum of branches'' of $g$
across $\gamma$ and along an auxiliary curve $S$. Then we define the notion
of the ``combinatorial monodromy'' and prove that it essentially defines the
usual monodromy of the analytic continuation of $I(t)$.

Next the necessary and sufficient conditions are given for algebraicity,
rationality, and vanishing of $I(t)$ in terms of the combinatorial monodromy.
Later in this section we obtain more explicit local necessary (sometimes also
sufficient) conditions for each of the properties in question.

Finally, a number of examples is considered in reasonable details. 

\noindent {\bf In Section 5} the results of Sections 3 and 4 are applied to
the study of the Polynomial Moment problem. In particular, we produce some new 
classes of definite polynomials. We also introduce and study a certain geometric
invariant of complex polynomials.

\noindent {\bf In Section 6} we study the double moments of rational functions
on a non-closed $\gamma$.
We show that in this case the vanishing of the double moments (in fact, only
algebraicity of the appropriate Moment generating functions)
is equivalent to a certain composition factorization of the integrand functions
which ``closes up'' the integration path, together with the Wermer and Harwey-Lawson
condition for the left factors. The last condition can be interpreted via
Theorem 1.1.3 above. 

%

%
\bigskip
\noindent {\twelveb 2. Partition of $\CC$ by $\gamma$. An example:
$I(g,\gamma,t)$ for $g$-polynomials and rational functions.}

\noindent Below we always assume the curve $\gamma$ to be
oriented, piecewise-smooth, and
to have only transversal self-intersections. In this section we also
assume that $\gamma$ is closed. A classical description of the geometry
of $\gamma$ given below closely follows [2,59,60]. 

The curve $\gamma$ subdivides $\CC$ into a finite number of open domains
$D_i$. One
of these domains which we denote by $D_0$ is unbounded and the rest are
bounded and simply connected. For a point $z\in\CC\backslash \gamma$
define $\mu(\gamma,z)$ as the rotation number of $\gamma$ around $z$.
Clearly, $\mu(\gamma,z)$ is constant on each $D_i$ and we will denote this
constant by $\mu_i$, $\mu_0=0$. Alternatively, $\mu(\gamma,z)$ can be
defined as the (signed) number of the
intersection points of $\gamma$ with any path joining $\infty$ to $z$
or as the linking number of the curve $\gamma$ and the point $z$.
According to this last definition for any complex one-dimensional chain $Z$
in $\CC$ with $\gamma=\partial Z$ the number $\mu(\gamma,z)$ is the
(signed) ``intersection number of $Z$ with $z$'' or, in other words, simply 
the number of times the chain $Z$ covers the point $z$. Figure 1 illustrates 
this construction. 

\medskip
\epsfxsize=12truecm
\centerline{\epsffile{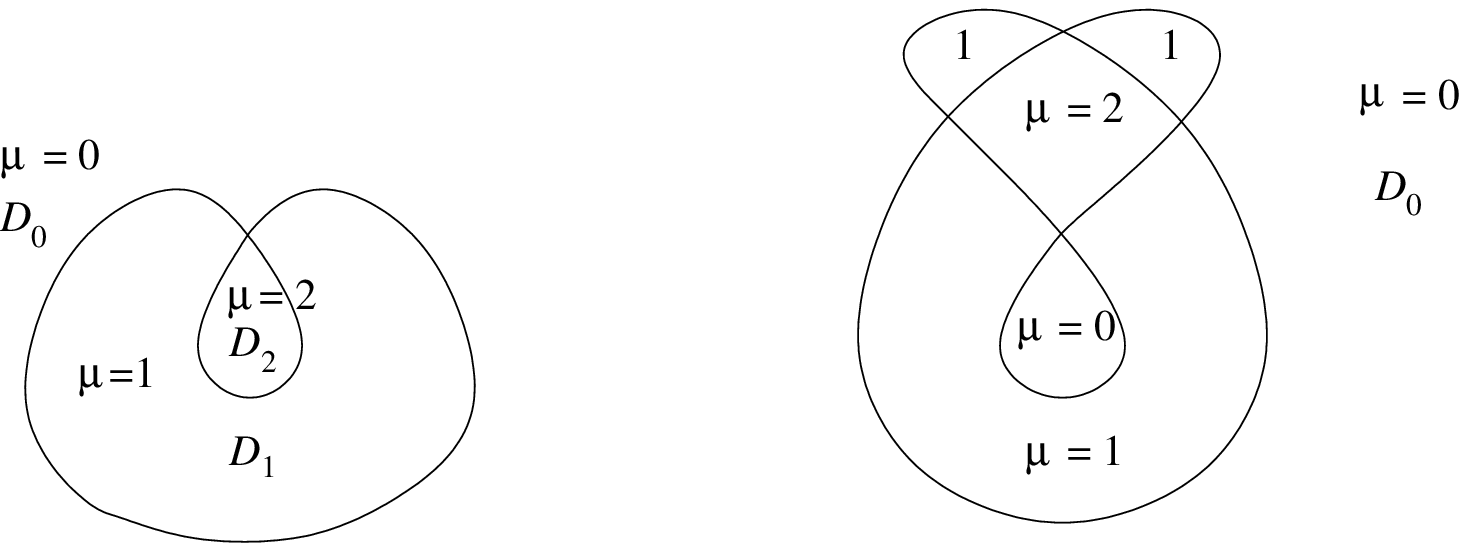}}
\smallskip
\centerline{Figure 1}
\medskip

It is natural to call the union of $D_j$ with $\mu_j=0$ ``an outside
part'' of the curve $\gamma$ and the union of $D_j$ with $\mu_j\ne 0$ ``an
inside part".

Let $\Delta=(z_1,\ldots, z_m)$ be a finite collection of points in $\CC$.

\noindent {\bf Lemma 2.1.} {\sl The curve $\gamma$ is homologous to zero
in $\CC
\backslash\Delta$ if and only if $\mu(\gamma,z_i)=0$ for each $z_i\in
\Delta$, i.e. all $z_i$ belong to the outside part of $\gamma$.}

\noindent PROOF. If $\gamma=\partial Z$ with $Z\subseteq\CC\backslash
\Delta$ then by the last definition of $\mu(\gamma,z)$ this function
is zero for any $z_i\in\Delta$. Conversely, we always have $\gamma=\partial
\left (\sum_{i=0}\mu_i\overline D_i\right )$ and if $\mu_i=0$ on any
domain $D_j$ containing the points of $\Delta$, then the 1-chain 
$\sum_{i=0}\mu_i\overline D_i$ is contained in $\CC\backslash\Delta$.

Now we are ready to describe Cauchy integrals of polynomials and rational
functions on closed curves.

\noindent{\bf Proposition 2.2.} {\sl For $g(z)$ - a polynomial in $z$ 
and for $\gamma$ - closed, $I(\gamma,g,t)=\mu(\gamma,t)\cdot g(t)$.}

\noindent PROOF. A function ${g(z)\over z-t}$ in $z$ has the only pole
at $z=t$ with a residue $g(t)$ while the curve $\gamma$ makes exactly
$\mu(\gamma,t)$ turns around $t$.

Now let $g(z)$ be a general rational function with the poles of the orders
$k_1, \dots ,k_l, k_i \ge 1$ at the points $z_1, \dots ,z_l$, respectively.
Let for $i=1, \dots, l$ the ``essential part'' $R_i(z)$ of the function $g(z)$ at 
$z_i$ (i.e the negative part of the Laurent polynomial of $g(z)$ at the 
pole $z_i$) be given by $R_i(z)=\sum^{k_i}_{j=1} {\alpha_{i,j}\over {(z-z_i)}^j}$.

\noindent{\bf Proposition 2.3.} {\sl The Cauchy integral $I(\gamma,g,t)$ is given by
$$
I(\gamma,g,t)=\mu(\gamma,t)g(t) -\sum^{\ell}_{i=1}\mu(\gamma,z_i)R_i(t). \e(2.1)
$$}
\noindent PROOF. Represent $g(z)$ as 
$g(z)=g_0(z)+\sum^{\ell}_{i=1}\sum^{k_i}_{j=1} {\alpha_{i,j}\over {(z-z_i)}^j}$
with $g_0(z)$ a polynomial. We have
$$
I(\gamma,g,t)=I(\gamma,g_0,t)+\sum^{\ell}_{i=1}\sum^{k_i}_{j=1} \alpha_{i,j}
I(\gamma,{1\over {(z-z_i)}^j},t)\ . \e(2.2)
$$
By Proposition 2.2, $I(\gamma,g_0,t)=\mu(\gamma,t)g_0(t)$. Now, representing the
integrand ${1\over (z-t){(z-z_i)}^j}$ in the Cauchy integral 
$I(\gamma,{1\over {(z-z_i)}^j},t)$ as the sum of the elementary fractions one 
obtains 
$$
{1\over (z-t){(z-z_i)}^j}={({1\over t-z_i})}^j({1\over z-t}-{1\over z-z_i}) + 
\sum^{j}_{s=2} \beta_{i,s} {1\over {(z-z_i)}^s}\ . \e(2.3) 
$$
Integrating along $\gamma$ with respect to $z$ and noticing that the last sum
in (2.3) does not contribute to the integral we get:
$$
I(\gamma,{1\over {(z-z_i)}^j},t) = ({1\over t-z_i})^j
(\mu(\gamma,t)-\mu(\gamma,z_i)). \e(2.4)
$$
Finally,
$$
I(\gamma,g,t) = \mu(\gamma,t)g_0(t)+\mu(\gamma,t)\sum^{\ell}_{i=1}
\sum^{k_i}_{j=1} {\alpha_{i,j}\over {(t-z_i)}^j}-\sum^{\ell}_{i=1}\mu(\gamma,z_i)
\sum^{k_i}_{j=1} {\alpha_{i,j}\over {(t-z_i)}^j} = 
$$
$$
=\mu(\gamma,t)g(t) -\sum^{\ell}_{i=1}\mu(\gamma,z_i)R_i(t).   
$$
This completes the proof of Proposition 2.3.

\noindent{\bf Corollary 2.4.} {\sl Let $g(z)$ be a rational function with
only the first order poles at $z_1,\ldots, z_{\ell}$, each with the residue
$\alpha_i$, $i=1,\ldots,\ell$, respectively.  Then
$$
I(\gamma,g,t)=\mu(\gamma,t)g(t)-\sum^{\ell}_{i=1} {\mu(\gamma,z_i)\alpha_i
\over t-z_i}\ . \e(2.5)
$$}

\noindent{\bf Corollary 2.5.} {\sl If all the poles of $g$ belong to the outside
part of $\gamma$ then $I(\gamma,g,t)=\mu(\gamma,t)g(t)$. In particular, this is
a necessary and sufficient condition for $I(\gamma,g,t)$ to vanish
identically for $t$ near infinity.}

\noindent{\bf Remark.} Computations of Propositions 2.2 and 2.3 did not use in any
essential form the rationality of the function $g$. Exactly in the same way we obtain
the corresponding results for regular and meromorphic functions, respectively. Let
$U \subset \CC$ be a simply-connected 
domain containing the closed integration curve $\gamma$.

\noindent{\bf Proposition 2.6.} {\sl For the integrand $g(z)$ being extendable to a 
holomorphic function in $U$ \ \ $I(\gamma,g,t)=\mu(\gamma,t)\cdot g(t)$.} 

Of course, for such $g(z)$ the Cauchy integral $I(\gamma,g,t)$ always vanishes
identically on the exterior domain $D_0$.

Now let $U$ be as above and let $g(z)$ be a meromorphic 
function in $U$ with the finite number of
poles at the points $z_1, \dots ,z_l$ with the orders $k_1, \dots ,k_l, \ \ k_i \ge 1$,
respectively. 

\noindent{\bf Proposition 2.7.} {\sl The Cauchy integral $I(\gamma,g,t)$ is given by
$$
I(\gamma,g,t)=\mu(\gamma,t)g(t) -\sum^{\ell}_{i=1}\mu(\gamma,z_i)R_i(t),
$$
where for $i=1, \dots ,l$\ \ $R_i(z)$ denotes the essential part of $g(z)$ at the
pole $z_i$. In particular, for such $g(z)$ the Cauchy integral $I(\gamma,g,t)$ is
a rational function on the exterior domain $D_0$.}

\bigskip
\noindent{\twelveb 3. Local structure of $I(t)$}

\noindent Let $\gamma$ be a curve (closed or non-closed) and let $z\in \gamma$.
We say that the integrand function $g$ {\it has a ``jump'' at $z$ if the branches $g_0$
and $g_1$ of $g$ on the two sides of $z$ on $\gamma$ cannot be obtained from one
another by a local analytic continuation (i.e. a continuation along a curve inside
any given neighborhood of $z$). Equivalently, the full local germs of $g_0$ and $g_1$
at a jump point $z$ do not coincide.}

Let us remind that we have denoted by $\Sigma$ the set containing the end-points of 
$\gamma$, all its multiple points, and all the points $z$ on $\gamma$ 
where the integrand function $g$ has either a jump or a ramification point. (In this 
paper we exclude the possibility for $g$ to have poles on $\gamma$).

Lemmas 3.1-3.3 below provide an elementary description of the behavior of the Cauchy
type integral near the integration curve $\gamma$. For a convenience of the reader
we give some proofs and explanations.

Consider first $z_0\in\gamma$ and $z_0\notin \Sigma$. In particular, $g$ is 
regular at $z_0$.

\noindent {\bf Lemma 3.1.} {\sl $I(\gamma,g,t)$ near $z_0$ is represented
by two regular analytic functions: $I_-(t)$ on the left side of $\gamma$
and $I_+(t)$ on the right side. Both $I_-$ and $I_+$ are analytically
extendible into an entire neighborhood $U$ of $z_0$ and $I_+=I_-+g$ in
$U$.}

\noindent PROOF. This is the usual property of Cauchy type integrals (see
[38]) taking into account that $g\vert \gamma$ is a restriction on
$\gamma$ of a regular analytic function $g$ in $U$. Indeed, using analytic
continuation of $g$ into the neighborhood $U$ of $z_0$ and deforming the
integration path we obtain the required extension of $I_-$ and $I_+$ just
by the original expression (1.1).

Let now $z_0$ be a double point of $\gamma$ and let $\gamma_0$ and $\gamma_1$
denote the two local segments of $\gamma$ crossing at $z_0$. We assume that
the restrictions of $g$ to $\gamma_0$ and $\gamma_1$ are both regular at $z_0$.

\noindent{\bf Lemma 3.2.} {\sl $I(\gamma,g,t)$ in a neighborhood $U$ of $z_0$ 
is represented as $I(t)=I^0_{\pm}+ I^1_{\pm}$,
where the combination of the signs is chosen according to the part of
$U\backslash\gamma$ considered, and $I^0_{\pm}$ and $I^1_{\pm}$ have  
with respect to $\gamma_0$ and $\gamma_1$ all the properties stated in
Lemma 3.1.}

\noindent PROOF. Up to a regular addition we have 
$I(\gamma,g,t)=I(\gamma_0,g,t)+I(\gamma_1,g,t)=I^0+I^1$.
Application of Lemma 3.1 to each of these integrals proves the lemma.
Figure 2 illustrates this construction.

\medskip
\epsfxsize=8truecm
\centerline{\epsffile{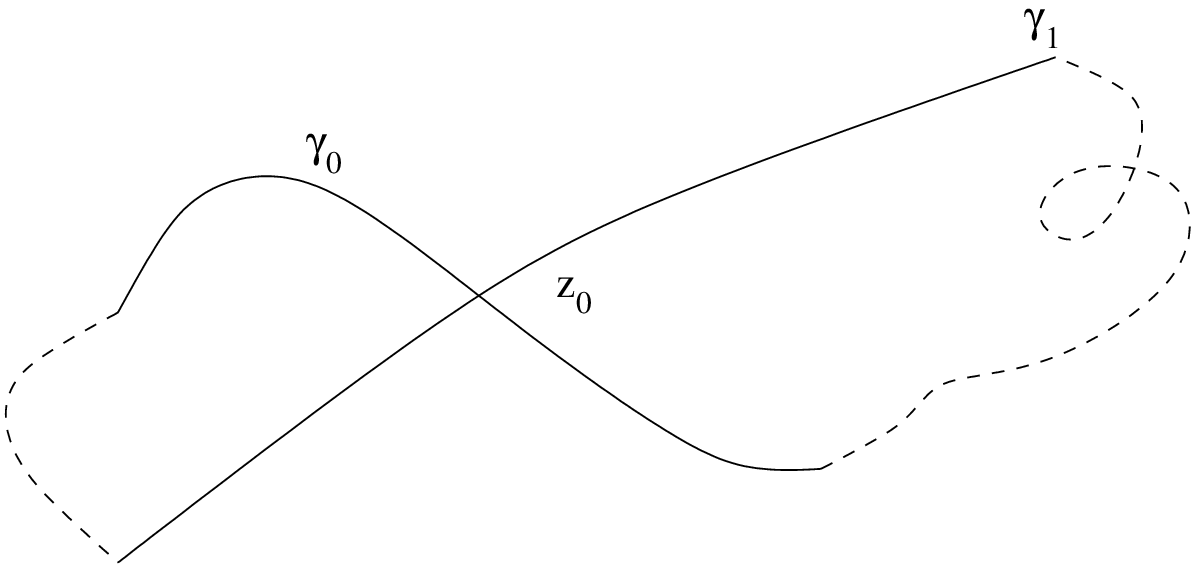}}

\centerline{Figure 2}
\medskip

Assume now that $z_0$ is the end-point of $\gamma$ with the positive
(integration) direction from $z_0$ along $\gamma$. Let $g$ be regular
near $z_0$.

\noindent{\bf Lemma 3.3.} {\sl In a neighborhood of $z_0$,
$$
I(t)=\tilde I(t)-{1\over 2\pi i}g(t)\log (t-z_0)\ 
$$
with $\tilde I(t)$ regular near $z_0$.}

\noindent PROOF. Once more, this is the standard property of the Cauchy type
integrals (see [38]). We give a proof in order to illustrate the approach
used below in a singular situation. 

Up to a regular addition, for $t$ near $z_0$ we can write 
$$
I(t) = {1\over 2\pi i}\int^{z_1}_{z_0} {g(z)dz\over z-t} 
= {1\over 2\pi i}\int^{z_1}_{z_0}{(g(z)-g(t))dz\over z-t} + {1\over 2\pi i}
\int^{z_1}_{z_0} {g(t)dz\over z-t}.
$$
But ${g(z)-g(t)\over z-t}=R(z,t)$ is a regular function of $z$ and $t$ for
$z,t$ near $z_0$. Hence, 
$$
I(t) = {1\over 2\pi i}\int^{z_1}_{z_0}R(z,t)dz +{g(t)\over 2\pi i}
\log {t-z_1\over t-z_0} 
= \tilde I(t)-{1\over 2\pi i}g(t)\log(t-z_0).
$$

The next result describes the local structure of $I(t)$ near the end point
$z_0$ of $\gamma$ which is also a ramification point of $g$. We believe
this result is new. It gives a very accurate description of
$I(t)$ near $z_0$ starting with a very accurate description
of $g$: its Puiseux series at $z_0$. Let
$$
g(z)=\sum^{\infty}_{k=0} a_k(z-z_0)^{k/n} \e(3.1)
$$
be a Puiseux series of $g(z)$ at $z_0$. We denote by $g_r(z)$ a ``regular
part'' of $g(z)$ at $z_0$:
$$
g_r(z)=\sum^{\infty}_{\ell =0} a_{n\ell}(z-z_0)^{\ell}\ . \e(3.2)
$$
Denote by $\tilde g(u)$ a regular function
$$
\tilde g(u) =\sum^{\infty}_{k=0} a_k u^k\ . \e(3.3)
$$
The expression (3.1) does not specify any individual branch of the
algebraic function $g$ at $z_0$, but rather represents the full local germ of
$g$ at $z_0$ (and in this way all its local branches).
However, in the definition (1.1) of the Cauchy type integral a certain
specific branch $g_0$ of $g$ on $\gamma$ near $z_0$ has been fixed. So let $U$ be
a sufficiently small simply-connected neighborhood of a part of $\gamma \setminus z_0$
near $z_0$. Denote by $h(z)$ those branch in $U$ of the
inverse function to ${(z-z_0)}^n$ for which $g_0$ on $\gamma$ near $z_0$ is given by
$$
g_0(z)=\tilde g(h(z)). \e(3.4)
$$
Let $U_0$ be a (simply-connected) part of $U$ lying on one side of 
$\gamma \setminus 0$ (say, the part in the counter-clockwise direction from $\gamma$).
Now we fix the branch $t_0$ of $(t-z_0)^{1/n}$ in $U_0$ given by
$t_0= h(t)$ and denote by $t_j$, $j=0,\ldots, n-1$, all the other
$n$-th roots of $t-z_0$, $t_j=\epsilon^jt_0$, $\epsilon=e^{2\pi i/n}$.
The functions $t_j=t_j(t)$ are univalued functions of $t\in U_0$. 
Remind that in the expression (1.1) of the Cauchy type integral $I(t)$ the
argument $t$ is not allowed to be in $\gamma$.

To formulate our result it remains to fix the branches of the logarithm appearing in
the expression (3.5) below. We fix in an arbitrary way a certain branch of the
$\log (t-z_0)$ for $t$ in the simply-connected domain $U_0$ and denote it by $Log$.
Let $c=h(z_1)$ with $z_1 \in \gamma \setminus z_0$. Fix the branch of the
$\log$ near $c$ satisfying $\log c = {1\over n} Log (z_1-z_0)$ 
for the branch $Log$ of $\log (t-z_0)$
chosen above, and extend it into a small simply-connected neighborhood $V$ of $c$.
We denote this extended branch of the logarithm in $V$ by $Log_1$. 
Finally, for $t$ small also $t_j=t_j(t)$ are small. Therefore $c-t_j \in V$ and we use 
the chosen branch $Log_1$ of $\log(c-t_j)$ in $V$ for each $j=0,\ldots, n-1$.

\noindent{\bf Theorem 3.4.} {\sl For any $t$ in the domain $U_0$,
$$
I(t)= R(t) -\sum^{n-1}_{j=0} ({j\over n}-{1\over 2})\tilde g(t_j) 
+{1\over 2\pi i}\sum^{n-1}_{j=0} \tilde g(t_j)Log_1(c-t_j)-
{1\over 2\pi i}g_r(t)Log (t-z_0). \e(3.5)
$$
Here $R(t)$ is a regular function and $c=h(z_1)$ with $z_1 \in \gamma \setminus 0$
is a nonzero complex constant. The branches $t_j=t_j(t)$ and the branches of the
logarithms are chosen as described above.}
 
\noindent PROOF. To simplify notations, we assume that $z_0=0$. As in the
proof of Lemma 3.3, up to addition of a regular function we have  
$$
I(t)={1\over 2\pi i}\int^{z_1}_0 {g_0(z)dz\over z-t}\ .
$$
Here $z_1 \in \gamma \setminus 0$.
Make a change of variables in this integral: $z=u^n$. Since $\gamma$ is
a smooth curve near $z_0$ the inverse transformation $u=z^{1/n}$ splits $\gamma$
into $n$ smooth curves ${\gamma}_i$ at $0$. Order these $n$ curves in such
a way that ${\gamma}_0$ be the one given by the branch $h(z)$ of $u=z^{1/n}$ defined
above and ${\gamma}_j= \epsilon^j{\gamma}_0$ for $j=0,\ldots, n-1$ (where
as above $\epsilon=e^{2\pi i/n}$). For any $t$ in $U_0$ the values of $t_j=t_j(t)$,
$j=0,\ldots, n-1,$ belong to the simply-connected domains $\Omega_j$ lying
on the ``left'' side of the curves ${\gamma}_j\setminus 0$. See Figure 3. 

\medskip
\epsfxsize=8truecm
\centerline{\epsffile{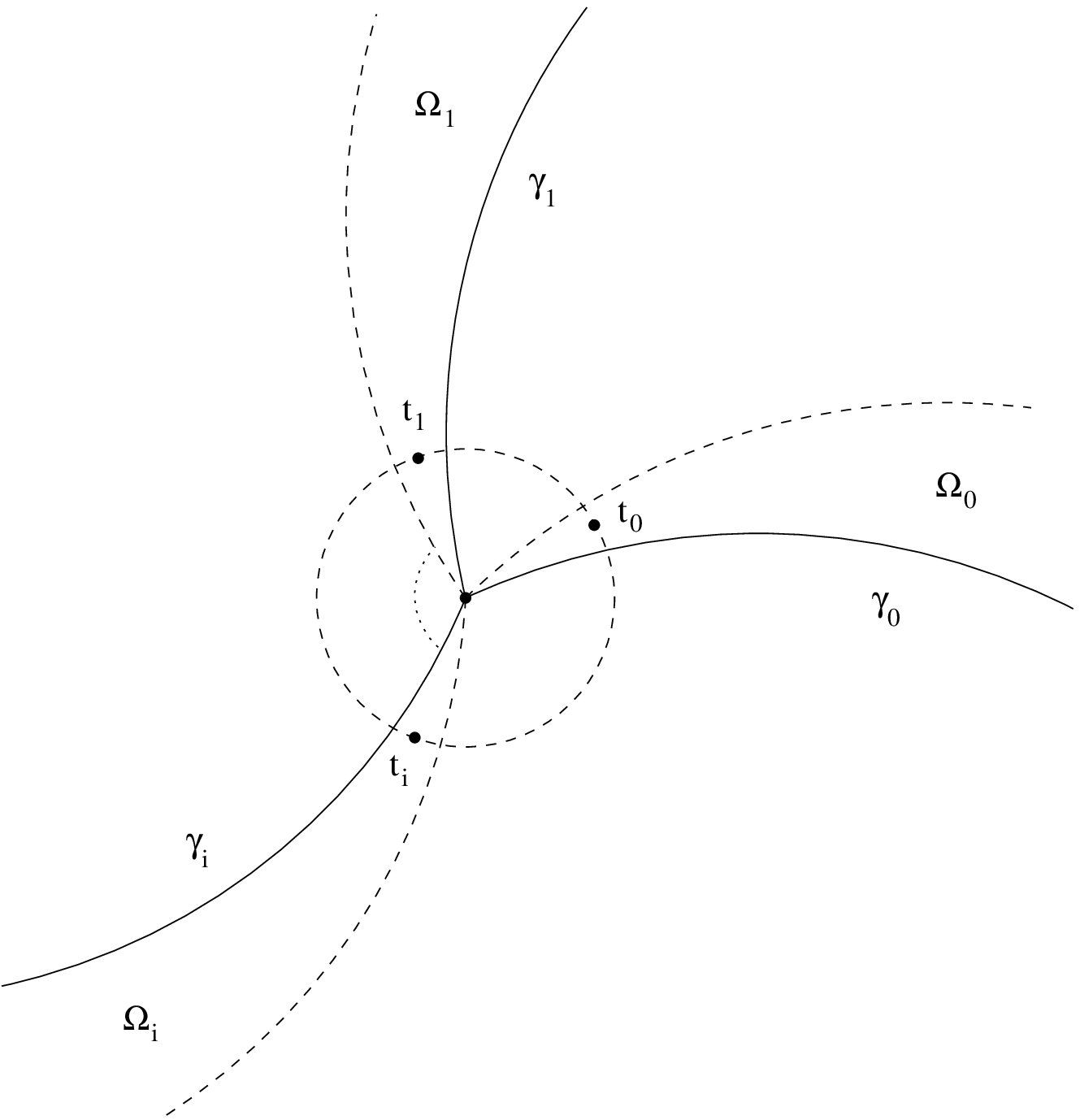}}
  
\centerline{Figure 3}
\medskip

So the integration after the change of variables goes along ${\gamma}_0$ from 
$0$ to $u_1=c=h(z_1)=z_1^{1/n}$. By (3.4) we get  
$$
I(t)={n\over 2\pi i}\int^{u_1}_0 {\tilde g(u)u^{n-1}du\over u^n-t}. \e(3.6)
$$ 
We have:
$$
{1\over u^n-t} =\sum^{n-1}_{j=0} {A_j\over u-t_j(t)}
$$
with $A_j={1\over nt^{n-1}_j(t)}$.  Hence (writing shortly $t_j(t)$ as $t_j$)
we obtain
$$
{2\pi i}I(t) = \sum^{n-1}_{j=0}\int^{u_1}_0{\tilde g(u)u^{n-1}
du\over t^{n-1}_j(u-t_j)} 
= \sum^{n-1}_{j=0} \int^{u_1}_0 \left [{\tilde g(u)
u^{n-1}-\tilde g(t_j)t^{n-1}_j\over t^{n-1}_j(u-t_j)} + {\tilde g(t_j)
\over u-t_j}\right ]du.\e(3.7)
$$
Integration and summation of the terms in (3.7) with only $u-t_j$ in the 
denominator gives (up to the factor ${1\over 2\pi i}$)
$$
\sum^{n-1}_{j=0} \tilde g(t_j)[\log(u_1-t_j)-\log(-t_j)]
= \sum^{n-1}_{j=0} \tilde g(t_j)\log(u_1-t_j)-
\sum^{n-1}_{j=0} \tilde g(t_j)\log(-t_j). \e(3.8)
$$
Now we transform the expression (3.8) taking into account the choice of the
branches of the logarithm fixed in the formula (3.5) of Theorem 3.4.
The differences $\log(u_1-t_j)-\log(-t_j)$ on the left hand side of (3.8) do
not depend on the choice of the branches of the logarithm but do depend on the
integration path (which becomes the continuation path of the logarithm). 
So let us take in each of these differences $\log(u_1-t_j)$ to be given
by the branch of the logarithm $Log_1(u_1-t_j)$ (where $Log_1$ is the branch
in the neighborhood $V$ of $u_1$ fixed in Theorem 3.4).
With this choice the first sum on the right hand side of (3.8) enters as the 
third term into the expression (3.5) of Theorem 3.4.

To compute the second sum on the right hand side of (3.8) (and thus to get the second 
and the fourth terms in (3.5)) let us specify for each summand in (3.8) the 
continuation path of the chosen branch of the logarithm. As $u$ goes from zero to
$u_1$ along ${\gamma}_0$, \ $u-t_j$ goes along the curve ${\gamma}^j_0$
obtained from ${\gamma}_0$ by the shift to $-t_j$. See Figure 4 which depicts
also an auxiliary path $\sigma$ obtained from ${\gamma}_0$ by the 
shift to $t_0$.

\medskip
\epsfxsize=8truecm
\centerline{\epsffile{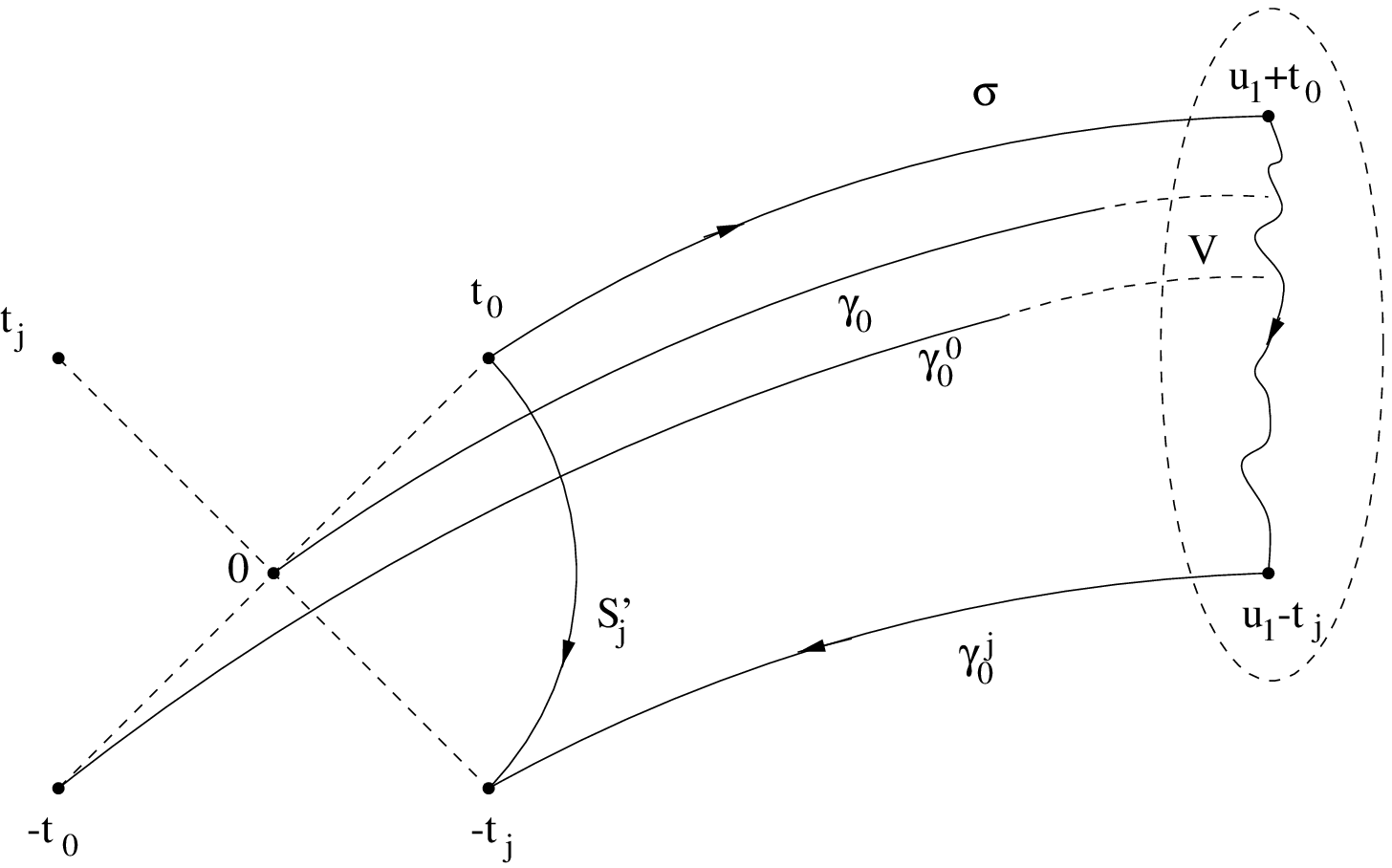}}

\centerline{Figure 4}
\medskip
Notice that the curves ${\gamma}^j_0$, $j=0, \dots ,n-1$, do not pass through $0$.
(Indeed, since $t$ is not allowed to be in $\gamma$, $t_j$ cannot be in any of the branches
${\gamma}_i$, in particular, in ${\gamma}_0$). Some of these curves pass (near
zero) on the same side of ${\gamma}_0$ as $\sigma$ and some pass on the 
other side. An easy computation (see Figure 4) shows that on the same side pass the curves
with $j= [{n-1\over 2}]+1,\dots ,n-1,$ and on the other side those with 
$j= 0\dots ,[{n-1\over 2}].$

Now, according to the integration path in (3.7) the branch of the logarithm in each of the 
summands $\log(-t_j)$ in (3.8) is obtained as follows: the branch $Log_1$ of the logarithm
in a neighborhood $V$ of $u_1$ is taken. Then this branch is continued from $u_1-t_j$
(which belongs to $V$) to $-t_j$ along the path ${\gamma}^j_0$. Notice that by the
choice of the branch $Log_1$ above its continuation to the domain $\Omega_0$ (which
we denote also $Log_1$) satisfies for each $t_0(t) \in \Omega_0$ the equality
$Log_1(t_0)= {1\over n} Log(t)$. 

Therefore we can obtain the branches of each of the summands $\log(-t_j(t))$ in (3.8) by 
continuing $Log_1$ from $t_0$ to $- t_j$ along the path $S_j$ obtained as follows:
we go from $t_0$ back to $u_1 + t_0 \in V$ along the path $\sigma$, then join
$u_1 + t_0$ to $u_1 - t_j$ inside $V$, and then continue from $u_1-t_j$ to $-t_j$ along 
the path ${\gamma}^j_0$. (See Figure 4). Now for $j= 0\dots ,[{n-1\over 2}]$ the path 
$S_j$ can be deformed into the part $S'_j$ of the small circle going from $t_0$ to $-t_j$
in a clockwise direction, while for $j= [{n-1\over 2}]+1,\dots ,n-1$ the path $S'_j$ goes from 
$t_0$ to $-t_j$ in a counter-clockwise direction. In each case the path $S'_j$ presents the
rotation from $t_0$ to the angle $-\pi + ({2\pi \over n})j$ (see Figure 4). Starting with 
the equality $Log_1(t_0)= {1\over n} Log(t)$ we finally obtain
$\log(-t_j)= {1\over n} Log(t)+ ({2\pi i \over n})j -\pi i$.
 
\noindent {\bf Remark.} The same result gives a ``naive'' computation
$$
\log(-t_j)=\log(-\epsilon^jt_0)=\log(-\epsilon^jt^{1/n})= 
{1\over n}\log t + ({2\pi i \over n})j -\pi i.
$$
Hence, the second sum on the right hand side of the expression (3.8) gives
$$
{1\over 2\pi in}Log (t)\sum^{n-1}_{j=0}\tilde g(t_j)+\sum^{n-1}_{j=0}
({j\over n}- {1\over 2})\tilde g(t_j). \e(3.9)
$$
To bring the last term to the form it takes in (3.5)
it remains to notice that the sum $\sum^{n-1}_{j=0}\tilde g(t_j)$ which is
the sum of all the local branches of $g$ at $0$ is equal to $n$ times the
regular part $g_r(t)$ of $g$. Indeed,
$$
\sum^{n-1}_{j=0}\tilde g(t_j) = \sum^{n-1}_{j=0}\sum^{\infty}
_{k=0}a_k(\epsilon^jt_0)^k =
\sum^{\infty}_{k=0} a_kt^k_0 \sum^{n-1}_{j=0}\epsilon^{jk} =
n\sum_{k=n\ell} a_kt^k_0 = ng_r(t). \e(3.10)
$$
This gives the second and the fourth terms in the expression (3.5) of
Theorem 3.4 and completes the computation of the sum of the terms in 
(3.7) with only $u-t_j$ in the denominator.

To complete the proof of Theorem 3.4 it remains to
show that the sum of the terms in (3.7) with the denominator 
$t^{n-1}_j(u-t_j)$ is a regular
function in $t$ near 0. Since $\tilde g(u)$ is given by a convergent power
series in nonnegative integer powers of $u$, it is enough to prove the
statement for each term of this series separately. We have
$$
\sum^{n-1}_{j=0}{u^m-t^m_j\over t^{n-1}_j(u-t_j)} = {1\over t}
\sum^{n-1}_{j=0} t_j(u^{m-1}+u^{m-2}t_j+\cdots +t^{m-1}_j) =
{1\over t}\sum^{n-1}_{j=0} \sum^m_{s=1} u^{m-s}t^s_j = $$
$$={1\over t} \sum^m_{s=1} u^{m-s} \sum^{n-1}_{j=0}t^s_j =
{1\over t}\sum_{s=\ell n,0<s\leq m}u^{m-s}\cdot nt^{\ell} =
n\sum_{0<\ell \leq [m/n]}u^{m-\ell n}t^{\ell -1}.
$$
This completes the proof of Theorem 3.4.

\noindent {\bf Remark 1.} The expression (3.5) of Theorem 3.4 with the branches
$t_j$ and the branches of the logarithms
specified as explained above, represent the actual value of the Cauchy integral
$I(t)$ given by (1.1) in a one-sided neighborhood $U_0$ of
$\gamma \setminus z_0$. However, the same formula gives also the local analytic
continuation of $I(t)$ as $t$ is allowed to vary on the full punctured neighborhood
of $z_0$. In this case the branches of $t_j$ and of the logarithms in (3.5) should
be interpreted as the appropriate analytic continuations of the original ones.

\noindent {\bf Remark 2.} If we fix each of the branches of the logarithms 
appearing in the expression (3.5) of Theorem 3.4 in an arbitrary way, this
expression takes the following form:
$$
I(t)= R(t) -\sum^{n-1}_{j=0} ({j\over n}-{1\over 2}+ m_j)\tilde g(t_j)
+{1\over 2\pi i}\sum^{n-1}_{j=0} \tilde g(t_j)Log_1(c-t_j)-
{1\over 2\pi i}g_r(t)Log (t-z_0),
$$
with $m_j$ -- integers, $j=0, \dots ,n-1$. Indeed, the addition of $(2\pi i)m$ to
$Log (t-z_0)$ brings a regular correction $(2m\pi i) g_r(t)$ which enters the
term $R(t)$. The additions of $(2\pi i)m_j$ to the branches of $Log_1(c-t_j)$
are reflected in the formula above. The computations in the proof of Theorem 3.4
show how to fix the branches of the logarithms in a coherent way in order to
eliminate the unspecified constants $m_j$.

\noindent{\bf Remark 3.} Theorem 3.4 can be considered as a generalization    
and sharpening of computations of C. Christopher [21] and M. Briskin [13].

\noindent {\bf Corollary 3.5.} {\sl In a neighborhood of $z_0$, $I(t)=\hat I
(t)-{1\over 2\pi i}g_r(t)\log(t-z_0)$, where the analytic continuation 
of the function $\hat I(t)$ has
a finite ramification of the order at most $n$ at $z_0$.}

\noindent PROOF. The term $R(t)$ in the expression of Theorem 3.4 is regular
in $t$ near $z_0$. The second and the third terms return to the original
branch as $t$ makes $n$ turns around $z_0$. Indeed, as $t$ is close to
$z_0$, $t_j$ are near 0 but $\log (c-t_j)$ is regular near 0 since
$c\ne 0$. $\tilde g$ is a regular function by definition. Hence as $t_j$ 
return to their original values after $t$ turns $n$ times around $z_0$,
the same do the second and the third terms in the expression of Theorem 3.4.

\noindent{\bf Corollary 3.6.} {\sl If $g_r$ at $z_0$ is not identically
zero then the analytic continuation of $I(t)$ has an infinite ramification 
around $z_0$. In particular, $I(t)$ cannot be an algebraic function.}

\noindent PROOF. For $g_r$ not identically zero the last term in Corollary 3.5
implies an infinite ramification of $I(t)$.

In particular, this happens if $g(z_0)\ne 0$ since in this 
case $g_r(z_0)=g(z_0)\ne 0$.

\bigskip

\noindent{\bf Example 1.} The following simple example illustrates the
statement and the proof of Theorem 3.4 (as well as some of the results
of the next section):

$$
I(t)=\int^1_0 {\sqrt z dz\over z-t}\ .
$$
The positive branch of $\sqrt z$ for $z$ on $[0,1]$ is chosen. Then
$$I(t)=2\int^1_0{u^2du\over u^2-t} =2\int^1_0du+2\int^1_0{tdu\over u^2
-t}.$$ Denote by $\sqrt t$ the branch which is positive for positive $t$
(so $t_0$ and $t_1$ in the proof of Theorem 3.4 are $\sqrt t$ and $-\sqrt t$).
We have
$$
{2t\over u^2-t} =t\left ({1/{\sqrt t}\over u-{\sqrt t}}
-{1/{\sqrt t}\over u+{\sqrt t}}\right )
= {{\sqrt t}\over u-{\sqrt t}}-{{\sqrt t}\over u+{\sqrt t}}\ .
$$
Thus 
$$
I(t)=2-[-{\sqrt t}\log(u-{\sqrt t})+{\sqrt t}\log(u+{\sqrt t})]\vert^1_0=$$
$$ 
= 2+{\sqrt t}\log\left ({1-{\sqrt t}\over 1+{\sqrt t}}\right ) 
-{\sqrt t}(\log(-{\sqrt t})-\log({\sqrt t}))=$$
$$= 2+{\sqrt t}\log\left ({1-{\sqrt t}\over 1+{\sqrt t}}\right ) -{\sqrt t}(\log(-1))
=2-\pi i\sqrt t+ \sqrt t \log \left ({1-\sqrt t\over 1+\sqrt t}\right ).
$$
The term with a logarithmic ramification around zero is absent since
$g(z)=\sqrt z$ has a regular part $g_r(z)$ equal to zero at $z=0$. Notice that
I(t) still has a logarithmic branching around $t=1$.

\medskip

We need also a
description of $I(t)$ near a simple interior point $z_0$ of $\gamma$ at
which $g$ may have a jump. Denote by $g_0$ and $g_1$ the branches of $g$
on $\gamma$ before and after $z_0$, respectively.

\noindent{\bf Theorem 3.7.} {\sl In a neighborhood of $z_0$ functions 
$I_{\pm}(t)$ can be
represented as
$$
I_{\pm}(t)=\hat I_{\pm}(t)+{1\over 2\pi
i}(g_{r_1}(t)-g_{r_0}(t))\log(t-z_0)\ ,
$$
with $\hat I_{\pm}(t)$ having a finite ramification at $z_0$. Here
$g_{r_0}$
and
$g_{r_1}$ denote the regular parts of $g_0$ and $g_1$ at $z_0$ respectively.}

\noindent PROOF. Up to a regular addition the functions $I_{\pm}(t)$ 
in a neighborhood of 
$z_0$ are given by the sum of the Cauchy integrals on the ``semi-curves''
$\gamma_0$ and $\gamma_1$ having $z_0$ as the end points. Now the required
representation follows directly from Corollary 3.5.

\noindent {\bf Remark 3.} One can write in the situation of Theorem 3.7 a
full expression completely analogues to the formula (3.5) of Theorem 3.4.
However, this expression becomes rather complicated since the ramification
orders of $g$ on the two sides of $z_0$ may be different and most of the
sums in (3.5) must appear twice. A simplified version of this formula given
by Theorem 3.7 is sufficient for our applications.

In the same way as Corollary 3.6 we now obtain the following result:

\noindent {\bf Corollary 3.8.} {\sl If $g_{r_0}\not\equiv g_{r_1}$ then
the analytic continuations of both $I_{\pm}(t)$ have an infinite
ramification around $z_0$. In particular, this happens if $g_0(z_0)\ne
g_1(z_0)$.}

The regular part of $g$ at its ramification point has been defined above
in terms of the Puiseax series of $g$ as the sum of all the 
terms in this series with the integer exponent. 
In the proof of Theorem 3.4 it was shown that in fact

$$
\eqalign{\sum^{n-1}_{j=0}\tilde g(t_j) &= ng_r(t)\ .\cr}
$$
In other words, the regular part $g_r(t)$ of $g(t)$ is an average 
(or a ``normalized sum") of all the local branches of $g$.

Let us summarize the results of Corollaries 3.6 and 3.8 as follows:

\noindent{\bf Corollary 3.9} {\sl Let $I(t)$ be algebraic. Then at each
interior point $z_0\in\Sigma$ the regular parts of the branches of $g$ on the  
two sides of $z_0$ coincide. In other words, the normalized sums of the local
branches of $g$ on the two sides of $z_0$ are equal to one another.
If $z_0$ is the endpoint of $\gamma$ then the sum of the local branches of $g$ 
at $z_0$ must be zero.}

Notice that in the statement of the results of Section 3 it is not
essential that the number of the local branches of $g$ is exactly
the denominator $n$ in the exponents of the Puiseux series. Some
of these branches may coincide between themselves - the normalized
sum of the branches remains the same. 

\bigskip
%
\noindent{\twelveb 4. Global structure of $I(t)$: continuation, algebraicity,
rationality, and vanishing}

\noindent Integral representation (1.1) defines $I(t)$ as a collection of
univalent regular functions $I_i(t)$ in each domain $D_i$ of the complement of
$\gamma$ in $\CC$. 
In this section we study the relation between $I_i(t)$
in the neighboring domains $D_i$ and on this base analyze their global analytic
continuation.

Denote by $\gamma_s$ the segments of $\gamma\setminus
\Sigma.$ So $g$ is regular in a
neighborhood of each interior point of $\gamma_s$. 
According to Lemma 3.1 for two adjoint domains $D_i$ and $D_j$
separated by their common segment $\gamma_s$ of the curve $\gamma,$ 
$I_j$ is obtained from $I_i$ as follows:

\item{a.} $I_i$ is analytically continued through $\gamma_s$ into a
certain neighborhood $\Omega$ of $\gamma_s$ in $D_j$.

\item{b.} An algebraic function $g_s$ in $\Omega$ (obtained by the 
analytic
continuation to $\Omega$ of the branch $g_s$ of $g$ on $\gamma_s$) is added to
$I_i$ (multiplied by $-1$ if the crossing orientation of $\gamma$ is
negative).

\item{c.} $I_i+g_s$ is analytically extended from $\Omega$ to the entire
domain $D_j$. $I_j$ is equal to this continued function $I_i+g_s$.

The operation of an addition of an algebraic function and its continuation
(as extended to several crossings of $\gamma$) is ``combinatorial''
in its nature. It depends only on the monodromy of $g$ and on the geometry of $\gamma$
and in principle it can be explicitly computed. To define this operation accurately
let us consider curves $S$ (or $S_{c,d}$)
starting at $c\in \CC \setminus \gamma$ and ending at $d\in \CC \setminus \gamma$. 
Say that $S$ is admissible if it avoids singularities of $g$ and
crosses $\gamma$ transversally and only at the interior points of the segments
$\gamma_s$. For any admissible curve $S_{c,d}$ denote by {\it $S^{*}_{c,d}$ the
operator of the analytic continuation along $S_{c,d}$ of the analytic germs at $c$
to the analytic germs at $d$}.

Now let $S_{c,d}$ be an admissible curve with $c\in D_i$ and $d\in D_j$.
Suppose that $S\cap \gamma =\{a_1,a_2, ... ,a_r\}$ and let
$\{g_1, g_2, ... , g_r\}$ be the germs of $g$ at $a_i$, $1\leq i \leq r.$

Define a {\it sum of branches $g(S_{c,d},\gamma)$ of $g$ along $S_{c,d}$ 
across $\gamma$} as follows: it is a regular algebraic germ at $d$ defined by 
$$
g(S_{c,d},\gamma)= \sum_{i=1}^r sgn(a_i) S^{*}_{a_i,d}(g_i), \e(4.1)
$$
where $S_{a_i,d}$ denotes the part of $S$ 
connecting $a_i$ and $d$ and $sgn(a_i)$ is equal to plus or minus one
according to the orientation of the crossing of $S$ and $\gamma$ at $a_i$. 

The following property of the sum of branches along $S$ is immediate:

\noindent {\bf Proposition 4.1.} {\sl Let the admissible curve $S_{c,d}$ be
the union of the admissible curves $S_{c,e}$ and $S_{e,d}$. Then 
$$
g(S_{c,d},\gamma)=S^{*}_{e,d}(g(S_{c,e},\gamma))+g(S_{e,d},\gamma). \e(4.2)
$$}
Denote by $\Sigma_1$ the set of those singular points of (all the branches) of $g$
that do not lie on $\gamma$ (and hence do not belong to $\Sigma$).
Denote by $U$ the complex plane $\CC$ with $\Sigma$ and $\Sigma_1$ deleted.
By definition, admissible curves $S$ lie entirely in $U$.

\noindent {\bf Proposition 4.2.} {\sl The sum of branches along $S_{c,d}$ depends
only on the homotopy class (with the fixed end-points) of this curve in $U$.}

\noindent PROOF. As we deform $S$ in $U$ preserving the transversality of the
intersection of $S$ and $\gamma$, each term in (4.1) remains the same, being
the analytic continuation of the same function along a continuously deformed
path. Now for a generic deformation of $S$ at each moment of a non-transversal
intersection of $S$ and $\gamma$ a couple of transversal intersections appears
(or disappears). The corresponding terms in (4.1) cancel one another.

\noindent {\bf Remark.} The result follows also from Lemma 4.3 below since the
analytic continuation of $I_i$ depends only on the homotopy class of $S$ in $U$.

The following lemma shows that the sum of branches along $S$ measures
the difference between the analytic continuation $S^{*}(I_i)$ of $I_i$ into the 
domain $D_j$ and the function $I_j$. Let $S_{c,d}$ be any admissible curve 
with $c\in D_i$ and $d\in D_j$.

\noindent {\bf Lemma 4.3.} {\sl The germ of $I_i$ at $c$ can be analytically 
continued along $S$. The resulting germ $S^{*}(I_i)$ at $d$ satisfies
$$
S^{*}(I_i)(t)=I_j(t)-g(S,\gamma)(t). \e(4.3)
$$}

\noindent PROOF. It follows by induction on the number $r$ of the intersection
points of $S$ and $\gamma$. If we write (4.3) in the form
$$
I_j(t)=S^{*}(I_i)(t)+g(S,\gamma)(t) \e(4.4)
$$
then for the first crossing of $S$ and $\gamma$ the equality (4.4) follows
directly from the description of the behavior of $I$ on $\gamma$
given in the steps a, b, c above. Assuming that (4.4) holds after $l$ crossings
of $S$ and $\gamma$ and combining the above description with the definition
of the sum of branches along $S$, we get that (4.4) is valid also after
$l+1$ crossings.

\medskip

Let $S_{c,c}$ be a closed admissible curve with $c\in D_i$. According to
Proposition 4.2, the sum of branches along $S$ across $\gamma$ depends only
on the element $\hat S$ of the fundamental group $\pi_1(U,c)$ defined by $S$. 
We define a {\it combinatorial monodromy of $I_i$ at $c$} as the mapping $A_i$ of
$\pi_1(U,c)$ to the analytic germs at $c$ which associates to each 
$\hat S \in \pi_1(U,c)$ the germ $A_i(\hat S)$ at $c$ equal to the sum of 
branches along $S$ across $\gamma$.

We say that {\it the combinatorial monodromy of $I_i$ at $c$ is finite} if
the image of $\pi_1(U,c)$ under $A_i$ is finite, and we say that
{\it the combinatorial monodromy of $I_i$ at $c$ is trivial} if
$A_i(\hat S)=0$ for any $\hat S \in \pi_1(U,c)$.

The combinatorial monodromy depends only on the monodromy of $g$ and on the
geometry of $\gamma$ and in principle it can be explicitly computed. In the present
paper we mostly restrict ourselves to the local behavior of the combinatorial 
monodromy. See Remarks 1- 3 at the end of this section where we outline a certain
global algebraic approach that captures naturally the combinatorial monodromy and 
simplifies its computation.

Now we are ready to prove the main results of this section. The following Theorem 4.4
provides a description of the complete analytic continuation of the Cauchy integral
$I_i(t)$ from the domain $D_i$ where it was initially defined by the expression (1.1). 
So fix a point $c\in D_i$ and let $I_i(t)$ be the function in $D_i$ defined by (1.1).
Remind that the usual {\it monodromy mapping} $MA_i$ of the fundamental group 
$\pi_1(U,c)$ to the analytic germs at $c$ is given by
$$
MA_i(\hat S)= S^{*}(I_i(t)),
$$
for any $\hat S \in \pi_1(U,c)$.

\noindent{\bf Theorem 4.4.} {\sl The function $I_i(t)$ allows for a complete
analytic continuation as a regular multivalued function $\hat I_i(t)$ in $U$.
For any admissible curve $S_{c,d}$ with $c\in D_i$ and $d\in D_j$ the analytic
continuation $S^{*}_{c,d}(I_i)$ of $I_i$ along $S_{c,d}$ is given by 
$S^{*}_{c,d}(I_i)=I_j-g(S,\gamma)$. In particular, the monodromy mapping $MA_i$
of $I_i(t)$ is given by $MA_i(\hat S)=I_i(t)-A_i(\hat S)$ for any $\hat S \in \pi_1(U,c)$,
where $A_i$ is the combinatorial monodromy of $I_i$.

At the singular points of $g$ in $\Sigma_1$ any leave of the function $\hat I_i(t)$ may have 
only the finite order ramifications. In fact, singularities of $\hat I_i(t)$ at the points of
$\Sigma_1$ are those of certain sums of the branches of $g$, up to a regular addition.

At each point of $\Sigma$ all the leaves of the function $\hat I_i(t)$ have simultaneously
either a finite or an infinite order of the local branching. The analytic representation
of each of the leaves of the function $\hat I_i(t)$ at these points (up to addition of an
algebraic germ) is given by Theorem 3.4.}

\noindent PROOF. The analytic continuation of $I_i(t)$ along any admissible curve
$S$ and its expression via the sum of branches of $g$ across $\gamma$ is provided by 
Lemma 4.3. Applying this expression to the closed curve $S=S_{c,c}$ we get the required
description of the monodromy action $MA_i$.

Now let us take $S$ with the end-point $d$ near a singular point $w_0 \in \Sigma_1$ of 
$g$. By the representation above we get the leave of the function $\hat I_i(t)$ obtained
by the analytic continuation along $S$ as the difference of the regular function $I_j(t)$
and a certain sum of the branches of $g$. This implies the required description of the
singularities of the leaves of $\hat I_i(t)$ at the points of $\Sigma_1$. 

Taking $S$ with the end-point $d$ near a point $w_1 \in \Sigma$ we get the corresponding
leave of $\hat I_i(t)$ as the difference of the possibly singular at $w_1$ function
$I_j(t)$ and a certain finite sum of the branches of $g$. Therefore the property of this
leave to have a finite or an infinite local ramification at $w_1$ depends only on the
local branching of $I_j(t)$. Hence all the leaves simultaneously have at $w_1$ either
a finite or an infinite order of the local branching. Since the analytic representation
of $I_j(t)$ at $w_1$ is given by Theorem 3.4, this provides the required analytic 
description of the singularities of of all the leaves of $\hat I_i(t)$ at the points of
$\Sigma$, up to addition of an algebraic germ. This completes the proof of Theorem 4.4.

\medskip

\noindent{\bf Remark 1.} Simple examples show that the combinatorial monodromy
of $I_i$ may be infinite. This happens for instance in Example 1 of Section 3 and
in Examples 2 and 5 below. In this case we still get certain finite sums of
branches of $g$ on each step of forming the sum of branches across $\gamma$.
However, it is exactly this step of forming sums that may lead ultimately to infinite
branching. 

\noindent{\bf Remark 2.} Let us return for a moment to the comparison between the
functions $I(t)$ and the Abelian integrals. One of the most important
analytic properties of the the Abelian integrals is the fact that they satisfy
certain Fuchsian linear differential equations with rational coefficients (see 
[6,7,26-29,32-35,43,46,56,57]). Theorem 4.4 allows us to show that the same is true
for the Cauchy type integrals of algebraic functions. Indeed, a necessary and sufficient
condition for a multivalued function $y(x)$ to satisfy a linear differential equation
with univalued coefficients is that the linear space spanned by the germs of all the
branches of $y(x)$ at each point $x$ is finite dimensional. Now by Theorem 4.4 all the
branches of the analytic continuation of $I_i(t)$ over the domain $D_j$ are obtained 
by adding to $I_j(t)$ certain sums of the branches of the algebraic function $g$.
Hence the linear space spanned by the germs of all the branches of the analytic
continuation of $I_i(t)$ over the domain $D_j$ always has the basis consisting of
$I_j(t)$ and of all the branches of $g$. More detailed analysis allows one to show
that $I(t)$ satisfies in fact a Fuchsian linear differential equations with rational 
coefficients. It would be important to construct such an equation explicitly and to
investigate its relation to the differential equation which is satisfied by the 
algebraic function $g$ itself. 

\noindent{\bf Remark 3.} There is another natural way to compute the monodromy of
the Cauchy type integrals of algebraic functions. As the argument $t$ of the function
$I_i(t)$ follows the loop $S$ and approaches $\gamma$ we start to deform $\gamma$ in 
order to avoid the crossing of $\gamma$ by $t$. Under these circumstances the integral 
expression (1.1) on the deformed curve $\gamma$ defines the analytic continuation of 
$I_i(t)$. After $t$ completes the full loop $S$, the integration contour $\gamma$
must be modified to ${\gamma}'$ by adding certain loops related to $S$. The analytic
continuation $S^{*}(I_i(t))$ is given by the Cauchy integral (1.1) over ${\gamma}'$.
Of course, the explicit computation of the resulting additions leads to the same 
formula with the sum of the branches of $g$ as in Theorem 4.4.

We would like to thank L. Gavrilov for pointing out the facts and questions mentioned
in Remarks 2 and 3. 

\noindent{\bf Remark 4.} The results above provide also a comparison of the analytic
continuations of the functions $I_i$ and $I_j$ from two different domains $D_i$ and $D_j$:
these continuations differ by a certain sum of the branches of $g$.

Theorem 4.4 allows us to give a necessary and sufficient condition for $I_i(t)$
to be an algebraic function. Let us remind shortly some basic definitions. 
A (multivalued) function $y(z)$ is called {\it algebraic} if it satisfies an
equation $y^d + a_{d-1}(z)y^{d-1} + \dots + a_1(z)y + a_0(z)= 0$ with $a_i(z)$ 
rational functions. A singular point $z_0$ of a multivalued analytic function $f(z)$ 
is called {\it algebroid} if $f(z)$ has a finite ramification at $z_0$ and an
absolute value of $f(z)$ near $z_0$ is bounded by a certain negative power of
$\vert z-z_0 \vert$. A multivalued analytic function $f(z)$ is called {\it locally 
algebroid} if all its singularities are algebroid. As we show in a moment, under some 
additional finiteness assumptions locally algebroid functions are in fact algebraic.

Let $V= \CC \setminus \{z_0, \dots , z_N\}$. Let $f(z)$ be a multivalued
regular analytic function in $V$. We say that {\it $f(z)$ has a globally finite branching
in $V$} if the number of the different univalued branches of $f(z)$ over each
simply-connected domain $\Omega$ in $V$ is finite. (This number does not depend on $\Omega$).
A basic fact here is the following: {\it a locally algebroid function $f(z)$ which has a
globally finite branching in $V$ is algebraic}. Let us give for completeness a very short proof
of this fact. Let $y_1(z), \dots , y_d(z)$ be the values of all the branches of
$f(z)$ at $z \in V$ (taken in any order). Consider the product
$(y-y_1(z)) \dots (y-y_d(z))= y^d + a_{d-1}(z)y^{d-1}+ \dots + a_1(z)y+ a_0(z)$.
The coefficients $a_{d-1}(z), \dots , a_0(z)$ are symmetric functions of
$y_1(z), \dots , y_d(z)$ so they are univalued on $V$. Since the singularities of all
the branches $y_1(z), \dots , y_d(z)$ of $f(z)$ at the points $z_0, \dots , z_N$ are
algebroid the same is true for $a_{d-1}(z), \dots , a_0(z)$. Now, {\it a univalued
function on $V$ with all the singularities at $z_0, \dots , z_N$ algebroid is rational.}
Indeed, if we multiply such a function by a polynomial with zeroes of a sufficiently
high order at all the finite singularities, we get a regular function on $\CC$ with
a polynomial growth at infinity, i.e. another polynomial. 

Now we are ready to give a criterion for algebraicity of $I_i(t)$.

\noindent{\bf Theorem 4.5.} {\sl The function $I_i(t)$ (and hence all the $I_j(t)$)
is algebraic if and only if the combinatorial monodromy of $I_i$ is finite.}   

\noindent PROOF. Assume that $I_i(t)$ is algebraic (which is the same as to assume
that its full analytic continuation $\hat I_i(t)$ is algebraic). In this case the image
of the monodromy mapping $MA_i$ is finite. By Theorem 4.4 the same is true for the
combinatorial monodromy of $I_i$. Conversely, assume that the combinatorial monodromy
of $I_i$ is finite. Then by Theorem 4.4 the monodromy mapping $MA_i$ has a finite image.
Hence $\hat I_i(t)$ has globally only a finite branching. The singularities of all the 
branches of $\hat I_i(t)$ at the points of $\Sigma_1$ are algebraic (up to a regular addition)
by Theorem 4.4. In our case of the finite branching 
the singularities of all the branches of $\hat I_i(t)$
at the points of $\Sigma$ are algebroid. Indeed, by Theorem 4.4 each of these singularities
at a certain point $z_0 \in \Sigma$ is given by the sum of the algebraic germ and the germ 
of the Cauchy integral $I_j$ at $z_0$,
where $D_j$ is the domain adjacent to $\gamma$ near $z_0$. Since $\hat I_i(t)$ has a finite
branching at $z_0$ the same is true also for $I_j(t)$ at $z_0$. Consider for simplicity the
case of the endpoint $z_0$ of $\gamma$. The case of the interior jump point $z_0$ is treated
exactly in the same way.
 
According to the representation of $I_j(t)$ given by Theorem 3.4, it has the form
$$
I(t)= R(t) -\sum^{n-1}_{j=0} ({j\over n}-{1\over 2})\tilde g(t_j)
+{1\over 2\pi i}\sum^{n-1}_{j=0} \tilde g(t_j)Log_1(c-t_j)-
{1\over 2\pi i}g_r(t)Log (t-z_0). \e(4.5)
$$
In the case of the finite branching of $I_j(t)$ the last term in (4.5) containing
$\log (t-z_0)$ must disappear. Now, the first three terms of (4.5) are bounded
near $z_0$ and have at this point a finite ramification (remind that the constant
$c$ in the third term is different from zero). Hence the singularity
of $I_j(t)$ at $z_0$ is algebroid. Therefore the same is true for the singularity of 
$\hat I_i(t)$ at $z_0$. 

Also at infinity the singularity of any branch of $\hat I_i(t)$
is algebroid. Indeed, $\hat I_i(t)$ has there a finite branching. Moreover, 
up to addition of an algebraic germ at infinity any branch of 
$\hat I_i(t)$ coincides there with the regular germ
$I_0(t)$ given by (1.1). Thus $\hat I_i(t)$ has globally a finite branching and all its
singularities algebroid.  By the classical description of algebraic functions given above
this implies that $\hat I_i(t)$ is algebraic. This completes the proof of Theorem 4.5.

\medskip

Using the same approach we obtain some {\it local} conditions on $g$ that are necessary
for algebraicity of $I(t)$. Below $\gamma$ may be open or closed.

\noindent{\bf Corollary 4.6.} {\sl Let $I(t)$ be algebraic. Then at each
``jump" point $z_0\in\Sigma$ (including each of the endpoints of $\gamma$)
a certain nontrivial integer linear combination of the local branches 
of $g$ must be zero.} 

\noindent PROOF. Consider the case of the endpoint. Let $z_0$ be one of the
endpoints of $\gamma$ and let the adjacent domain to $\gamma$ near $z_0$ be $D_i$.
Denote by $\sigma$ a small closed loop going around $z_0$ in the counter-clockwise
direction from a certain point $c$ near $\gamma$, and let 
$\sigma^{*}$ denote as above the operator 
of the analytic continuation along $\sigma$. Let $g_1$ be the branch of $g$ on $\gamma$
near $z_0$. By Proposition 4.1 the sum of the branches along $\sigma$ repeated
$n$ times is given by
$$
g_1 + \sigma^{*}(g_1)+{\sigma^{*}}^2(g_1)+ \dots +{\sigma^{*}}^n(g_1). \e(4.5)
$$
Now $I_i(t)$ is algebraic and its ramification at $z_0$ is finite. Fix the smallest
$n$ for which ${\sigma^{*}}^n(I_i(t)=I_i(t)$. Then the formula (4.3) of Lemma 4.3
shows that for this $n$ the sum (4.5) is zero. This provides the required
relation. The proof for the interior point $z_0$ is essentially the same. We get an
equality of between certain sum of the branches of $g$ on the two sides of $z_0$ on
$\gamma$. Since $z_0$ was assumed to be a ``jump" point of $g$, the branches of $g$
on the two sides of $z_0$ on $\gamma$ cannot be transformed one into another by the
local monodromy of $g$ and hence the resulting sum of the branches of $g$ is non-trivial.
(We would like to thank the referee for suggesting the above calculation).

\noindent{\bf Remark 1.} Corollary 3.9 is formally stronger than Corollary 4.6 since
it provides {\it the specific vanishing sum} of the local branches of $g$. A modification
of the arguments above allows one to get the same specific vanishing sum via the
approach of this section. 

\noindent{\bf Remark 2.} Example 5 below shows that in general the vanishing of the sums
of the local branches of $g$ given by Corollary 3.9 {\it does not imply the global 
finiteness} of the combinatorial monodromy. One can hope that the algebraic approach to the
representation and computing the combinatorial monodromy described in the remark at the end
of this section can provide a unified way to producing all the necessary and sufficient
finiteness conditions in terms of the vanishing of certain sums of the branches of $g$
(local and global). 

\medskip

Let us continue by providing conditions for $I_i(t)$ to be a rational function.
Notice that in contrast with the algebraicity case these conditions depend
on the specific function $I_i(t)$ (and the domain $D_i$) we start with. Indeed, if
$I_i(t)$ is rational then the functions $I_j(t)$ for $i\ne j$ will be usually only
algebraic and not rational (unless $g$ itself is rational).

\noindent{\bf Theorem 4.7.} {\sl $I_i(t)$ is a rational function if and only
if the combinatorial monodromy of $I_i$ is trivial.}

\noindent PROOF. Assume that $I_i(t)$ is rational. This is the same as to say that
its full analytic continuation $\hat I_i(t)$ is rational.
In this case the image of the monodromy mapping $MA_i$ is the germ $I_i(t)$.
Hence by Theorem 4.4 the image of the combinatorial monodromy of $I_i$
is zero. Conversely, assume that the combinatorial monodromy of $I_i$ is trivial.
Then by Theorem 4.4 the image of the monodromy mapping $MA_i$ is the germ $I_i(t)$.
In particular, this implies that $\hat I_i(t)$ is univalued over $U$.
It was shown in the proof of Theorem 4.5 that all the singular points of $\hat I_i(t)$
(including infinity) are algebroid. By the basic result on rational functions presented 
above in this case $\hat I_i(t)$ must be rational.

Of course, the local conditions of Corollary 4.6 are satisfied if $I_i(t)$ is rational.
Let us present some stronger local conditions for rationality (which in contrast to
the algebraicity local conditions are both necessary and sufficient). These conditions
just express the non-branching of $\hat I_i(t)$ at each of its singularities. We consider 
separately three cases: the endpoint of $\gamma$, the interior ``jump'' point, and
the singular point of $g$ not in $\gamma$.

Let $z_0$ be the endpoint of $\gamma$, belonging to the domain $D_j$, and let
$S_{c,d}$ be an admissible curve with $c\in D_i$ and with $d\in D_j$ close to
a certain point $w$ on the curve $\gamma$. (See Figure 5). The sum of branches
$g(S_{c,d},\gamma)$ of $g$ (along $S$ across $\gamma$) is an algebraic germ at 
$d$ which we extend to the algebraic function $F$ (multivalued in general) 
defined in $D_j$. Denote by $\sigma$ a small closed loop going around $z_0$
in the counter-clockwise direction from the point $w$ to itself. Finally,
let $g_1$ be the branch of $g$ on $\gamma$ near $z_0$.

\noindent {\bf Proposition 4.8.} {\sl If $I_i(t)$ is a rational function then
in a neighborhood of $w$ we have $g_1= F -{\sigma^*}(F)$. In particular, if $D_j=D_i$ 
then $F\equiv 0$ and therefore $g_1\equiv 0$ on $\gamma$.}
  
\noindent PROOF. Denote by $\tilde S$ the curve following $S$ and then going from
$d$ to $w$ and let $T=T_{c,w}$ be the curve following $\tilde S$ and then following
$\sigma$. If $I_i(t)$ is rational then by Lemma 4.3 the sum of branches along
any two admissible curves leading to the same point $w$ must be the same. Hence
$g(T,\gamma)=g(\tilde S,\gamma)=F$. Applying Proposition 4.1 we obtain:
$g(T,\gamma)=g_1+{\sigma^*}(F)$ or $g_1=F - {\sigma^*}(F)$.

\medskip
\epsfxsize=8truecm
\centerline{\epsffile{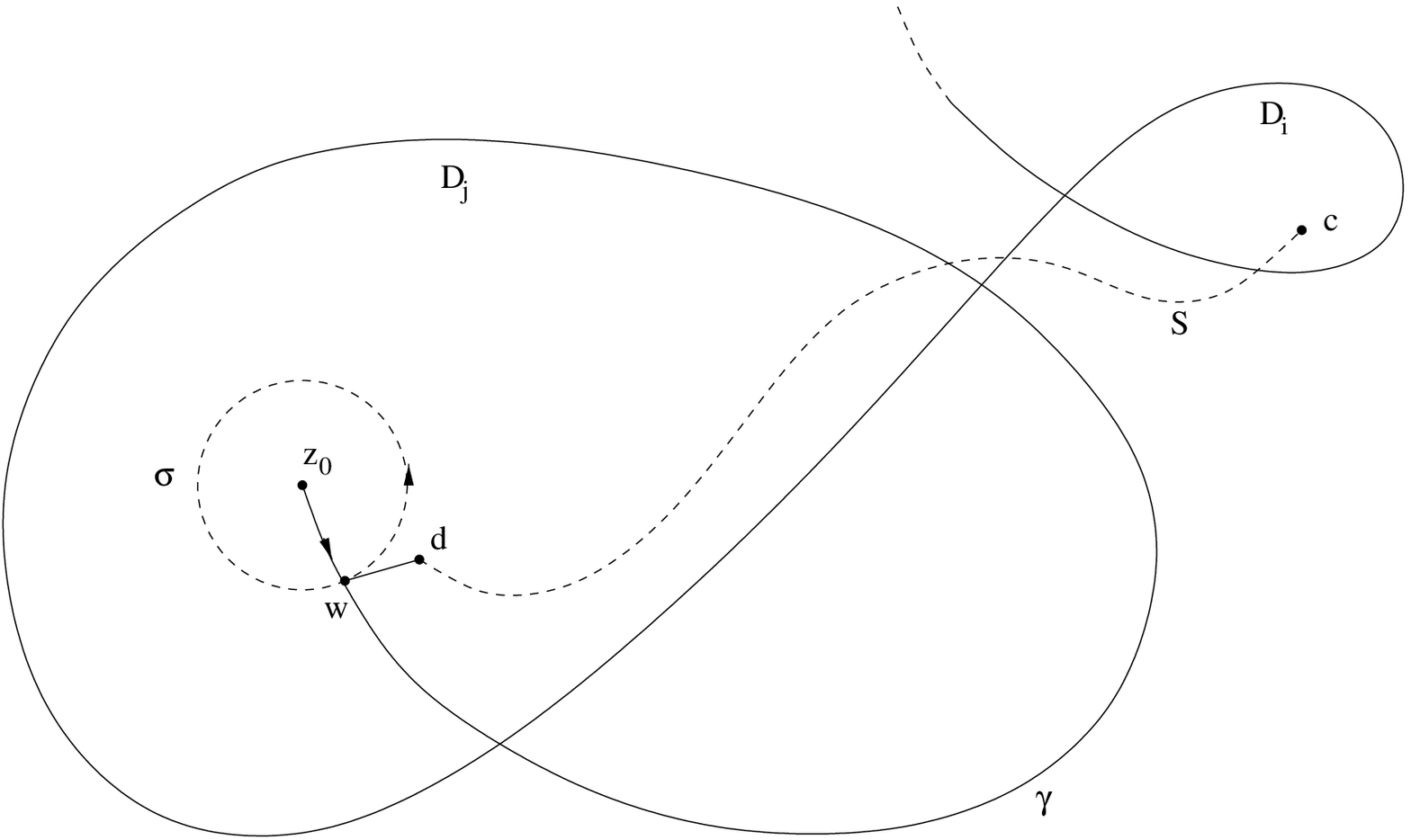}}

\centerline{Figure 5}
\medskip

\noindent{\bf Remark.}
If the endpoint $z_0$ of $\gamma$ belongs to the initial domain $D_i$ and $I_i$ is
rational then Proposition 4.8 implies that $g_1\equiv 0$ on $\gamma$.
This fact follows immediately also from the possibility to rich both sides of $\gamma$
near $z_0$ from the same point $c$ in $D_i$. Indeed, $g_1$ on $\gamma$ is the difference
of the continuations of $I_i(t)$ on the two sides of $\gamma$. See Figure 6.

\medskip
\epsfxsize=8truecm
\centerline{\epsffile{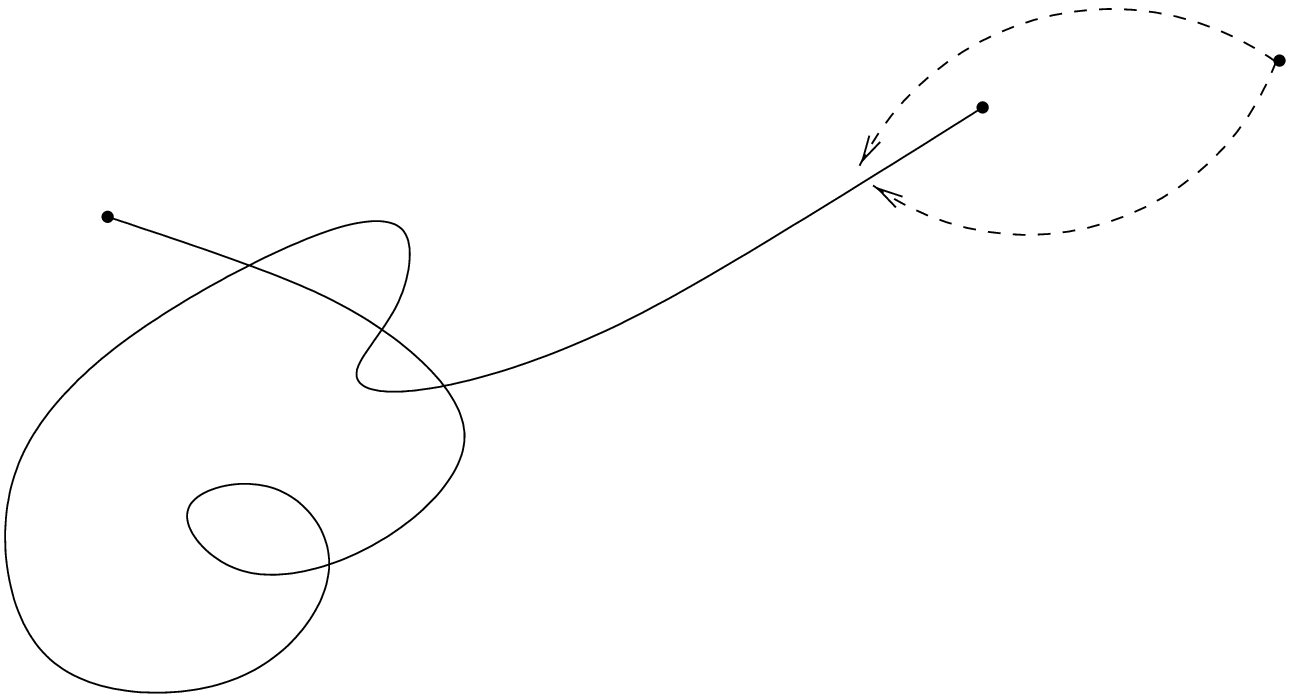}}

\centerline{Figure 6}
\medskip
Now let $z_0 \in \Sigma$ be a simple interior point of $\gamma$. According to the
definition of $\Sigma$ the function $g$ at $z_0$ 
may have a jump and/or a branching point. Denote by $g_0$ and $g_1$
the branches of $g$ on $\gamma$ before and after $z_0$, respectively. See Figure 7.

\medskip
\epsfxsize=8truecm
\centerline{\epsffile{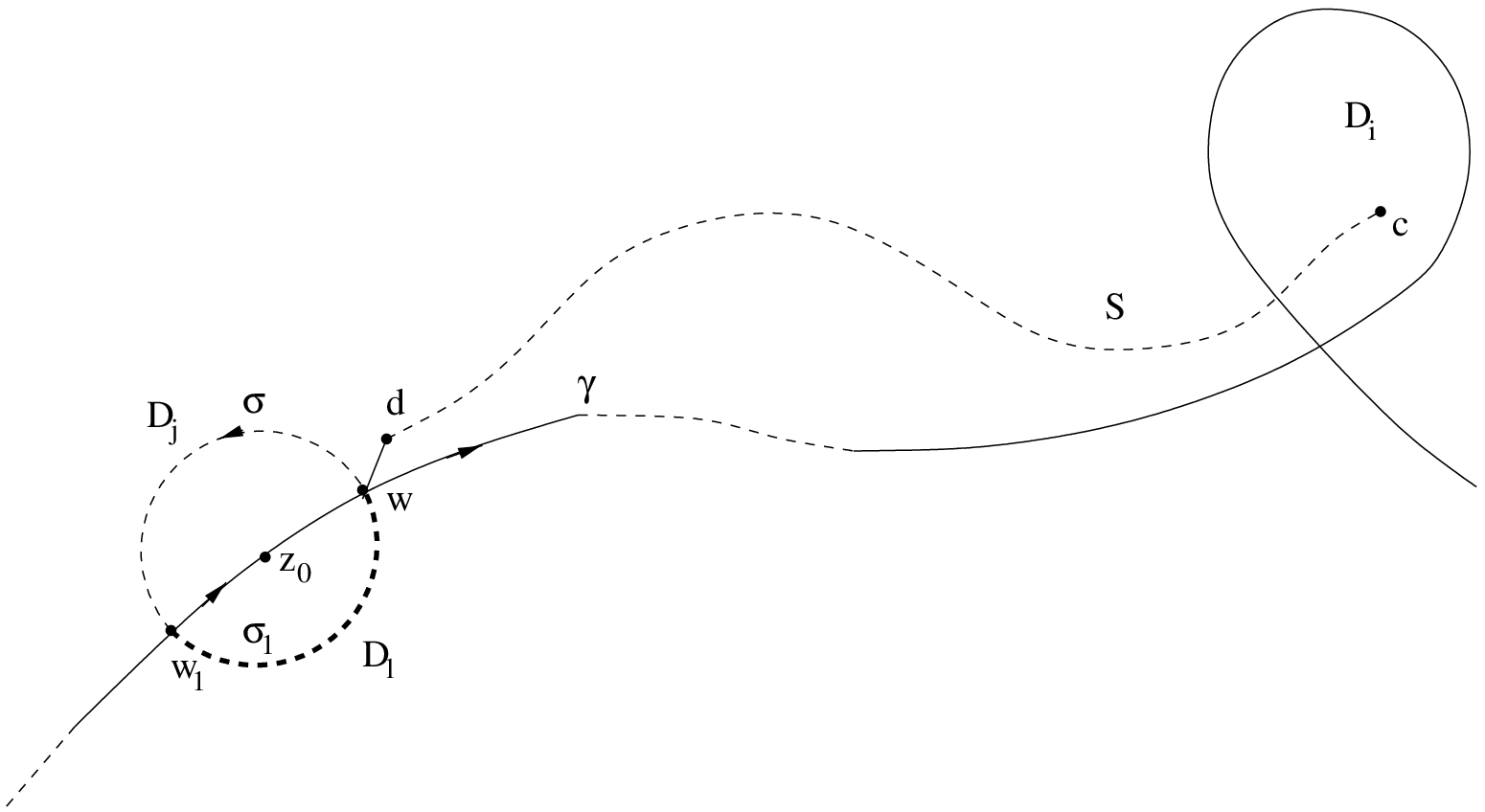}}

\centerline{Figure 7}
\medskip

Let the partition domains on the left and on the right of $\gamma$ near $z_0$ be  
$D_j$ and $D_l$, respectively. Let $S_{c,d}$ be an admissible curve with $c\in D_i$
and with $d\in D_j$ close to a certain point $w$ on 
the curve $\gamma$ after $z_0$ (see Figure 7).
The sum of branches $g(S_{c,d},\gamma)$ at $d$ we extend to the algebraic function $F$
defined in $D_j$. Denote by $\sigma$ a small closed loop going around $z_0$
in the counter-clockwise direction from the point $w$ to itself and let $\sigma_1$ be
the part of $\sigma$ in $D_l$ (so $\sigma_1$ goes in $D_l$ from a certain point
$w_1\in \gamma$ before $z_0$ to the point $w\in \gamma$ near $d$ (see Figure 7).

\noindent {\bf Proposition 4.9.} {\sl If $I_i(t)$ is a rational function then
in a neighborhood of $w$ we have $g_1-{\sigma_1}^*(g_0)= F -{\sigma^*}(F)$. 
In particular, if $D_j=D_i$ then $F\equiv 0$ and therefore 
$g_1\equiv {\sigma_1}^*(g_0)$ on $\gamma$.}

\noindent PROOF. We use the same auxiliary curves $\tilde S$ and $T=T_{c,w}$ as 
in the proof of Proposition 4.8. In our case application of Proposition 4.1 gives
$g(T,\gamma)={\sigma^*}(F)-{\sigma_1}^*(g_0)+g_1$ and from the equality
$g(T,\gamma)=g(\tilde S,\gamma)=F$ we get $g_1-{\sigma_1}^*(g_0)= F -{\sigma^*}(F)$.

\noindent {\bf Remark.} An important special case of the situation described in 
Proposition 4.9 occurs when $i=0$, $I_0(t)\equiv 0$ and $z_0$ is on the boundary of the
infinite domain $D_0$. Since $D_j=D_i=D_0$ the proposition gives
$g_1\equiv {\sigma_1}^*(g_0)$. In the specific case considered this follows directly
from the fact that both $g_0$ and $g_1$ are boundary values of the function $I_l$ in
$D_l$. See Figure 8. 

\medskip
\epsfxsize=8truecm
\centerline{\epsffile{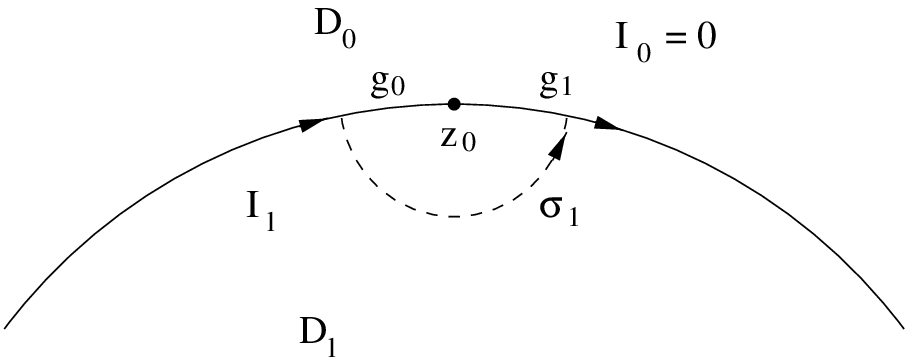}}

\centerline{Figure 8}
\medskip

Finally, let $z_0 \in \Sigma_1$ be a singular point of $g$ inside the domain $D_j$.
Let $\sigma$, $S_{c,d}$ and $F$ be as above, with $d$ near $z_0$. See Figure 9.

\noindent {\bf Proposition 4.10.} {\sl If $I_i(t)$ is a rational function then
$F -{\sigma^*}(F)=0$. In other words, $F$ does not ramify at $z_0$.}

\noindent PROOF. We use the original curve $S$ and the auxiliary curve $T=T_{c,d}$ 
obtained as $S$ followed by $\sigma$. As above, we denote $g(S,\gamma)$ by $F$. In our 
case $\sigma$ does not cross $\gamma$ and by Proposition 4.1 $g(T,\gamma)={\sigma^*}(F)$. 
Hence the equality $g(T,\gamma)=g(S,\gamma)=F$ proves the required result.

\medskip
\epsfxsize=8truecm
\centerline{\epsffile{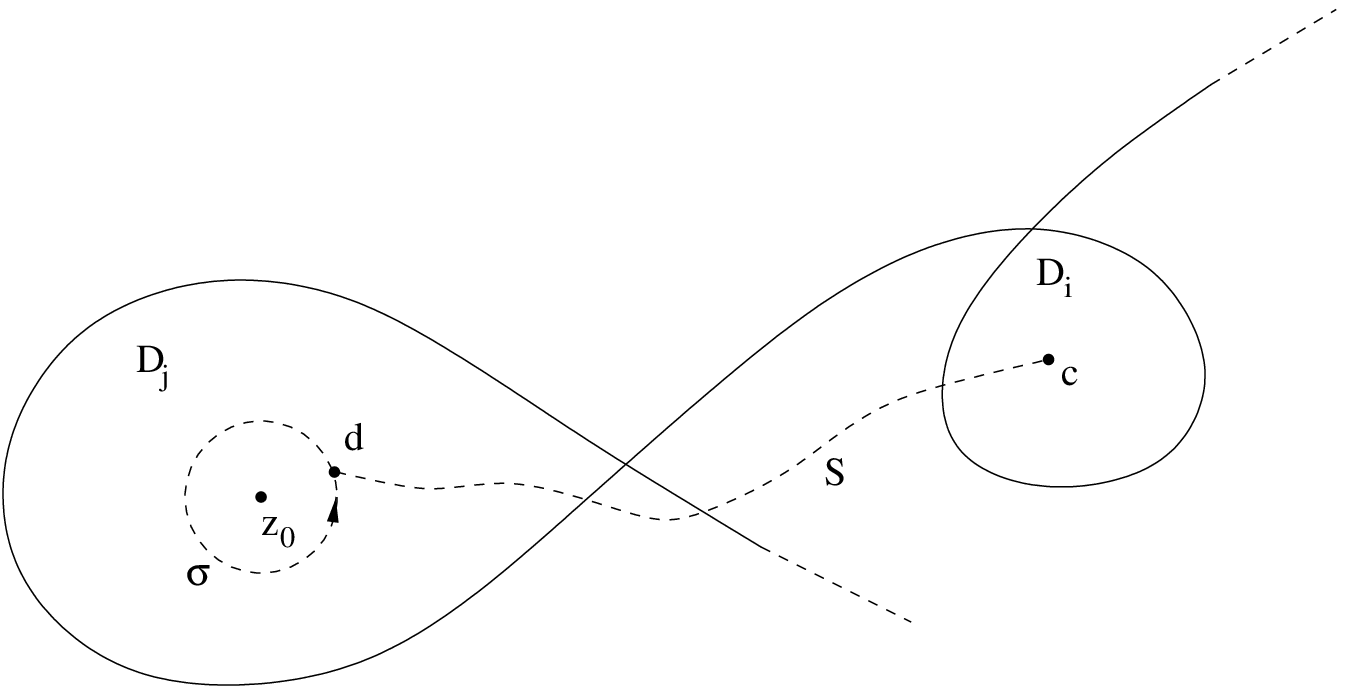}}

\centerline{Figure 9}
\medskip

\noindent {\bf Remark 1.} Of course, this result follows directly from Lemma 4.3:
$F= S^{*}(I_i(t))-I_j(t)$ with $I_j(t)$ regular at $z_0$ and $S^{*}(I_i(t))$ rational
(being the analytic continuation of the rational function $I_i$). In particular, the
algebraic function $F$ in $D_j$ may have only poles as singularities.

\noindent {\bf Remark 2.} The condition of a non-ramification of $F$ at $z_0$ can be
translated into the vanishing of certain sums of the branches of $g$ (as well as most
of the conditions in this section). Indeed, applying $\sigma^*$ to each of the summands 
$F_r$ in $F$ and equating the result to $F$ we get the vanishing of a certain sum of
the branches of $g$. 

\noindent {\bf Remark 3.} A natural question is under what conditions we get here
a {\it nontrivial} sum of the branches of $g$. In general, it would be important to find
the mutual dependencies between the local conditions of Propositions 4.8-4.10. 

Now the necessary conditions for rationality of $I_i(t)$ given by the three
propositions above turn out to be also sufficient. Let $z_k, \ k=1,\dots,N,$ be
all the points of $\Sigma$ (including the endpoints of $\gamma$) and of $\Sigma_1$,
taken in a certain order. Let us fix a point $c\in D_i$ and for each $z_k$ let us
fix an admissible curve $S^k$ leading to a neighborhood of $z_k$ and a small loop
$\sigma^k$ around $z_k$ (in each case as described in Propositions 4.8-4.10, 
respectively). 

\noindent {\bf Theorem 4.11} {\sl Assume that at each $z_k, \ k=1,\dots,N,$ (and
for the chosen $S^k$ and $\sigma^k$) the conclusion of Proposition 4.8 
(respectively, 4.9 or 4.10) is satisfied. Then the function $I_i(t)$ is rational.}

\noindent PROOF. Denote by $\tilde S^k={\tilde S^k}_{c,c}$ the loop following $S^k$
then $\sigma^k$ and then returning via $S^k$ in the opposite direction. If
the conclusion of the appropriate proposition 
above is satisfied then the sum of branches
along the loop $\tilde S^k$ is zero. Indeed, this conclusion expresses the fact
that the sum of branches along $S^k$ and along $S^k$ followed by $\sigma^k$ is the
same. But then the sum of branches along the loop $\tilde S^k$ is the same as
for $S^k$ passed forward and then back, and the last path is
homotopic to the constant one. Now the loops $\tilde S^k, \ k=1,\dots,N,$ generate the
fundamental group $\pi_1(U)$. By Proposition 4.1 if the sum of branches along the
loops $\tilde S^k, \ k=1,\dots,N,$ is zero the same is true for the products of these
loops. Therefore, the combinatorial monodromy of $I_i$ is trivial and by Theorem 4.7
the function $I_i(t)$ is rational.

\medskip

Finally we come to the conditions for the identical vanishing of $I_i(t)$. It is more
convenient to characterize first the property of $I_i(t)$ being identically constant.
We shall consider not the complex plane $\CC$ but the Riemann sphere $\CC P^1$.

\noindent{\bf Proposition 4.12.} {\sl $I_i(t)\equiv \rm Const$ in $D_i$ if and only if
the combinatorial monodromy of $I_i$ is trivial and the following additional
condition is satisfied: for any $j \ne i$ and for any admissible curve $S_{c,d}$
with $c\in D_i$ and $d\in D_j$ the sum of branches $F_j=g(S,\gamma)$ is regular 
in $D_j$. For $i=0$ (i.e. for the exterior domain $D_0$ and the Cauchy integral
$I_0(t)$ on it) the same conditions are necessary and sufficient for $I_0(t)\equiv 0$.}

\noindent PROOF. In one direction the result follows from Theorem 4.7
and Lemma 4.3. Indeed, $I_i(t)$ being identically constant is rational. Hence
the combinatorial monodromy of $I_i$ is trivial. Since the analytic continuation
$\hat I_i(t)$ is the same constant, the relation $F_j=\hat I_i(t)-I_j(t)$
of Lemma 4.3 with $I_j(t)$ regular implies regularity of $F_j$ in $D_j$.
In the opposite direction, if the combinatorial monodromy of $I_i$ is trivial then
$I_i(t)$ is rational by Theorem 4.7. Since $F_j$ are regular in $D_j$ for any $j \ne i$
the same relation of Lemma 4.3 implies regularity of $\hat I_i(t)$ in each $D_j$.
But in $D_i$ itself $\hat I_i(t)=I_i(t)$ is regular by definition. 
Hence the global rational
function $\hat I_i(t)$ is regular everywhere on $\CC P^1$ so it is constant. Since by
(1.1) $I_0$ is equal to zero at infinity $I_0$ is constant if and only if it is 
identically zero.

Of course, the condition of Proposition 4.12 essentially coincides with the classical
vanishing condition for the Cauchy-type integrals (i.e. that $g\vert\gamma$ bounds a 
holomorphic one-chain). This chain is provided by the sums of branches $F_j$: the
definition of the sum of branches shows immediately that $g\vert\gamma$  is a boundary
of $\sum F_j$.

\noindent {\bf Remark.} Regularity of the sum of branches $F_j=g(S,\gamma)$ in $D_j$
is equivalent to the cancellation of the negative Laurent terms of $F_j$ at each singular 
point of $g$ in $D_j$. This provides a set of local conditions at the singularities
of $g$ in $D_j$ expressed by certain linear equations on the branches of $g$. 
By Proposition 4.12 these conditions (together with the requirement that the combinatorial 
monodromy of $I_i$ be trivial) are equivalent to $I_i(t)$ being identically constant.
However, these conditions (in contrast to the ``sum of branches'' vanishing conditions)
are not {\it a priori} invariant under the monodromy action on $g$. It would be important
to understand the role of these ``regularity conditions'' and their relation to the
rest of the properties investigated above.

One way to explicitly verify conditions of Proposition 4.12 is to check the 
position of the poles of the algebraic function $g$.

\noindent{\bf Corollary 4.13.} {\sl If the combinatorial monodromy of $I_i$ is trivial
and all the poles of $g$ are in $D_i$ then $I_i(t)\equiv \rm Const$.}

\noindent PROOF. If the combinatorial monodromy of $I_i$ is trivial then for 
any $j \ne i$ the sum of branches $F_j=g(S,\gamma)$ is a univalued algebraic
function in $D_j$. Since $g$ has no poles in $D_j$ the same is true also for
$F_j$ (which is the sum of certain branches of $g$ in $D_j$). Hence $F_j$ is
regular in $D_j$.

\noindent{\bf Remark.} It is interesting to compare this result with the direct
computations for rational functions given in Section 2 above. For $g$ rational
and for any $\gamma$ the sum of branches along any $S$ starting in the exterior
domain $D_0$ and ending in some $D_j$ is equal to $\mu_j g$ (in notations of
Section 2). Consequently, the condition of Proposition 4.12 is satisfied if and
only if all the poles of $g$ belong to the ``outside" of $\gamma$ i.e. to the
domains $D_j$ with $ \mu_j$ = $0$. (In Corollary 2.4 above this result was obtained
by a direct computation). A natural question is whether it is possible to relax
accordingly the conditions of Corollary 4.13.  

\noindent{\bf Corollary 4.14.} {\sl Let all the poles of $g$ belong to the exterior
domain $D_0$. Then the complete analytic continuation $\hat I_0(t)$ cannot be
univalued on $U$ unless it is identical zero. In particular, $I_0(t)$ at infinity 
cannot be a nonzero germ of a polynomial, a rational or a meromorphic in $\CC$
function.}

\noindent PROOF. Theorem 4.4 implies that if $\hat I_0(t)$ is univalued on $U$
then the combinatorial monodromy of $I_0(t)$ is trivial. But since all the poles
of $g$ belong to the exterior domain $D_0$, by Corollary 4.13 this implies
that $I_0(t)\equiv 0$ at infinity.

\medskip

Let us consider now some examples. Returning to the Example 1 given in Section 3,
we see that at the endpoint $0$ of $\gamma = [0,1]$ the sum of local branches
of $g(z)= \sqrt{z}$ is zero while at the endpoint $1$ the sum of local branches
is not zero.
Consequently, $I(t)$ has a logarithmic ramification at $z=1$ while each branch of
$I(t)$ has a finite ramification (of order $2$) at $z=0$. On can see this also
from the explicit expression for $I(t)$ given in Section 3.

\noindent{\bf Example 2.}
Let $\gamma$ be the unit circle $S^1$ and let $g(z)=\sqrt{z}$, with $g(1)$ = $1$,
analytically continued along $S^1$ in a counter-clockwise direction. So $g$ has
a jump at $1$. The sum of the branches of $g$ across $S^1$  and along the curve $S$ 
shown by a dotted line on Figure 10-a, is twice the positive branch of the $\sqrt{z}$
at $z=2$. In the same way we can see that the sum of branches along $S$ repeated
$n$ times is $+2n \sqrt{z}$. Hence the combinatorial monodromy of $I(t)$ is infinite
in this example. The analytic continuation of $I(t)$ both from inside and from outside
$S^1$ has a logarithmic branching at $z=1$.

\noindent{\bf Example 3.}
For the same $g(z)=\sqrt{z}$ and a curve $\gamma$ going twice around 0,
$I_0(t)\equiv 0$ at $\infty$. Indeed, after the substitution $z = w^2$ we get
$I(t)={1 \over {\pi i}}\int_{S^1}{w^2 dw\over {w^2 - t}}$ and for $t$ near 
infinity the integrand is regular inside $S^1$. Hence in this example
the combinatorial monodromy is trivial. The sum of branches starting in the exterior
domain $D_0$ gives $+\sqrt{z}$ in the domain $D_1$ and the identical zero in 
the domain $D_2$ containing the origin (see Figure 10-b). Accordingly, the Cauchy
integral which coincides in our case with the sum of branches gives
$I_0(t) \equiv 0$, $I_1(t) = +\sqrt{t}$, and $I_2(t) \equiv 0$.

\noindent{\bf Example 4.}
Let $g(z)=\sqrt{z(z-1)}$ and let the curve $\gamma$ go around 0 and 1 in an
``$\infty$'' shape (see Figure 10-c). We continue along $\gamma$ the positive branch
of $g$ at $z = 2$. Here the continuation ``closes up" and $g$ does not
have jumps on $\gamma$. Nevertheless, the sum of branches across $\gamma$
and along a curve $S$ shown by a dotted line on Figure 10-c gives twice the positive
branch of the $\sqrt{z(z-1)}$ at $z=3$. To simplify the notations we denote the germ of 
this branch at $z=3$ by $a$. Since $a\ne 0$ the combinatorial monodromy of $I_0$
is not trivial and $I_0(t)$ is not rational. Here the obstruction is not in the jump
point of $g$ on $\gamma$ but rather in the branching of $g$ at $z=0$ and $z=1$.

Let us show that in this example $I(t)$ is in fact an algebraic function. To do this
consider a second loop $S'$ going from $z=3$ around the origin as shown on 
Figure 10-c. An easy computation gives $g(S',\gamma)=2a$. Applying Proposition 4.1
we obtain that the sum of branches along $SS$, $S'S$, $SS'$, and $S'S'$ is zero
as well as along other products of any two loops $S$, $S'$ or their inverses. 
(We use the notation for the product of the loops in the fundamental group of
$\CC \setminus \{0,1\}$ where the loops in the product are passed in the order they are 
written from the left to the right).
Remind that the monodromy of $g$ along both $S$ and $S'$ is given by a multiplication
by $-1$. Now application of Proposition 4.1 to any product of $S$ and $S'$ in the
fundamental group of $\CC \setminus \{0,1\}$ shows that
the sum of branches along this product is either $2a$ or zero. Hence the combinatorial 
monodromy of $I_0$ is of order $2$ and by Theorem 4.5 \ $I_0$ is algebraic (as well as the
other $I_j$) but not rational.

Let us stress once more that in this example the curve $\gamma$ is closed and the
function $g$ is regular at each point of $\gamma$. In particular, $g$ has no 
``jump point'' on $\gamma$. So the set $\Sigma$ contains exactly one point -- the double
point of $\gamma$. But since $g$ is regular near this points on both the crossing 
segments of the curve $\gamma$, any branch of $I_0(t)$ at this point is regular (by
Lemma 3.2 and Theorem 4.4). Therefore, in this example all the conditions of the
Propositions 4.8 and 4.9 are automatically satisfied. This is the condition of
Proposition 4.10 that is violated at the singular points $0$ and $1$ of $g$ and that
prevents $I_0(t)$ from being rational.

\medskip
\epsfxsize=14truecm
\centerline{\epsffile{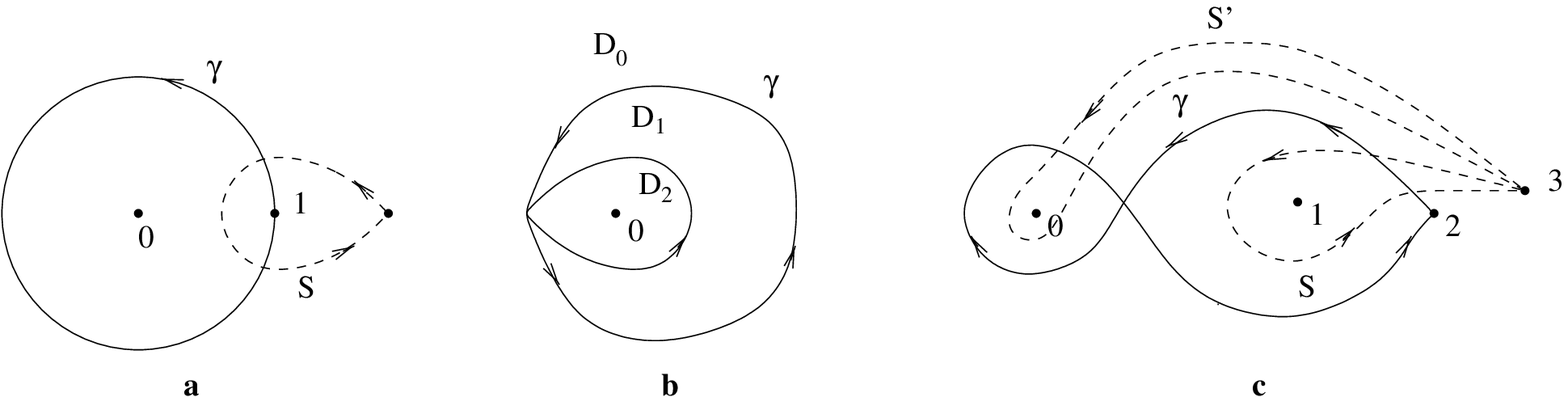}}
\smallskip
\centerline{Figure 10}

\medskip              

\noindent{\bf Example 5.} Let $\gamma$ be the interval $[-1,1]$ and let $g$ on
$[-1,1]$ be given by the positive branch of ${(1-z^2)}^{1\over r}$. We shall show
that for $r=2$ the Cauchy integral $I(t)$ is a non-rational algebraic function,
while for any integer $r \ge 3$ the function 
$I(t)$ is a non-algebraic locally algebroid function.
First of all let us notice that $g(z)$ satisfies the equation $g^r - (1-z^2)=0$.
Hence, the sum of all the branches of $g(z)$ is identically zero (being equal to the
$(r-1)$-th coefficient of the equation defining $g(z)$). Since the local
germs of $g(z)$ at $-1$ and $1$ contain {\it all its branches} the sum of the local
branches of $g(z)$ at each of its ramification points $-1$ and $1$ is zero. By
Theorem 3.4 and Theorem 4.4 this implies that for any leave of the full analytic 
continuation $\hat I(t)$ of the Cauchy integral $I(t)$ its branching
at $-1$ and $1$ is finite and the growth is bounded. Therefore all the singularities
of all the leaves of $\hat I(t)$ are algebroid and hence $\hat I(t)$ is locally algebroid.  

Now consider two loops $S$ and $S'$ going in a counter-clockwise direction 
around $-1$ and $1$, respectively, from a fixed point $c$ on the negative part
of the imaginary axis near the interval $[-1,1]$ (see Figure 11).

\medskip
\epsfxsize=8truecm
\centerline{\epsffile{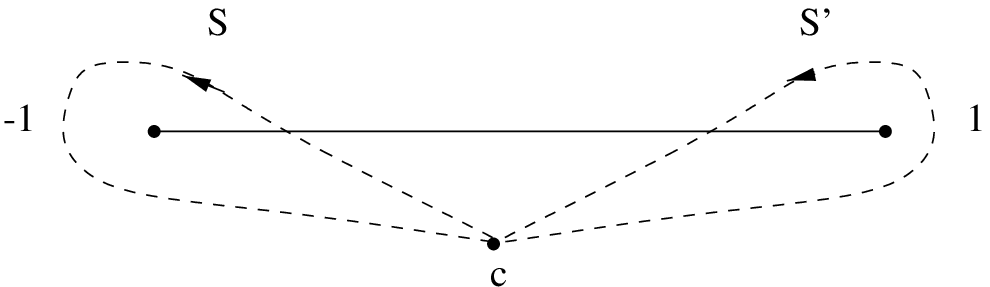}}
\smallskip
\centerline{Figure 11}

\medskip

The monodromy of $g$ along both $S$ and $S'$ is given by a multiplication by
$\epsilon = \exp ({{2\pi i}\over r})$. Denote by $a$ the germ at $c$ of the branch
of ${(1-z^2)}^{1\over r}$ that takes positive values on the interval $[-1,1]$.
On can easily see that the sum of branches along
$S$ and $S'$ is equal to $\epsilon a$ and $-a$, respectively. 

Now let us fix $r=2$. In this case $\epsilon = -1$ and the 
sum of branches gives $-a$ along both $S$ and $S'$. 
Proposition 4.1 then shows that for the sums of branches along the products $SS'$ and
$S^{-1}S'$ we have $g(SS',\gamma)=0$ and
$g(S^{-1}S',\gamma)=0$. Also the sums of branches along the products of any other two 
loops $S$, $S'$ or their inverses gives zero.
Using this fact we show exactly as in Example 4
above that the sum of branches along any product of $S$ and $S'$ in the
fundamental group of $\CC \setminus \{-1,1\}$ is either $-a$ or zero. Hence in
this case the combinatorial monodromy of $I_0$ is of order $2$ and by Theorem 4.5
the function $I_0(t)$ is a non-rational algebraic function.

Finally, to show that $I_0(t)$ is not algebraic for $r \ge 3$ consider the loop
$S^{-1}S'$ (see Figure 11). The monodromy of $g$ 
along this loop is trivial. An easy computation
shows that the sum of branches $g(S^{-1}S',\gamma)$ of $g$ along $S^{-1}S'$ is equal to 
$b= -a(1+\epsilon) \ne 0$. Let us formulate the last step of our computation as
a lemma.

\noindent{\bf Lemma 4.15.} {\sl Let $S=S_{c,c}$ be an admissible loop. Assume that
the monodromy of $g$ along $S$ is identical. If the sum of branches $g(S,\gamma)=a$
then the sum of branches $g(nS,\gamma)= na$ for any natural $n$.}

\noindent PROOF. This is a direct consequence of Proposition 4.1.

Applying Lemma 4.15, we see that the sum of branches along the loop $S^{-1}S'$
passed $n$ times is $nb$. Since $b\ne 0$ we conclude that the combinatorial monodromy
of $I$ is infinite and therefore by Theorem 4.5 \ $I(t)$ is not algebraic for $r \ge 3$.

\noindent{\bf Remark}. Example 5 presents a sequence of functions 
$$
I^r(t)=\int^1_{-1}{{(1-z^2)}^{1\over r} dz \over z-t}  \e(4.6) 
$$
with $I^2$ algebraic and $I^r$ transcendental but ``locally algebroid'' (according
to the definition above) for any natural $r \ge 3$. In fact, the functions $I^r$
also for $r \ge 3$ possess a number of remarkable properties which put them very close
to the algebraic ones. First of all, the full analytic 
continuation $\hat I^r$ of $I^r_0$ is 
regular in the domain $U=\CC \setminus \{-1,1,\infty\}$. Secondly, each of the infinite
number of the leaves of $\hat I^r$ has at the points $-1,1,\infty$ algebroid 
singularities with the branching of order $r$. Finally, the monodromy action of the
fundamental group ${\pi}_1(U)$ on the branches of $\hat I^r$ can be represents in a
relatively simple way via the combinatorial monodromy of $I^r$ (see Theorem 4.4 and
the computations in Example 5 above).

The functions with the same properties as $I^r$ appear as Cauchy Integrals of algebraic 
functions in many important cases. In general, by Theorem 4.4 this happens if and only
if the local sums of the branches of $g$ (given by Corollary 3.9) vanish at each of the 
jump points of $g$ on $\gamma$ (including the end points of $\gamma$). In particular,
this is the case if $\gamma$ is closed and $g$ is regular on $\gamma$. Especially important
example of this sort is provided by the rational Moment generating function (1.3) on the
closed path $\Gamma$ -- the case which corresponds to the classical Center-Focus problem.  

It turns out that the functions as above are closely related to certain Kleinian groups
and automorphic functions. Indeed, the ramification properties described above are rather
similar to that of the inverse to the factorization mapping by certain Kleinian groups.
Composing the functions as above with this factorization mapping we get a single-valued
meromorphic function whose behavior on the shifts of the fundamental domain is described
via their combinatorial monodromy of the original Cauchy integral. We plan to present 
separately the rigorous results in this direction.

We would like to thank S. Natanzon for pointing out the relations mentioned in the 
above remark.

\noindent{\bf Example 6.} Let us return once more to the case of $g$ itself  
being a rational function. The direct computations 
given in Section 2  show that in this case $I_0(t)$ as well as each $I_i(t)$
are rational functions. As it was mentioned in the remark above, for $g$ rational
and for any $\gamma$ the sum of branches along any $S$ starting in the exterior
domain $D_0$ and ending in some $D_j$ is equal to $\mu_j g$ (in notations of
Section 2). As we can expect, the combinatorial monodromy is trivial (although the
nontrivial sums of the branches of g do appear. These sums reflect just the geometry
of $\gamma$). So $I_0(t)$ is always rational. Now for each domain $D_j$ the sum of 
branches $F_j$ (obtained along any curve going from$D_0$ to $D_j$) is equal in
$D_j$ to $\mu_j g$. By Proposition 4.12 \ $I_0(t)\equiv 0$ if and only if the
functions $F_j$ are regular in $D_j$ for any $j\ne 0$. Now for $\mu_j= 0$ this
is automatic and for $\mu_j \ne 0$ the only way for $F_j$ not to have poles
in $D_j$ is that $g$ itself does not have poles. Consequently, $I_0(t)\equiv 0$
if and only if all the poles of $g$ belong to the ``outside" of $\gamma$ i.e. to
the domains $D_j$ with $ \mu_j$ = $0$.  

\noindent{\bf Example 7.}
In this example $\gamma$ is the union of the interval $[-1,1]$ and of the circle
$S_2$ of radius $2$ centered at the origin (Figure 12-a). 
The function $g$ on $[-1,1]$ is the
same as in Example 5 (for $r=2$) i.e. the positive branch of $\sqrt{(1-z^2)}$. 
The function $g$ on $S_2$ is the branch of ${1\over 2}\sqrt{(1-z^2)}$ taking the
values with the positive real part above the interval $[-1,1]$. Now the direct 
computation shows that the conclusion of Proposition 4.8 is satisfied at each
endpoint $-1$ and $1$ of $[-1,1]$ (which are also the only singular points of $g$). 
There are no jump points of $g$ on $\gamma$. Hence the conclusions of Proposition 4.9
and 4.10 are automatically satisfied. By Theorem 4.11 we obtain that $I_0(t)\equiv 0$.
We see also directly that $g$ on $\gamma$ bounds the chain $F$ equal to
the branch of ${1\over 2}\sqrt{(1-z^2)}$ taking the values with the positive real
part above the interval $[-1,1]$. (This statement can be interpreted also by replacing
the interval $[-1,1]$ with its two copies passed in the opposite directions and with
$g$ equal to the positive branch of ${1\over 2}\sqrt{(1-z^2)}$ on the ``upper''
interval and equal to the negative branch of ${1\over 2}\sqrt{(1-z^2)}$ on the 
``bottom'' one).

One can modify slightly the construction above and get a {\it connected} 
non-closed curve ${\gamma}'$ for which $g$ is not zero on ${\gamma}'$ but
$I_0(t)\equiv 0$. The curve ${\gamma}'$ starts at $-1$, goes along $[-1,1]$ till $0$, 
then goes up till $S_2$, encircles $S_2$ in the clock-wise direction, returns to the
interval $[-1,1]$ at $0$, and finally goes along $[-1,1]$ till $1$ (see Figure 12-b). 
The piecewise algebraic function $g$ on ${\gamma}'$
is defined as above on the parts of ${\gamma}'$ belonging to $[-1,1]$ and $S_2$ and
it is defined as a linear interpolation of the end values on the inserted parts.
As the integration on the two inserted intervals goes in the opposite directions, the
Cauchy integral of the extended $g$ on ${\gamma}'$ is the same as the original integral
on $\gamma$. An important question here is whether such an example (with 
$\gamma=P(\Gamma)$ non-closed, $g=Q(P^{-1})$, and $I_0(t)\equiv 0$) can appear in the
polynomial Moment problem.

\medskip
\epsfxsize=5truecm
\centerline{\epsffile{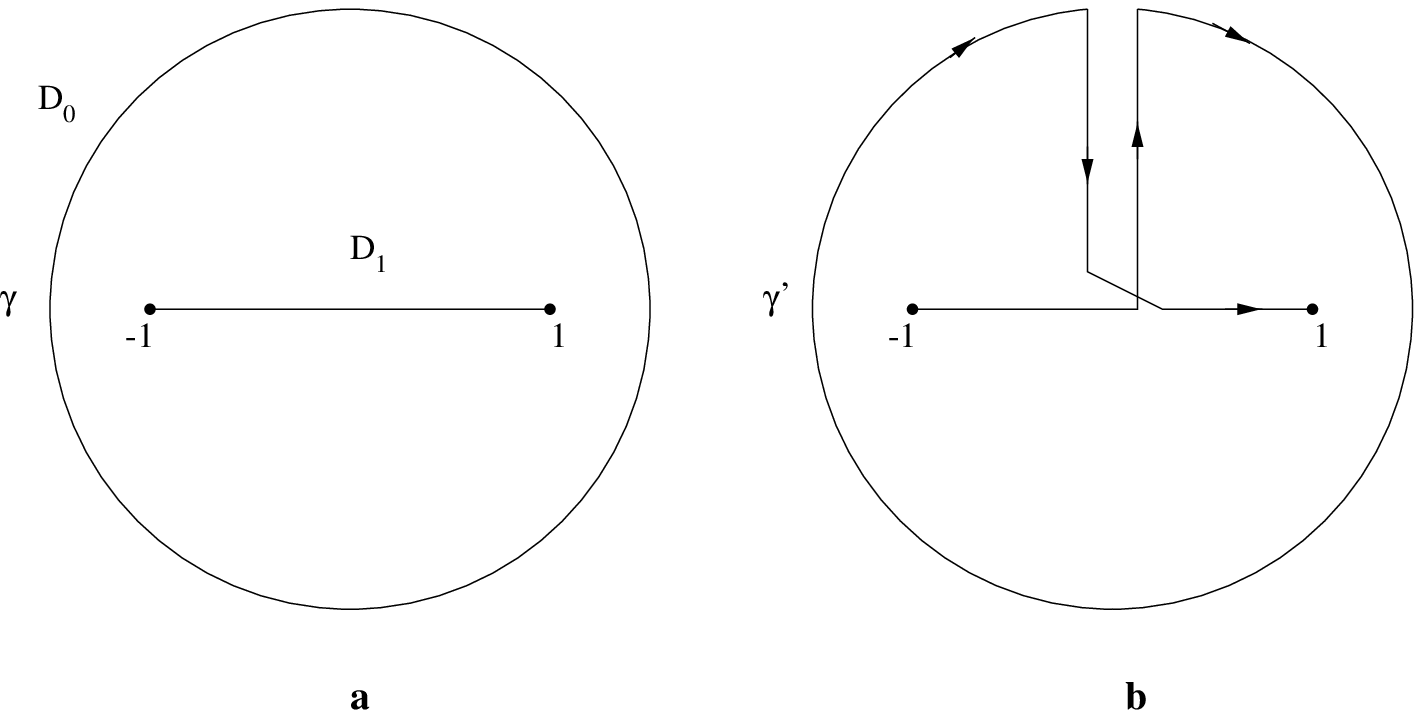}}
\smallskip
\centerline{Figure 12}

\medskip

\noindent{\bf Example 8.}
This rather unexpected example presents a curve $\gamma$ and an
algebraic function $g$ which does have a ``jump" on $\gamma$ and still  
$I(\gamma,g,t)\equiv 0$ for $t$ near $\infty$. As it was shown above 
under this condition the regular parts (or the normalized sums of the branches)
of $g$ on both sides of the jump point must coincide. However, in the situation
below the branches themselves on the two sides cannot be transformed one into
another by any local analytic continuation. Equivalently one can say that the
full local germs of $g$ on both sides of the jump point are different.

From the point of view of Proposition 4.9 above, we see that in this example
the difference between the branch $g_0$ of $g$ on one side of $0$ and the
local analytic continuation ${\sigma^*_1}(g_1)$ of the branch $g_1$ on the other 
side of $0$ is nonzero. However, it is compensated by the ``monodromy shift''  
$F_1 - {\sigma^*}(F_1)$ at $0$ of the Cauchy integral $I(t)=F_1$ (as defined
in the domain $D_1$ next to the exterior domain $D_0$).

This example seems to us to be quite unexpected. It is based on a recent
counterexample ([39]) to the ``Moment Composition conjecture'' (which
asserted that the vanishing of the moments (1.2) is equivalent to the Composition
condition (PCC). See Section 1.1.2 above).

\noindent{\bf Theorem 4.16.} {\sl Let $\gamma$ be the curve in $\CC$ shown
in Figure 13 below, and let $g(z)=Q(P^{-1}(z))$ with $Q(x)=T_2(x)+T_3(x)$
and $P(x)=T_6(x)$, where $T_n(x)=\cos({n\arccos}(x))$ is the $n$-th
Chebyshev polynomial. The function $g$ is analytically continued from
$0$ along $\gamma$ in the positive direction, starting with $Q(-\sqrt 3/2)$.
Then $I(\gamma,g,t)\equiv 0$
for $t$ near $\infty$, while the branches of $g$ on the two sides of $0\in
\gamma$ cannot be obtained from one another by any local analytic
continuation.}

\noindent PROOF. First of all, one can easily check that the curve $\gamma$
shown in Figure 13 is equal to $P(\Gamma)$, with $\Gamma$ obtained from  
the real interval $[-\sqrt 3/2, \sqrt 3/2]$ by a small shift (preserving the
ends) into the positive imaginary direction. Hence $I(t)$ is given by
$$
{2\pi i}I(t)=\int^b_a{Q(x)p(x)dx\over P(x)-t}\ \e(*),
$$
with $a = -\sqrt 3/2$, $b = \sqrt 3/2$ and $P$ and $Q$ as above.
Therefore $I(\gamma,g,t)\equiv 0$ for $t$ near $\infty$. Indeed, for $Q = T_2$
and for $Q = T_3$ the corresponding integrals vanish since the composition
condition (PCC) is satisfied: $T_6=T_2(T_3)=T_3(T_2)$ and all these three
polynomials take equal values at $a$ and $b$. Then for the sum $Q = T_2 + T_3$
the integral $I(t)$ vanishes by linearity in $Q$.
Now Claim 2 from [39] implies that $P$ and $Q$ do not have a common composition
factor. Therefore, by a well-known characterization of the composition factors
(see Lemma 5.1 and the ``Gluing Condition'' in Section 5 below) 
the branches of $g$ on the two sides of $0\in\gamma$ cannot be obtained from
one another by any local analytic continuation. This completes the proof.

\medskip
\epsfxsize=8truecm
\centerline{\epsffile{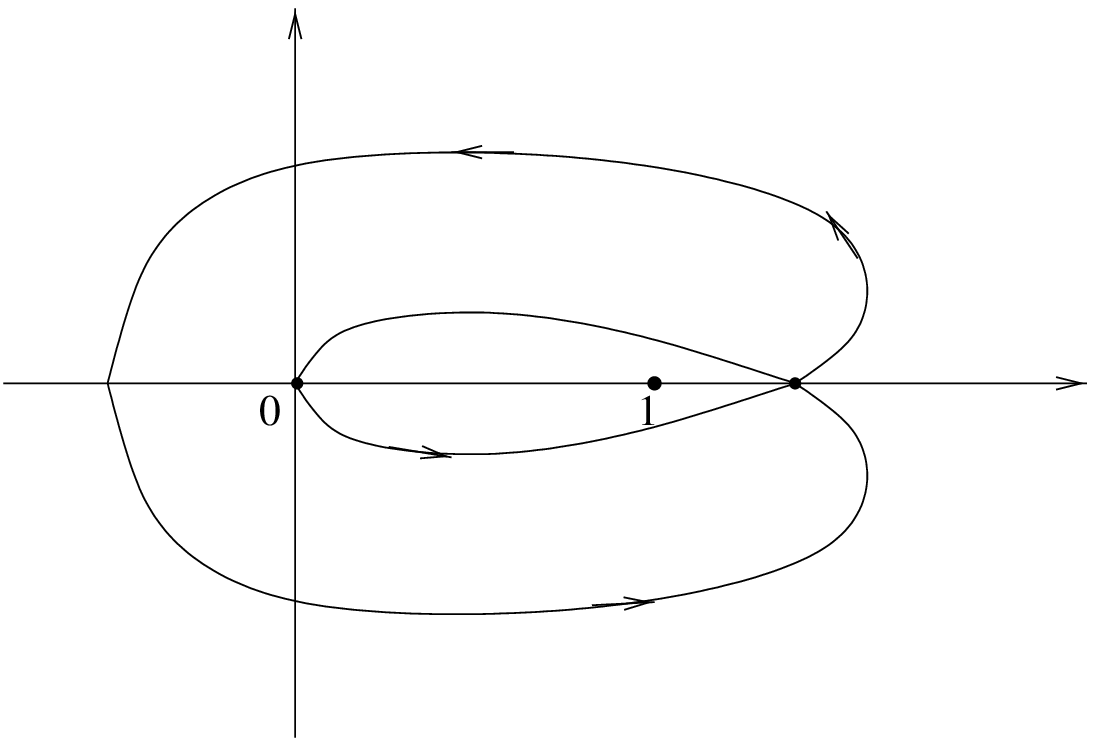}}

\centerline{Figure 13}
\medskip

Notice that in the example above $g$ = $g_1$ + $g_2$, with   
$g_1=T_2(T^{-1}_6)$ and $g_2=T_3(T^{-1}_6)$, respectively.
For each of $g_1$ and $g_2$ separately, their branches on the two
sides of $0\in\gamma$ can be obtained from one another by a local analytic 
continuation, but along two different local paths. In fact, these two paths are 
the two parts of the small circle around $0$. Each one joins the two pieces
of $\gamma$ on the two sides of $0\in\gamma$, but in different domains $D_j$.

Let us conclude this section with three remarks outlining some further
development directions for the tools introduced above.

\noindent{\bf Remark 1.} In the present paper we consider the ``sum of
branches'' of $g$ across $\gamma$ and along an auxiliary admissible curve
$S$. This notion can be further extended to the ``analytic continuation across
$\gamma$ and along $S$'' which can be applied to any analytic germ $u$ at the
starting point $c$ of $S$ (not necessarily to the germ defined by the Cauchy
integral). This operation consists in a continuation of $u$ along $S$ till 
the first crossing of $S$ with $\gamma$. At this crossing we add the germ of 
$g$ at the crossing point (with the sign corresponding to the orientation of the 
crossing). The resulting sum is analytically continued along $S$ till the next
crossing of $S$ with $\gamma$. At this crossing we add the germ of $g$ at the
crossing point, and so on till the endpoint $d$ of $S$. Lemma 4.3 shows that the
functions $I_i$ are obtained from one another via the analytic continuation across
$\gamma$ (and along any auxiliary admissible curve $S$). The notion of the
continuation across $\gamma$ allows one to extend the ``combinatorial monodromy''
(defined above as a {\it mapping} from the fundamental group $\pi_1(U)$ into
the germs at a certain point) into the {\it action} of $\pi_1(U)$ on these germs.

\noindent{\bf Remark 2.} In some cases it may be convenient to define the
continuation across $\gamma$ as a combinatorial process related to the functions
in the domains $D_j$ and not involving auxiliary curves $S$. For two adjacent
domains $D_i$ and $D_j$ separated by their common segment $\gamma_s$ of the curve
$\gamma$ and for a function $u_i$ in $D_i$ define its continuation across $\gamma$
into the domain $D_j$ as follows:
 
\item{a.} $u_i$ is analytically continued through $\gamma_s$ into a
certain neighborhood $\Omega$ of $\gamma_s$ in $D_j$.

\item{b.} An algebraic function $g_s$ in $\Omega$ (obtained by the
analytic continuation of the branch $g_s$ of $g$ on $\gamma_s$) is added to
$u_i$ (multiplied by $-1$ if the crossing orientation of $\gamma$ is negative).

\item{c.} $u_i+g_s$ is analytically extended from $\Omega$ to the entire
domain $D_j$. The function $u_j$ is equal to this continued function $u_i+g_s$.

\noindent Let us call the process consisting of the steps a, b, c an ``elementary
combinatorial continuation across $\gamma$'' and let us call any chain of subsequent
elementary continuations a ``combinatorial continuation across $\gamma$''. Assuming
that all the analytic continuations in the 
steps a, b, c are possible, the combinatorial
continuation across $\gamma$ can be applied to any analytic function defined in 
one of the domains $D_i$. In some cases this notion simplifies significantly the
description of the combinatorial monodromy.

\noindent{\bf Remark 3.} In the present paper we do not use a possibility to
deform $\gamma$ and to bring it to a certain simple standard form (without
affecting the Cauchy integral $I(t)$ near infinity). We also do not try to push
forward a unified algebraic framework where the combinatorial monodromy of $I$
can be naturally represented. Both these tools are important for the study of the
Cauchy-type integrals, and we plan to present them in detail separately, giving
here only a very short and informal outline of our approach.

In the process of the deformation of $\gamma$ it is natural in our setting
to keep fixed the ``jump'' points $\Sigma$ of $g$ on $\gamma$.
Besides this, $\gamma$ can be deformed in an arbitrary way
until it stays in $U$ (defined as above as $\CC$ with $\Sigma$ and all the
singularities of $g$ deleted). In fact, in the process of the deformation
certain crossings by $\gamma$ of the points of $\Sigma$ and of the singularities
of $g$ are also permitted (those which do not affect the local branch of $g$ on
$\gamma$ at the crossing). This allows us to bring $\gamma$ to a simple standard
form which we call the Diagram of $I(t)$. For the Cauchy integrals $I(t)$ coming
from the Polynomial moment problem, the Diagram of $I(t)$ can be computed quite
effectively using the methods of the ``Topological Theory of polynomials'' involving
a graphical representation of the monodromy group by means of the graphs obtained as 
the preimages of certain intervals, etc. Some applications of these methods to the 
Moment problem can be found in [40-41].

The sums of branches across $\gamma$ (and more generally, the result of the 
continuation across $\gamma$ as defined in the remarks 1 and 2 above) can be easily 
read off the Diagram. In these terms simple necessary and sufficient conditions for
the identical vanishing of the Cauchy integrals $I(t)$ near infinity can be given.

Now let us describe a global algebraic framework for representing the combinatorial
monodromy. It was shown above that the continuation across $\gamma$ and along a given
curve depends only on the homotopy class of $S$ in $U$. Therefore, it can be 
described through a certain action of the fundamental group $\pi_1(U)$ on the branches
of $g$. From the algebraic point of view, the object which naturally appears here
is the $Z(G)$-module $M$ consisting of all the ``formal'' finite sums of the branches
of the algebraic function $g$ at a given point $c\in U$. Here $G$ is the monodromy
group of $g$. There is also an ``evaluation homomorphism'' from $M$ into the 
$Z(G)$-module of the germs of the analytic functions at $c$.

In terms of the continuation across $\gamma$ a certain action $A$ of the fundamental 
group $\pi_1(U,c)$ on the module $M$ can be described. This action $A$ provides an
algebraic representation of the combinatorial monodromy of $I$. Consequently, the 
study of the algebraic structure of the module $M$ and of the action $A$ on it 
provides a natural and strong tool for a description of the global ramifications
of the Cauchy integral $I(t)$ and of its singularities. In particular, various
``vanishing sums of the branches'' relations given in this section can be naturally
expressed in this language. Indeed, all the vanishing sums of the branches of $g$
form a $Z(G)$-submodule $M_0$ of $M$ (the kernel of the evaluation homomorphism).
The study of the algebraic structure of $M_0$ allows us to represent the ``vanishing 
sums of the branches'' relations in a uniform way and to determine their mutual
dependence.

For the Cauchy integrals $I(t)$ coming from the Polynomial moment problem an
important algebraic information on $M_0\subset M$ and $A$ can be obtained
explicitly using the methods of the ``Topological Theory of polynomials''
mentioned above.

%
\bigskip
\noindent {\twelveb 5. Polynomial moments on an interval}

\noindent In this paper we investigate the polynomial moments (1.2)
$$
m_k=m_k(P,Q)=\int^a_b P^k(x)Q(x)p(x)dx  \e(5.1)
$$
and the Moment generating function (1.3)
$$
H(y)=\int^b_a{Q(x)p(x)dx\over t-P(x)}=\sum^{\infty}_{k=0}m_ky^{-k-1} .\e(5.2)
$$
Here $P(x)$ and $Q(x)$ are polynomials in $x\in\CC$,\ $a,b\in\CC$. As above,
we denote by $p(x)$ and $q(x)$ the derivatives of $P(x)$ and $Q(x)$,
respectively.

However, in most of the preceding papers [9-18,21,39-41,47,61] a slightly
different definition for the moments was used:
$$
\tilde m_k =\int^a_b P^k(x)q(x)dx , \ \ k \ge 0 .  \e(5.3)
$$
The setting of (5.1) and (5.2) is more convenient for the purposes of the 
present paper since it leads to the Cauchy integral (1.1) with the function
$g(z)=Q(P^{-1}(z))$ which does not have pole singularities at the finite points
of $\CC$ (see (5.5) below). (In the Cauchy integral for the
generating function $\tilde H(y)$ of the moments (5.3) the function $g$ is 
given by $g(z)= ({q\over p})(P^{-1}(z))$ and it may have poles at the finite
points of $\CC$, in particular, on $\gamma$.  

Let us show that the problems of the vanishing of the moments $m_k$ and 
$\tilde m_k$ are essentially the same. We do not assume a priori that $P(a)=P(b)$
or $Q(a)=Q(b)$. Notice also that adding a constant in $Q$ does not affect the
moments $\tilde m_k$, while it may affect the moments $m_k$.

\noindent {\bf Claim.} {\sl The vanishing of the moments $\tilde m_k$ 
implies $Q(a)=Q(b).$ Assuming that the primitive function $Q= \int q$
is chosen to satisfy $Q(b)=Q(a)=0$ all the moments $m_k$ also vanish.
In the opposite direction, if $Q(b)=Q(a)=0$ then the vanishing of the
moments $m_k$ implies that of $\tilde m_k$.}

\noindent {PROOF} Set, as above,
$$
H(t)=\sum_{k=0}^{\infty}m_kt^{-(k+1)}, \  \rm {and \ let}\
\tilde H(t)=\sum_{k=0}^{\infty}\tilde m_kt^{-(k+1)}.
$$ 
Then
$$
H(t)= \int_a^b {{Q(z)P^{\prime}(z)dz}\over {t-P(z)}}, \ \ \ \ \
\tilde H(t)= \int_a^b {{q(z)dz}\over {t-P(z)}}.
$$
We have: 
$$
{{dH(t)}\over {dt}}=
-\int_a^b {{Q(z)P^{\prime}(z)dz}\over {(t-P(z))^2}}=
-\int_a^bQ(z)d ({{1}\over {t-P(z)}})=$$
$$=-{{Q(z)}\over {t-P(z)}}\vert_a^b+\int_a^b{{q(z)dz}\over {t-P(z)}}=
{{Q(a)}\over {t-P(a)}}-{{Q(b)}\over {t-P(b)}}+\tilde H(t). \e(5.4)
$$
Suppose now that $\tilde m_i=0$ for all $i\geq 0.$ Then in particular 
$\tilde m_0= Q(b)-Q(a)=0$
and hence for any choice of $Q(z)=\int q(z)dz$ the equality
$Q(a)=Q(b)$ holds. Choose now $Q(z)$ such that $Q(a)=Q(b)=0.$ 
Then by (5.4) \ \ ${{dH(t)}\over {dt}}= {\tilde H(t)}$. Therefore 
$\tilde m_i=0,$ $i\geq 0,$ implies that $m_i=0,$ $i\geq 0.$ In the opposite
direction, under the assumption $Q(a)=Q(b)=0$ the formula (5.4) shows that
the vanishing of the moments $m_k$ implies that of $\tilde m_k$.

Let us return now to our original expressions (5.1) and (5.2).
A change of variables $P(x)=z$ brings (5.2) to the form
$$
H(y)= -{2\pi i}I(y)=-\int_{\gamma} {g(z)dz\over z-y}\ , \eqno(5.5)
$$
with $\gamma=P([a,b])$ and $g(z)=Q(P^{-1}(z))$. Notice that
the requirement that for $z \in \gamma$ the point $P^{-1}(z)$ belongs to $\Gamma$
defines the branch of $P^{-1}$ on $\gamma$ uniquely at any simple point of $z\in \gamma$.
Therefore the above expression $g(z)=Q(P^{-1}(z))$ correctly defines a
piecewise-algebraic function $g$ on $\gamma$ which satisfies all the requirements 
of (1.1). 

Notice also that for $\vert t\vert \gg 1$ we can take in (5.2) any integration
path joining $a$ and $b$. We shall use this later.

To investigate the relation between the vanishing of the moments and the
Composition condition we need the following property of the algebraic function
$Q(P^{-1}(z))$ (see [21,39-41,44,45,47,52]): 

\noindent{\bf Lemma 5.1.} {\sl Let $P$ and $Q$ be two rational functions.
There exist rational $\tilde P$, $\tilde Q$, and $W$, $\deg W > 1$ such that 
$$
P(x)=\tilde P(W(x)), \ \ Q(x)=\tilde Q(W(x))
$$
if and only if in a certain simply-connected domain $\Omega$ not containing
critical values of $P$ for two different branches $P^{-1}_0(z)$ and $P^{-1}_1(z)$
the following equality is satisfied: $Q(P^{-1}_0(z)) \equiv Q(P^{-1}_1(z))$.

Under the additional assumption that $P(a)=P(b)=z_0$ there exist rational 
$\tilde P$, $\tilde Q$, and $W$ with
$$
P(x)=\tilde P(W(x)), \ \ Q(x)=\tilde Q(W(x))\ \rm and \ W(a)=W(b)
$$
if and only if the full local germs of $g_0=Q(P^{-1}_0)$ and $g_1=Q(P^{-1}_1)$ 
at $z_0$ coincide. Here the branches $P^{-1}_0$ and $P^{-1}_1$ of $P^{-1}$ take
at $z_0$ the values $a$ and $b$, respectively.}

We shall call the property of the coincidence of the two local germs
$g_0=Q(P^{-1}_0)$ and $g_1=Q(P^{-1}_1)$
at $z_0$ {\it the Gluing condition.} In the setting of the Cauchy integral (5.5) it
is equivalent to the fact that {\sl the branches $g_0(z)$ and $g_1(z)$ of 
$g(z)=Q(P^{-1}(z))$ on the two sides of $z_0$ on $\gamma$ (corresponding to the 
branches $P^{-1}_0$
and $P^{-1}_1$ of $P^{-1}$ on $\gamma$ taking at $z_0$ the values $a$ and $b$, 
respectively), can be obtained from one another by a local analytic continuation
near $z_0$.}

Below in this section we always assume that $P$ and $Q$ are polynomials and
$P(a)=P(b)=z_0$. In this case the Gluing condition is equivalent to the polynomial
Composition condition (PCC) of Section 1.1.2: 
$$
P(x)=\tilde P(W(x)), \ \ Q(x)=\tilde Q(W(x)), \ \ W(a)=W(b),
$$
with $\tilde P$, $\tilde Q$, and $W$ polynomials.

\medskip

As it was mentioned above, the ``sum of the branches'' condition provided
by Corollary 3.9 (and necessary already for algebraicity of $I(t)$) plays a
central role in the investigation of the Moment vanishing in [14-17,21,40-41,47,61]. 
The Gluing condition can be considered as a special case of the vanishing of the 
sum of the branches. Indeed, it corresponds to the case where there are exactly
two branches of $g$ in the sum (on each side of $z_0$ in $\gamma$) with the
signs $1$ and $-1$. For $a$ and $b$ -- regular points of $P$ this is the case
already for the initial relation produced by Corollary 3.9 (see [21]).
The main approach of [14-17] in the case of real polynomials and of [21,40-41,47,61]
in general case is to produce the Gluing condition starting with a more complicated
initial vanishing sum of the branches and using some additional algebraic (or analytic) 
considerations. A similar approach is used also in Section 6 below.

In the present section we concentrate on the consequences of the stronger (than the
vanishing of the sum of the branches in Corollary 3.9) condition provided by 
Proposition 4.9. This condition is necessary for the rationality of $I(t)$ and
consequently for its identical vanishing. It turns out that under some additional
geometric assumptions on $P$ this stronger condition leads directly to the Gluing
condition for $Q(P^{-1})$.  

In particular, we shall describe on this base some natural classes of ``definite 
polynomials''.
Let us remind that {\it definite} polynomials $P$ have been defined in Section 1.1.2 
above as those for which the vanishing of the one-sided moments (1.2) implies (and 
hence is equivalent to) the composition condition (PCC) for any $Q$. The role  
of definite polynomials in the local Center-Focus problem is shown in [9,18,61].
Even more important role they play in the global study of the Center equations
near infinity as presented in [17]. (At the end of Section 1.1.2 above we outline very
shortly these applications). In the present section we describe some classes of definite 
polynomials $P$ specified by the geometry of the curve $\gamma=P([a,b])$. In this
connection a simple geometric invariant of complex univariate
polynomials is introduced and some results and problems concerning this 
invariant are stated.

The starting result here is the following:

\noindent {\bf Theorem 5.2.} {\sl Let $P(x)$ be a complex polynomial,
$P(a)=P(b)=z_0$. Assume that there exists a path $\Gamma\subset \CC$
joining $a$ and $b$ such that $z_0$ is a simple point of $\gamma=P(\Gamma)$
and $z_0$ is on the boundary of the exterior domain $D_0$. Then for any 
polynomial $Q$ the moments $m_k(P,Q)$ defined by (1.2) vanish if and only
if the composition condition (PCC) is satisfied and hence $P(x)$ is definite
on $[a,b]$.}

\noindent PROOF. We use Proposition 4.9. In a special case where
$z_0$ is on the boundary of the exterior domain $D_0$, it implies that
$g_0 = {\sigma^*_1}(g_1)$ for $g_0$ and $g_1$ the branches of $g$ on
the two sides of $z_0$ in $\gamma$ and $\sigma_1$ a small path joining these two
sides. Hence $g_0$ and $g_1$ can be obtained from one another by a local 
analytic continuation near $z_0$ and the Gluing condition is satisfied.

\noindent Figure 14 shows two types of $\gamma$ with respect to
$z_0=P(a)=P(b)$. Another case of $P(a)=P(b)=z_0$ being ``strongly inside''
$\gamma=P(\Gamma)$ (and for any $\Gamma$ joining $a$ and $b$, as we shall
see below) is given by Example 8 of Section 4 above (see Figure 13).

\medskip
\epsfxsize=9truecm
\centerline{\epsffile{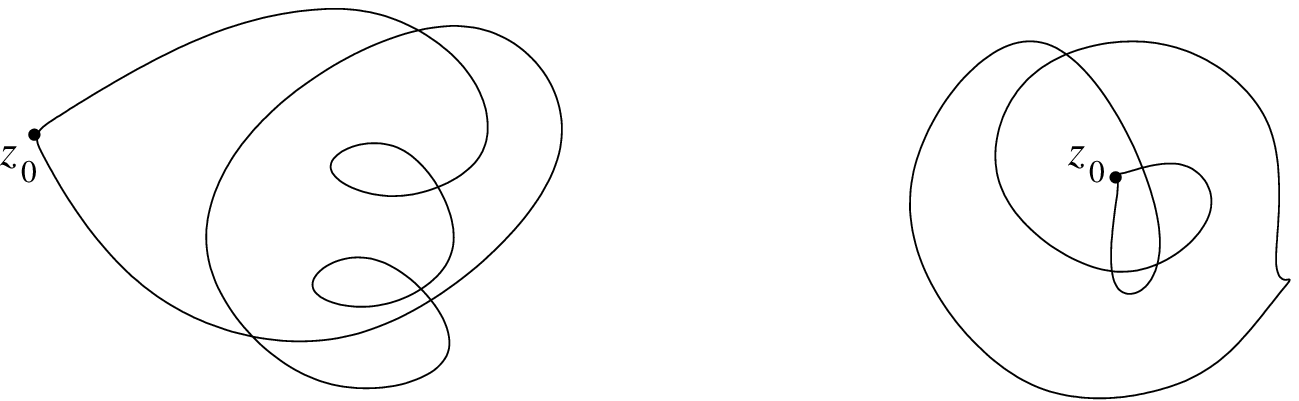}}
\smallskip
\centerline{Figure 14}

It would be important to characterize explicitly those polynomials $P$ 
(and couples $a,b\in\CC$, $z_0=P(a)=P(b)$) for which {\it there is a path $\Gamma$
joining $a$ and $b$ such that $z_0=P(a)=P(b)$ is a simple point of 
$\gamma=P(\Gamma)$ and this point is on the boundary of the exterior domain
$D_0$.} Let us call this {\it Property (E).}

In this context it is natural to generalize slightly this property and to introduce
a certain invariant of complex polynomials measuring how far from the
exterior domain $D_0$ lies the image of the endpoints. More accurately,
{\it the invariant $D(P,a,b)$ of a polynomial $P(x)$ with respect to 
$a,b, \ P(a)=P(b)=z_0$, is the minimal ``depth" of the point $z_0$ with respect
to the curve $\gamma=P(\Gamma)$ for various $\Gamma\subset \CC$ joining $a$ and $b$.
(The ``depth" here is the minimal number of crossings $\gamma$ necessary
to join $z_0$ to infinity.)}

Let us summarize some properties of $D(P,a,b)$.

\noindent{\bf 1.} For some $P$, $D(P,a,b)$ may be strictly positive.
For example, for $P=T_6$ - the 6-th Chebyshev polynomial, $a= -\sqrt 3/2$,
$b=\sqrt 3/2$ the invariant $D(P,a,b)$ is equal to one. 
Indeed, as it was shown in [39] for $P=T_6$
and $Q=T_2+T_3$ \ \ $I(t)\equiv 0$ at $\infty$ but (PCC) does not hold (see also
Example 8, Section 4). If there exists a path $\Gamma\subset \CC$ joining 
$a$ and $b$ for which $0=P(a)=P(b)$ lies on the boundary of the exterior domain 
$D_0$ with respect to $\gamma=P(\Gamma)$, then by Theorem 5.2 the Composition
condition (PCC) must be satisfied for any polynomial $Q$ for which the moments
(1.2) vanish. In particular, this must be true for $Q=T_2+T_3$ --
a contradiction. On the other hand, Figure 13 after Example 8 shows explicitly
the curve $\gamma=P(\Gamma)$ with the depth of $P(a)=P(b)=z_0$ equal to one. 

\noindent{\bf 2.} In many cases one can deform the path $\gamma=P(\Gamma)$ in
order to reduce the depth of $z_0$ in such a way that this deformation is covered
via $P$ by the corresponding deformation of $\Gamma$. One of the possible constructions
of such deformations is given below.

Let $S_{c,d}$ be a simple (without self-intersections) admissible 
curve with $c\in D_0$ and $d\in D_j$ where $D_j$
is one of the domains containing $z_0$ in its boundary. We assume also that $d$
is sufficiently close to $z_0$. Suppose that $S\cap \gamma =\{a_1,a_2, ... ,a_r\}$.
Let $a_i= P({\alpha}_i),$ $1\leq i \leq r.$ Notice that since $S$ is an admissible curve
then $a_i$, $1\leq i \leq r$, are simple points of $\gamma$. Therefore ${\alpha}_i$ are
defined in a unique way. Let $\{u_1, u_2, ... , u_r\}$ be the branches of $P^{-1}$ at $a_i$,
$1\leq i \leq r,$ taking at these points the values ${\alpha}_i$. 

Denote by $S_{a_i,d}$ the part of $S$ connecting $a_i$ and $d$ and consider 
the analytic continuations $h_i= S^{*}_{a_i,d}(u_i)$ of the germs $u_i$ along
$S$ to $d$. So $h_i$ represent certain branches of $P^{-1}$ at $d$ and they can be
analytically extended to $z_0$. Assume that the germs  $h_{i_1},\dots , h_{i_l}$ 
are regular at $z_0$ while the remaining $r-l$ germs $h_{i_{l+1}},\dots , h_{i_r}$
have singularities at $z_0$. We denote by $\nu (P,S,\gamma)$ the difference $r-l$.

Now assume that the polynomial $P$, the points $a,b\in\CC$ with $z_0=P(a)=P(b)$, and
the path $\Gamma$ joining $a$ and $b$ are given, and $\gamma=P(\Gamma)$.

\noindent{\bf Proposition 5.3} {\sl $D(P,a,b)$ does not 
exceed the minimum of $\nu (P,S,\gamma)$ over all the simple
admissible curves $S_{c,d}$ with $c\in D_0$ and $d\in D_j$ where $D_j$
is one of the domains containing $z_0$ in its boundary.}

\noindent PROOF. Fix an admissible curve $S_{c,d}$ as above. The deformation of 
$\gamma$ (covered via $P$ by the corresponding deformation of $\Gamma$) which
reduces the depth of $z_0$ to $\nu (P,S,\gamma)= r-l$ is constructed as follows:
consider one of the crossing points $a_{i_s}$, $s= 1, \dots ,l$, for which the germ
$h_{i_s}$ is regular at $z_0$.  Now we deform $\gamma$ along the curve $S_{a_{i_s},d}$ 
in such a way that the deformation is contained in a small neighborhood $\Omega$ of 
$S_{a_{i_s},d}$ and that in the final stage  the deformed curve $\gamma$ passes on 
another side of $z_0$ (see Figure 15). Since $z_0$ is a regular point of the branch
$h_{i_s}$ of $P^{-1}$ and since $S$ by assumptions does not contain other critical
values of $P$ then in fact the function $h_{i_s}$ is regular and univalued in a whole 
simply-connected neighborhood $\Omega$ of $S_{a_{i_s},d}$ (and it coincides with $u_{i_s}$
near the crossing point $a_{i_s}$). Moreover, since the curve $S_{a_{i_s},d}$ (and hence
the domain $\Omega$) do not have self-intersections, the function $h_{i_s}$ is one-to-one
on $\Omega$ (as an inverse function to $P$). Let us define the domain ${\Omega'}_{i_s}$
as the image ${h_{i_s}}(\Omega)$ (or as the preimage $P^{-1}(\Omega)$ for the appropriate
branch of $P^{-1}$). We obtain that $P$ restricted to ${\Omega'}_{i_s}$ is a regular 
covering over $\Omega$. Therefore the above deformation of $\gamma$ can be lifted by
$h_{i_s}=P^{-1}$ to the corresponding deformations of $\Gamma$. Repeating this operation
for each crossing point $a_{i_1},\dots , a_{i_l}$ we get a new curve $\gamma$ that
crosses $S$ only at the points $a_{i_{l+1}},\dots , a_{i_r}$. This completes the proof.

\medskip
\epsfxsize=5truecm
\centerline{\epsffile{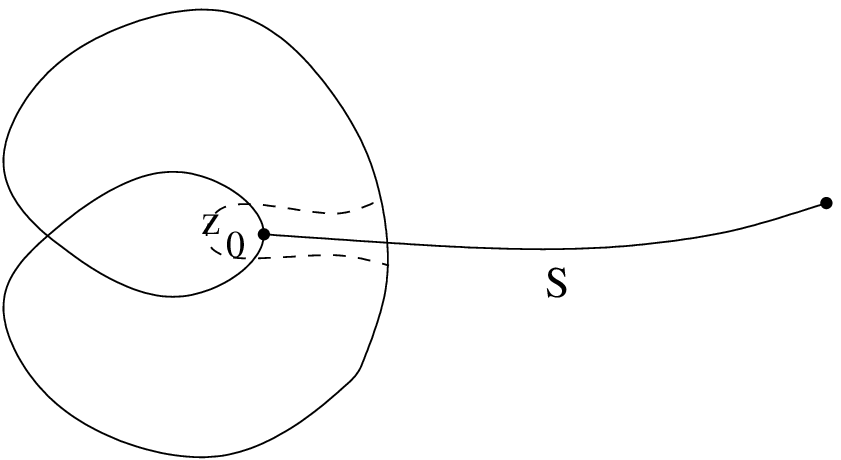}}
\smallskip
\centerline{Figure 15}

\noindent{\bf 3.} In general the deformations of $\gamma$ reducing the 
depth of $z_0=P(a)=P(b)$ to the minimum can be naturally analyzed in terms of the
``diagrams'' (shortly discussed in the concluding remarks in Section 4). As in
other cases this analysis becomes quite explicit via the methods of the ``Topological
Theory of polynomials'' mentioned above. We plan to present these results separately. 

Let us return now to definite polynomials. Using the approach of Proposition 5.3
in many cases one can show that a given polynomial $P$ possesses Property (E)
and hence it is definite. As the firest example let us give the following corollary:

\noindent{\bf Corollary 5.4.} {\sl If $z_0$ is a regular value of $P$ and 
$P(a)=P(b)=z_0$ then 
$D(P,a,b)=0$ and $P$ possesses Property (E). In particular, $P$ is a definite
polynomial on $[a,b]$.}

\noindent PROOF. If $z_0$ is a regular value of $P$ we get for any 
admissible curve $S_{c,d}$ as above $\nu (P,S,\gamma) = 0$. Proposition 5.3 now
implies that $D(P,a,b)= 0$ and hence $P$ possesses Property (E).

Corollary 5.4 follows also from the result of [21] that any $P$
is definite with respect to any two its regular points $a$ and $b$ with 
$P(a)=P(b)$.

Consider now reals polynomials on the real line. 
An application of Proposition 5.3 provides the following result:

\noindent{\bf Corollary 5.5.} {\sl Let $a,b\in\RR$ and let $P(x)$ be a real
polynomial with $P(a)=P(b)=0$. Assume that all the real zeroes $x_i$,
$i=1,\dots,s$, belonging to the open interval $(a,b)$ are simple. Then $P(x)$
possesses Property (E). In particular, it is definite on $[a,b]$.}

\noindent PROOF. Define $\Gamma$ by  shifting slightly the real interval $[a,b]$
into the upper half plane (and fixing $a$ and $b$). See Figure 16-a. (The corresponding
$\gamma$ is shown in Figure 16-c. Notice that $\gamma$ crosses the real axis near the 
the critical values $d_j$ of $P$, on the side of each $d_j$ which is determined by the
sign of the second derivative of $P$ at the corresponding critical point of $P$). 
We take as $S$ the part of the imaginary axes going
from $c= iD$, $D$ real and sufficiently big, to zero. Each 
crossing $a_{i_l}$ of $S$ and $\gamma$, $l=1, \dots, r \leq s$, corresponds to the  
point ${\alpha}_{i_l}$ of $\Gamma$ lying above and near one of the zeroes $x_{i_l}$ 
of $P$ where $P'(x_{i_l})>0$. (The parts of $\Gamma$ lying above and near the zeroes
$x_j$of $P$ with $P'(x_j)<0$ are mapped by $P$ into the pieces of $\gamma$ lying below
zero and hence $S$ does not cross these pieces of $\gamma$. Here we use the assumption
that {\it all the real zeroes} $x_i$ of $P(x)$ in the open interval $(a,b)$ are simple,
not only those where $P(x)$ changes sign from - to +).

Let $u_{i_l}$ be the branch of $P^{-1}$ taking at $a_{i_l}$ the value ${\alpha}_{i_l}$
and let $h_{i_l}$ be the analytic continuation of $u_{i_l}$ along $S$ from $a_{i_l}$ to
zero. The function $h_{i_l}$ maps $0$ into the root $x_{i_l}$ of $P$ and hence by the 
assumptions $h_{i_l}$ is regular at $0$ for $l=1, \dots, r$. We conclude that 
$\nu (P,S,\gamma) = 0$. Proposition 5.3 now implies that $D(P,a,b)= 0$ and hence $P$ 
possesses Property (E).   

\medskip
\epsfxsize=8truecm
\centerline{\epsffile{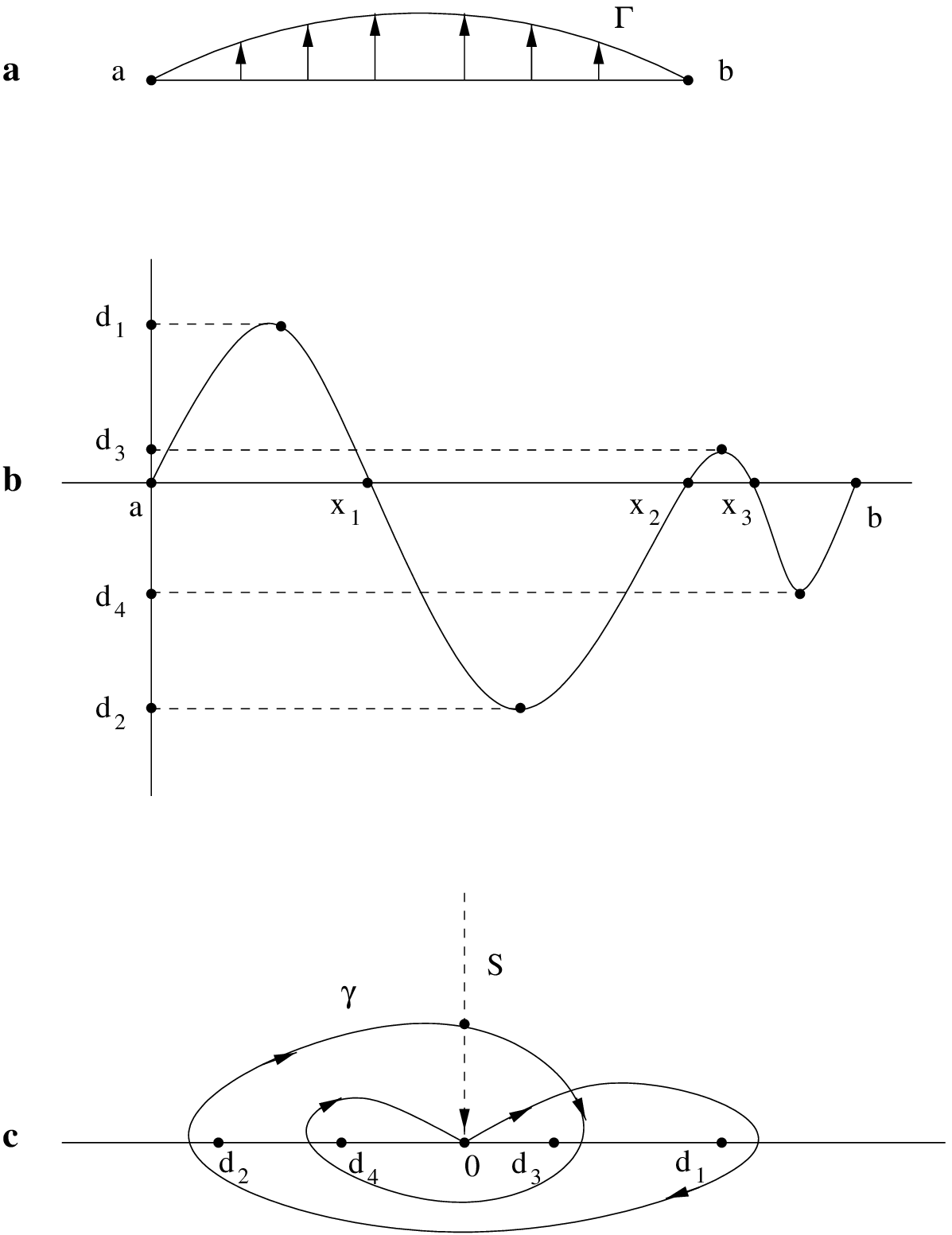}}
\smallskip
\centerline{Figure 16}
\medskip
\noindent{\bf Remark 1.} Apparently the result of Corollary 5.5 is not implied by the
other known characterizations of definite polynomials. Indeed, it involves conditions
only on the {\it real roots of $P(x)$} between $a$ and $b$, while most of the other
results work for general complex polynomials and cannot take an advantage of $P(x)$
being real and of a special properties of its real roots.

\noindent{\bf Remark 2.} It seams plausible that the assumptions of Corollary 5.5 may
be relaxed (if we want to show just that $P(x)$ is definite, without insisting on the
``geometric'' Property (E)). This is because in general the fact that if for a curve
$S$ as in Proposition 5.3 we have $\nu (P,S,\gamma)=0$ then $P$ is definite on $[a,b]$
follows also from Proposition 4.9. Indeed, the sum of branches $F$ along $S$ across
$\gamma$ for any integrand $g$ of the form $g= Q(P^{-1})$ is given by 
$$
g(S_{c,d},\gamma)= \sum_{l=1}^r sgn(a_{i_l})S^{*}_{a_{i_l},d}(Q(u_{i_l}))=
\sum_{l=1}^r sgn(a_{i_l})Q(h_{i_l}),
$$
where the branches $u_{i_l}$ and $h_{i_l}$ of $P^{-1}$ at $a_{i_l}$ and at d,
respectively, have been defined in the proof
of Proposition 5.3. Since the branches $h_{i_l}$ are regular at $z_0$ for $l=1,\dots,r$ 
we conclude that $F$ is regular for any polynomial $Q$. 
Remind now that Proposition 4.9 claims that if $I(t)\equiv 0$ at $\infty$
then the following equality between the branches $g_0$ and $g_1$ of $g= Q(P^{-1})$
on $\gamma$ before and after $z_0$ is satisfied:
$$
g_1 - {\sigma^*_1}(g_0) = F - {\sigma^*}(F),
$$
with $\sigma$ a local loop around $z_0$ and $\sigma_1$ is the part of $\sigma$
joining the segments of $\gamma$ before and after $z_0$ (see Figure 7).
An immediate consequence is that if $F$ is regular at $z_0$ then the branches
$g_0$ and $g_1$ can be obtained from one another by a local analytic continuation
along $\sigma_1$. Hence, the Gluing condition is satisfied which implies
the Composition condition on $P(x)$ and $Q$. 

However, the regularity of $F$ at $z_0$ may follow from a certain cancellation
effect and not just from the regularity of each of the branches $h_{i_l}$ of $P^{-1}$.
It would be interesting to find weaker conditions on $P$ providing regularity
of $F$ for any polynomial $Q$, besides the conditions given in Corollary 5.5
(and generally in Proposition 5.3). As it was mentioned above, these last conditions 
imply the
property (E) for $P$ which is presumably stronger than just the regularity of $F$ 
for any polynomial $Q$. A closely related conjecture is the following:

\noindent{\bf Conjecture} {\sl Let $P$ be a complex polynomial, $P(a)=P(b)=0$.
If all the roots of $P$ besides $a$ and $b$ are simple then $P$ is definite.}

\medskip

Let us now describe some classes of polynomials $P$ which possess the Property (E)
by ``geometric'' reasons. The following simple observation works in many specific
situations: 

\noindent{\bf Proposition 5.6.} {\sl Let $P$ be a complex polynomial
with $P(a)=P(b)=0$, $a,b \in \CC$, let $\Gamma$ be a piecewise-analytic curve
in $\CC$ joining $a$ and $b$ and let $\gamma = P(\Gamma)$. Assume that the open part
$\gamma \setminus 0$ is contained in an open domain $\Omega$ with piecewise-analytic
boundary and assume that $0$ belongs to the exterior boundary of $\Omega$.
Then $P$ possesses Property (E).}

\noindent PROOF. By the assumptions $0$ already belongs to the exterior boundary 
of $\Omega$ and hence to the boundary of the domain exterior to $\gamma$.
The only difficulty is that (as it happens in the examples below and in many other 
natural examples) for a 
specific given $\Gamma$ the curve $\gamma = P(\Gamma)$ may be not in general position.
So we need to perturb $\gamma$ to provide all its self-intersections transversal. If
we can do it in such a way that the point $0$ remains fixed and that for the perturbed 
curve $\gamma'$ the open part $\gamma' \setminus 0$ is still in $\Omega$, then $0$
belongs to the exterior boundary of $\gamma'$ and the result follows. We can restrict
the consideration to an arbitrarily small neighborhood of $0$. Indeed, outside such a
neighborhood $\gamma$ is at a certain positive distance from the boundary of $\Omega$
and hence it can be brought there to a general position by any sufficiently small 
generic smooth perturbation.

Near zero we use the assumption that $\gamma$ is a piecewise-analytic curve and that
the boundary of $\Omega$ is also piecewise-analytic. We obtain that $\gamma$ at $0$
has two local analytic branches $\gamma_0$ and $\gamma_1$ ``on the two sides of $0$''
(these branches may coincide with one another).
Since the boundary of $\Omega$ is piecewise-analytic and $\gamma_0 \setminus 0$ and
$\gamma_1 \setminus 0$ are contained in the open domain $\Omega$, we can perform an
analytic deformation of one of these branches (say, of $\gamma_0$) in such a way that 
the point $0$ remains fixed, the deformed curve $\gamma'_0 \setminus 0$ remains in
$\Omega$, and $\gamma'_0 \setminus 0$ and $\gamma_1 \setminus 0$ do not intersect
in a certain neighborhood of $0$ (see Figure 17). Indeed, if $\gamma_0$ and $\gamma_1$
do not coincide with one another then no perturbation is necessary. If 
$\gamma_0=\gamma_1$ then the required deformation of $\gamma_0$ can be achieved, for
example, by adding to it (in an appropriate coordinate system) an analytic germ with
a sufficiently high order of the vanishing at the origin.

Then we extend this deformation in a $C^{\infty}$
way to the rest of the curve $\gamma$ making all its self-intersections transversal.
If this deformation is small enough then $\gamma' \setminus 0$ is still in $\Omega$.
This completes the proof.

\medskip
\epsfxsize=6truecm
\centerline{\epsffile{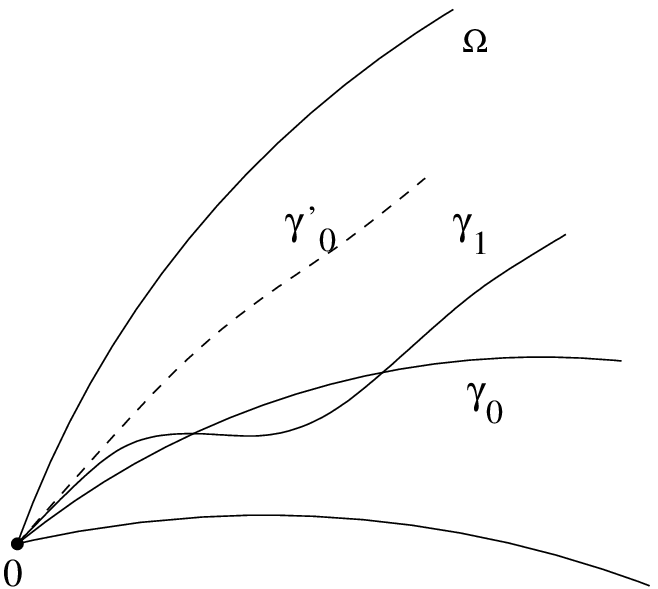}}
\smallskip
\centerline{Figure 17}
\medskip

As the first application of Proposition 5.6 we prove the following corollary
(which is also a special case of the Corollary 5.5 above):

\noindent{\bf Corollary 5.7.} {\sl Let $a,b\in\RR$ and let $P(x)$ be a real
polynomial with $P(a)=P(b)=0$ and $P(x)>0$ for $a<x<b$. Then
property (E) holds for $P$.}

\noindent PROOF. Apply Proposition 5.6 with $\Gamma = [a,b]$ and $\Omega$
an open cone $-\delta < Arg(z) < \delta$ for some $\delta > 0$. Notice that in this
situation $\gamma=P(\Gamma)$ is a real interval covered several times. So a perturbation
is indeed necessary to bring $\gamma$ into a general position. 

\noindent{\bf Remark.} In this specific case we can get $\Gamma$ with 
$\gamma = P(\Gamma)$ in general position and still inside $\Omega$ also by shifting 
slightly the real interval $[a,b]$ into the upper half plane while fixing $a$ and $b$.
See Figure 18. 

\medskip
\epsfxsize=15truecm
\centerline{\epsffile{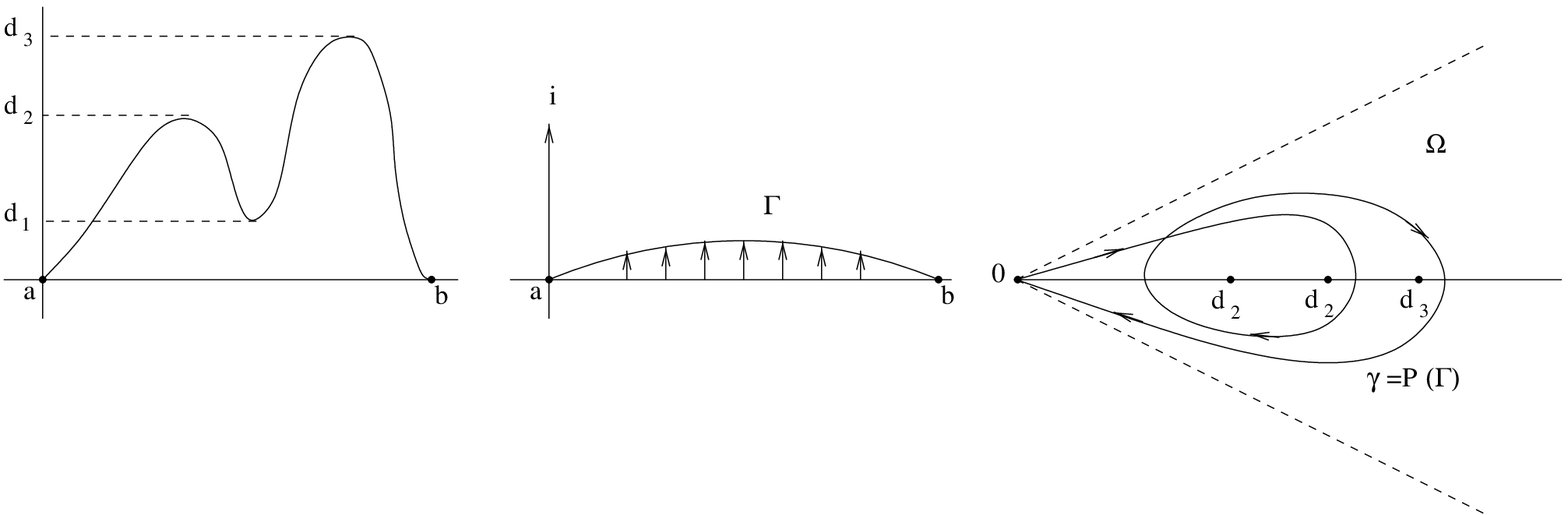}}
\smallskip
\centerline{Figure 18}
\medskip
The following corollary provides a class of definite polynomials which are characterized
directly by the geometry of their coefficients: these 
coefficients are assumed to belong to a certain
convex set in $\CC$ not containing zero. This assumption is not too restrictive from the
algebraic point of view, so the polynomials $P$ satisfying it may have zeroes of various  
multiplicities as well as various composition factorizations. In this sense the fact
of these polynomials $P$ being definite does not follow from the results of
[14-16,21,40-41,47] and from the rest of the results of the present paper.

\noindent{\bf Corollary 5.8.} {\sl Let $a,b\in\RR$, $0\leq a<b$ and let $P(x)=
(x-a)(b-x)P_1(x)$, with $P_1(x)=\sum^n_{k=0}a_kx^k$, $a_k\in\CC$. If the
convex hull $CH$ of the coefficients $a_k$ does not contain 0 then $P$
has Property (E).}

\noindent PROOF. Apply Proposition 5.6 with $\Gamma = [a,b]$ and $\Omega$ some
open cone $\alpha < Arg(z) < \beta$ containing the closed cone with the vertex
at $0\in \CC$ generated by $CH$. For any $x\in (a,b)$, $P(x)\in \kappa (CH)$, 
with $\kappa>0$.
Hence for $\gamma = P(\Gamma)$ the open part $\gamma \setminus 0$ is in $\Omega$.
See Figure 19.
 
\medskip
\epsfxsize=6truecm
\centerline{\epsffile{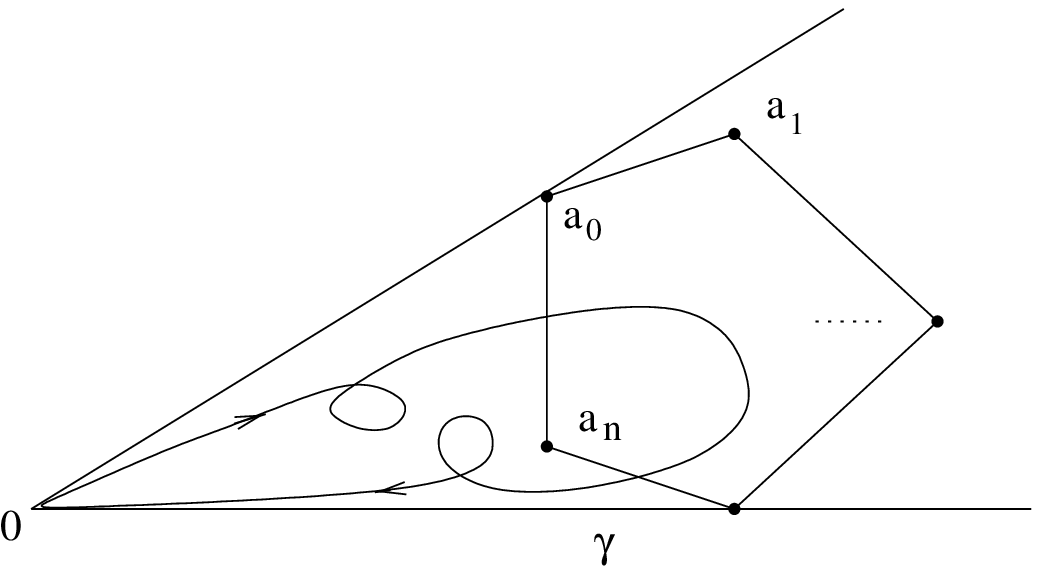}}
\smallskip
\centerline{Figure 19}
\medskip
\noindent {\bf Remark.} The results of the present section provide some explicit
classes of definite polynomials. Other classes have been described in 
[14-16,21,40-41,47] in terms of the monodromy group, indecomposability, or through 
certain
restrictions on their critical points. Yet another class of definite polynomials
can be produced via the algebraic methods of [15,16] (see [17]). All these classes
have almost no apparent intersections (besides the polynomials with regular zeroes
at $a$ and $b$, for the first two approaches). The only examples of the non-definite
polynomials we know at present are provided by [39]. In particular, the polynomial 
$P=T_6$ is not definite on $[-\sqrt 3/2, \sqrt 3/2]$. Each counterexample to the 
``Moment Composition conjecture'' given in [39] provides also an example of a 
non-definite polynomial $P$, and all these counterexamples are based on the composition
relations $A(B)=C(D)$ classified in [44,45] (see also [52]). Once more, the appearance
of these non-definite $P$ is not easy to relate with the properties used in 
[14-17,21,47]. It would be important to understand the nature of definite polynomials
and in particular to ``unify'' the approaches of the present paper, of [14-16,21],
and of [39-41]. The first steps in this direction are given in [9,17]. 
 
We conclude this section with a discussion of the following problem: is it possible
to relax the condition of the vanishing of all the one-sided moments (1.2) in the
Polynomial Moment problem (requiring, for example, the vanishing of only a part of them)
and still to get ultimately that all the $m_k$ vanish?
Some initial results in this direction can be derived from Corollary 4.14 above.
This corollary states that if all the poles of the integrand $g$ of the Cauchy integral
$I(t)$ given by (1.1) belong to the exterior domain $D_0$ then the complete analytic 
continuation $\hat I_0(t)$ of the ``exterior branch'' $I_0$ of $I$ cannot be univalued
unless it is identical zero. In particular, $I_0(t)$ at infinity cannot be a nonzero
germ of a polynomial, a rational or a meromorphic in $\CC$ function.  

The algebraic function $g=Q(P^{-1})$ on $\gamma$ that appears in the Cauchy
integral coming from the Polynomial Moment problem, does not have poles in the
finite part of $\CC$. Therefore by Corollary 4.14 the Moment generating function
$H(y)$ defined by (1.3) (or by (5.1)) cannot be rational unless it vanishes 
identically. This proves the following: 

\noindent {\bf Proposition 5.9.} {\sl If the moments $m_k=m_k(P,Q,\gamma)$
satisfy a linear recurrence relation
$$
m_j = \sum^N_1 \alpha_i m_{k-i}
$$
for each $j\ge N$ (in particular, if $m_j=0$ for each $j\ge N$)
then in fact $m_j=0$ for each $j\ge 0$.} 

The next interesting question in this direction is whether $H(y)$  may be a non-rational 
algebraic function. In the Examples of Section 4 where $I(t)$ is a non-rational 
algebraic function the Cauchy integral does not come from the Polynomial Moment 
problem (i.e. not from the integral (5.2)). As it was mentioned in the remark after
Example 5 in Section 4, for the closed integration path $\Gamma$ the rational Moment
generating function $H(y)$ is either algebraic or locally algebroid (i.e. has a finite
ramification and a polynomially bounded growth at each of its singular point). On the 
other hand, the following proposition shows that ``generically'' the polynomial 
Moment generating function on a {\it non-closed} interval cannot be locally algebroid. 

\noindent {\bf Proposition 5.10.} {\sl Let $P$ and $Q$ be polynomials and let $a$ and 
$b$ be regular points of $P$ with $P(a)=P(b)=z_0$. If one of the branches of the analytic 
continuation $\hat H(y)$ of the Moment generating function $H(y)$ defined by (5.1)
has a finite ramification at $z_0$ then in fact $H(y)\equiv 0$ near $\infty$. In
particular, $H(y)$ near $\infty$ cannot be a nonzero algebraic function.}

\noindent PROOF. By Theorem 4.4 if one of the branches 
of $\hat H(y)$ has a finite ramification 
at $z_0$ then this is true for all the branches, and in particular, for the one
represented by the expression in Theorem 3.7 of Section 3. Therefore, the logarithmic
term in this expression must vanish, i.e the regular parts of $g=Q(P^{-1})$ on the two
sides of $z_0$ in $\gamma=P([a,b])$ must coincide. Since $a$ and $b$ are regular points
of $P$ these regular parts are the corresponding regular branches of $g$ themselves.
We get a local coincidence of the branches of $g$ i.e. the Gluing condition. This 
implies Composition condition (PCC) that in turn implies the vanishing of all the
one-sided moments (1.2) and the identity $H(y)\equiv 0$ near $\infty$. 

\noindent {\bf Remark 1.} Proposition 5.10 can be formally restated as follows: it shows
that $H(y)$ cannot be a nonzero algebraic function near $\infty$ for any $P$
which is definite by the result of [21] (i.e. because $a$ and $b$ are regular points of 
$P$). Exactly in the same way we can show that $H(y)$ near $\infty$ 
cannot be a nonzero algebraic function for any $P$ which is definite by the results of
[40-41]. Indeed, in these papers the validity of the Composition condition (PCC) is
shown under certain restrictions on $P$, starting with the 
vanishing of a certain sum of the branches of $g$ which in turn
is implied just by the algebraicity of $H(y)$ (or even by the property of $H(y)$ to be
locally algebroid). It is not clear whether the same conclusion is true for the 
polynomials $P$ which are definite by the results of [15,16,47] or of the present 
section. 

\noindent {\bf Remark 2.} The condition of algebraicity of $H(y)$ can be expressed
in terms of its Taylor coefficients $m_k$, but not in a straightforward way (see
[24,48,49]). It would be interesting to find an analogue of Proposition 5.10
with the assumptions given explicitly in terms of the moments $m_k$.

\bigskip

\noindent {\twelveb 6. Rational Double moments}

\noindent In this section we investigate the ``Rational Double Moment problem''
on the non-closed curve. This problem consists in providing necessary and
sufficient conditions for the vanishing of the double moments (1.4)
$$
m_{i,j}=\int_{\gamma} P^i(x) Q^j(x)p(x) dx, \ \  i,j=0, 1, \ldots .
$$
We assume that $P(x)$, $Q(x)$ are rational functions and the integration path
$\Gamma$ is {\it non-closed}.
Remind that in the case of the closed integration curve the answer is given
by the classical result of Wermer and Harwey-Lawson: double moments vanish
if and only if the image curve $\delta=(P,Q)(\Gamma) \subset {\CC}^2$ of the path $\Gamma$
under $(P,Q)$ bounds a compact complex 1-chain in ${\CC}^2$ (see [2,23,30,59,60] and
Section 1.1.3 above). We show that on a non-closed curve $\gamma$ the vanishing of the 
double moments (and in fact just an algebraicity of the appropriate generating functions)
is equivalent to a certain composition factorization of the integrand functions
which ``closes up'' the integration path, combined together with the Wermer and 
Harwey-Lawson condition for their ``left factors''.

Next we show that under the additional assumption that the monodromy group
of $P(x)$ is doubly transitive the vanishing of the one-sided moments only   
implies a composition factorization which closes up the integration path. Moreover,
this composition factorization has a very special form: $Q(x)=\tilde Q(P(x))$. 

The results of this section generalize the results of [40].

\vskip 0.2cm

\noindent{\bf Theorem 6.1.} {\it
Let $P(x),Q(x)$ be rational functions and let
$\Gamma$ be a non-closed curve containing no poles of $P(x),Q(x)$ which
starts at the point $a$ and ends at the point $b.$ Suppose that
$$m_{i,j}=\int_{\Gamma} P^i(x)Q^j(x)p(x)\,dx=0$$
for $0\leq i\leq \infty,$ $1\leq j \leq d_a+d_b-1,$
where $d_a$ (resp.
$d_b$) is the multiplicity of the point $a$ (resp. $b$) with respect to
$P(x).$ Then there exist rational functions $\tilde P,$
$\tilde Q,$ $W$ such that $P(x)=\tilde P(W(x)),$
$Q(x)=\tilde Q(W(x)),$ and $W(a)=W(b).$}

\vskip 0.2cm

\noindent Note that if $a,b$ are not critical points of $P(z)$
(that is if $d_a=d_b=1$) then conditions of the theorem reduce to
the vanishing of single moments and therefore theorem 6.1
can be considered as a wide generalization of the result of C.
Christopher ([21]).

\vskip 0.2cm
 
\noindent{PROOF.} Suppose first that $P(a)=P(b).$
Let $U$ be a simply connected domain which
contains no critical values of $P(x)$ such that $P(a)=P(b)= z_0\in \partial
U.$ Denote by $P_{u_1}^{-1}(z),$ $P_{u_2}^{-1}(z), ... ,
P_{u_{d_a}}^{-1}(z)$ (resp. $P_{v_1}^{-1}(z),$ $P_{v_2}^{-1}(z), ... , 
P_{v_{d_b}}^{-1}(z)$) the branches of $P^{-1}(z)$ defined in $U$
which map
points close to $z_0$ to points close to $a$ (resp. $b$).
Then for any $j,$ $1\leq j \leq d_a+d_b-1,$
Corollary 3.9 applied to the function $g(z)=Q^j(P^{-1}(z))$
implies that
$$d_b\sum_{s=1}^{d_a} Q^j(P_{u_s}^{-1}(z)) \equiv d_a\sum_{s=1}^{d_b}  
Q^j(P_{v_s}^{-1}(z)).
\eqno(6.1)$$
Clearly, this equality holds also for $j=0$. Our assumption that the curve
$\Gamma$ is not closed (i.e. that $a\ne b$) implies that the branches of
$P^{-1}(z)$ on the two sides of (6.1) are different. Hence the correspondent 
branches of $Q^j(P^{-1}(z))$ are different and therefore (6.1) 
provides a {\it non-trivial} relation between the branches of $Q^j(P^{-1}(z))$.
(For $a=b$ the two sides of (6.1) are identically equal to one another).

Consider
a Vandermonde determinant $D=\parallel Q^j(p^{-1}_i(t)) \parallel,$
where $0\leq j \leq d_a+d_b-1$ and $i$ ranges the set of indices
$\{u_1,u_2, ... , u_{d_a}, v_1, ... , v_{d_b}\}.$
Since system (6.1) implies that $D=0$ we conclude that
$Q(P^{-1}_i(t)) \equiv Q(P^{-1}_j(t))$ for some $i\neq j,$ $1\leq i,j \leq n.$
(Here $n$ is the degree of the rational function $P$).
By Lemma 5.1 above the last condition is equivalent to the
condition that there exist rational functions $\tilde P,$
$\tilde Q,$ $W$ with $\deg W >1$ such that $P(x)=\tilde P(W(x)),$
$Q(x)=\tilde Q(W(x))$. Furthermore, without loss of generality we can suppose 
that $\tilde P$ and $\tilde Q$ have no nontrivial (of degree greater than
$1$) common right divisor in the composition algebra. Indeed, otherwise we 
compose $W$ with this common right divisor and get a new $W$ of a higher 
degree. The fact that $W$ satisfies additionally the equality $W(a)=W(b)$
we prove below.

Let us suppose now that $P(a)\neq P(b).$ In this case
instead of (6.1) Corollary 3.9 gives two systems
$$\sum_{s=1}^{d_a} Q^j(P_{u_s}^{-1}(z))=0,
\ \ \ \sum_{s=1}^{d_b}
Q^j(P_{v_s}^{-1}(z))=0, \ \ \ 1\leq j \leq d_a+d_b-1,
\eqno(6.2)$$ where
$P_{u_1}^{-1}(z),$ $P_{u_2}^{-1}(z), ... ,
P_{u_{d_a}}^{-1}(z)$ (resp. $P_{v_1}^{-1}(z),$ $P_{v_2}^{-1}(z), ... ,
P_{v_{d_b}}^{-1}(z)$) denote the branches of $P^{-1}(z)$ defined in
some neighborhood of $P(a)$ (resp. $P(b)$)
which map
points close to $P(a)$ (resp. $P(b)$) to points close to $a$ (resp. $b$). 
Now the same reasoning as above applied to the system
$$\sum_{s=1}^{d_a} Q^j(P_{u_s}^{-1}(z))=0, \ \ \ \ 1\leq j \leq d_a,$$
(taking into account that $d_a+d_b-1\geq d_a$)
shows that $Q(P^{-1}_i(t))=Q(P^{-1}_j(t))$ for two
different branches $P^{-1}_i(z),$ $P^{-1}_j(z)$
of $P^{-1}(z).$ Once more, this implies the existence of
rational functions $\tilde P,$
$\tilde Q,$ $W$ with $\deg W >1$ such that $P(x)=\tilde P(W(x)),$
$Q(x)=\tilde Q(W(x))$ and such
that $\tilde P$ and $\tilde Q$ have no nontrivial (of degree greater than
$1$) common right divisor in the composition algebra. 

Let us show that such a $W$ must satisfy $W(a)=W(b).$  
Indeed, otherwise after the change of variable $x\rightarrow w= W(x)$
we get
$$m_{i,j}=\int_{\delta}{\tilde P}^i(w){\tilde Q}^j(w){\tilde P}'(w)\,dw=0$$
for $0\leq i\leq \infty,$ $1\leq j \leq d_a+d_b-1,$
where $\delta = W(\Gamma)$ and $d_a$ (resp. $d_b$) is (as above) the 
multiplicity of the point $a$ (resp. $b$) with respect to $P(x).$
Taking into account that for any $c\in \CC$ the multiplicity
of $c$ with respect to $P(x)=\tilde P(W(x))$ is greater or  
equal than the multiplicity of $W(c)$ with respect to $\tilde P(w)$,
in the same way as above we would conclude that $\tilde P(w)=\bar
P(U(w))$, $\tilde Q(w)=\bar Q(U(w))$
for some rational functions $\bar P, \bar Q, U$ with
$\deg U >1$. This contradicts the assumption that $\tilde P,$ $\tilde Q$
have no common right divisor in the composition algebra.

\vskip 0.2cm

\noindent{\bf Corollary 6.2.} {\it
Let $P(x),Q(x)$ be rational functions and let
$\Gamma$ be a non-closed curve containing no poles of $P(x),Q(x)$ which
starts at the point $a$ and ends at the point $b.$ Then
$$m_{i,j}=\int_{\Gamma} P^i(x)Q^j(x)p(x)\,dx=0$$
for $0\leq i\leq \infty,$ $0\leq j \leq \infty$ if and only if
there exist rational functions $\tilde P,$
$\tilde Q,$ $W$ such that $P(x)=\tilde P(W(x)),$
$Q(x)=\tilde Q(W(x)),$ $W(a)=W(b),$ and all the poles of $\tilde P$ and
$\tilde Q$ lie on one side of the closed curve $\delta=W(\Gamma)$.}

\noindent{PROOF.} Sufficiency of these conditions follows from Theorem 1.1.3 
above, after we perform a change of variables $x\rightarrow w= W(x)$.
Necessity is obtained as follows: assuming that the moments $m_{i,j}$ vanish,
we apply Theorem 6.1 and get the factorization $P(x)=\tilde P(W(x)),$
$Q(x)=\tilde Q(W(x))$ with $W(a)=W(b)$. Performing a change of variables
$x\rightarrow w= W(x)$ we get the vanishing of the moments
$$\int_{\delta}{\tilde P}^i(w){\tilde Q}^j(w){\tilde P}'(w)\,dw$$
on the closed curve $\delta=W(\Gamma)$. Finally we apply Theorem 1.1.3.

\vskip 0.2cm

\noindent{\bf Corollary 6.3.} {\it If the moments $m_{i,j}$ vanish for
$0\leq i\leq \infty,$ $1\leq j \leq d_a+d_b-1$ then $P(a)=P(b), \ Q(a)=Q(b)$.} 

Notice that the results of Section 5 leave open the question whether the vanishing
of the {\it one-sided} moments implies $P(a)=P(b)$ and $Q(a)=Q(b)$.
It turns out that under the additional assumption that the monodromy group 
of $P(x)$ is doubly transitive the vanishing of the one-sided moments does
imply the equality $P(a)=P(b)$ as well as a composition factorization of a very
special form: $Q(x)=\tilde Q(P(x))$ (which of course closes up the integration 
path for $P(a)=P(b)$).

\vskip 0.2cm

\noindent{\bf Theorem 6.4.} {\it
Let $P(x),Q(x)$ be non-zero rational functions and let
$\Gamma$ be a non-closed curve containing no poles of $P(x),Q(x)$ which
starts at the point $a$ and ends at the point $b.$
Suppose that $$\int_{\Gamma} P^i(x)Q(x)p(x)\, d x=0$$
for $i\geq 0.$ If, additionally, the monodromy group of $P(x)$
is doubly transitive then the functions $P(x),Q(x)$ 
must satisfy $P(a)=P(b), \ Q(a)=Q(b)$ 
and there exists a rational function $\tilde Q$ such that $Q(x)=\tilde Q(P(x))$.}

\vskip 0.2cm

\noindent{PROOF.} 
Lemma 2 of [40] (see also [8]) states that if the monodromy group of $P(x)$ is
doubly transitive and if the branches $Q(P^{-1}_i(z))$ satisfy
$$
\sum^n_{i=1} a_i Q(P^{-1}_i(z)) \equiv 0
$$
for some $a_i \in \CC$ not all equal between themselves then there exists a
rational function $\tilde Q$ such that $Q(x)=\tilde Q(P(x))$. (Here as above
$n$ is the degree of the rational function $P$).

Now as in the proof of theorem 6.1, the vanishing of the one-sided
moments implies via Corollary 3.9 that the branches $Q(P^{-1}_i(z))$
are related either by relation (6.1) or by relation (6.2) (with $j=1$). Notice
that in each of these relations the coefficients are not all equal between 
themselves. This is immediate for (6.1). The only case where both the
sums in (6.2) contain all the branches of $Q(P^{-1}_i(z))$ is when $P(x)$ 
can be reduced to $x^n$ by the transformation $P(x) \rightarrow A(P(B(x)))$,
where $A, B$ are rational functions of the first degree. This possibility
is excluded by the assumption that the monodromy group of $P$ is doubly 
transitive. 

It remains to show that $P$ satisfies $P(a)=P(b)$. Let us perform a change of 
variable $x\rightarrow z=P(x)$. We get
$$
\int_{\gamma} z^i \tilde Q(z) dz = 0 \e(6.3)
$$
for $i\geq 0$, where $\gamma=P([a,b])$. If $P(a)\ne P(b)$ then the curve $\gamma$
is non-closed. In this case each of the relations (6.2) takes the form
$\tilde Q(z) \equiv 0$ (at the points $P(a)$ and $P(b)$, respectively). In other words,
the vanishing of the moments (6.3) for a non-closed curve $\gamma$ is possible only
for $\tilde Q(z) \equiv 0$. Since by the assumptions $Q(x)\ne 0$ also $\tilde Q$
cannot vanish identically. Hence $P(a)=P(b)$. This completes the proof of Theorem 6.4.

\vskip 0.2cm

\noindent{\bf Remark.} Corollary 3.9 requires only algebraicity of the Cauchy
integral $I(t)$ (and not necessarily its identical vanishing) to get the
vanishing of the local sum of the branches of $g$. Accordingly, we can replace
the assumption of the vanishing of the double (the one-sided) moments in Theorem 6.2
(Theorem 6.4, respectively) by the assumption of the algebraicity of the 
corresponding moment generating functions.

\bigskip

References

\item {1.} N.I. Akhiezer, {\it The classical moment problem},
Oliver and Boyd, Edinburgh, London, 1965.

\item{2.} H. Alexander, J. Wermer, {\it Several Complex Variables and
Banach Algebras}, Graduate Texts in Mathematics, {\bf 35}, Springer-Verlag.

\item {3.} M.A.M. Alwash, N.G. Lloyd,  {\it Non-autonomous equations related
to polynomial two-dimensional systems}, Proceedings of the Royal Society of Edinburgh,
{\bf 105A} (1987), 129--152.

\item {4.} M.A.M. Alwash, {\it On a condition for a centre of cubic
non-autonomous equations}, Proc. Royal Soc. of Edinburgh, {\bf 113A}
(1989), 289-291.

\item {5.} M.A.M. Alwash, {\it Periodic solutions of a quartic differential
equation and Groebner bases}, J. Comp. Appl. Math. {\bf 75} (1996), 67-76.

\item{6.} V.I. Arnold, {\it Problems on singularities and
dynamical systems}, In: Developments in Mathematics: The Moscow School,
eds. V.I. Arnold, M. Monastyrsky. London: Chapman \& Hall, 1993, 251-274.

\item{7.} V.I. Arnold, Yu. Il'yashenko, Ordinary differential
equations, {\it Encyclopedia of Mathematical Sciences} {\bf 1},
(Dynamical Systems - I), Springer, Berlin, 1988.

\item{8.} K. Girstmair, {\it Linear dependence of zeros
of polynomials and construction of primitive elements}, Manuscripta
Math. {\bf 39} (1982), no. 1, 81--97

\item{9.} M. Blinov, M. Briskin, Y. Yomdin, {\it Local Center Conditions
for Abel Equations}, preprint, 2003.

\item{10.} M. Blinov, N. Roytvarf, Y. Yomdin, {\it Center and Moment
Conditions for Rational Abel Equations}, Funct. Diff. Equations,
{\bf 10} (2003), No. 1-2, 95-106.

\item{11.} M. Blinov, Y. Yomdin , {\it Center and Composition Conditions for
 Abel Differential Equation, and rational curves},
Qualitative Theory of Dynamical Systems, {\bf 2} (2001), 111--127.

\item{12.} M. Blinov, {\it Center and Composition conditions for Abel equation},
Ph.D Thesis, the Weizmann Institute of Science, Rehovot (2002).

\item{13.} M. Briskin, Private communication.

\item{14.} M. Briskin, J.-P. Fran\c{c}oise, Y. Yomdin,
{\it Center conditions, compositions of polynomials and moments on
algebraic curves}, Ergodic Theory Dynam.
Systems {\bf 19} (1999), no. 5, 1201--1220.

\item{15.} M. Briskin, J.-P. Fran\c{c}oise, Y. Yomdin,
{\it Center conditions II: Parametric and model center problems }, Israel
J. Math. {\bf 118} (2000), 61--82. 

\item{16.} M. Briskin, J.-P. Fran\c{c}oise, Y. Yomdin,
{\it Center conditions III: Parametric and model center problems }, Israel
J. Math. {\bf 118} (2000), 83--108.

\item{17.} M. Briskin, N. Roytvarf, Y. Yomdin, {\it Center-Focus problem 
``at infinity'' for Abel equation, Moments and Compositions}, in preparation.

\item{18.} M. Briskin, Y. Yomdin, {\it Tangential Hilbert problem for Abel equation},
preprint, 2003.

\item{19.} A. Brudnyi, {\it On the Center Problem for Ordinary Differential
Equations}, preprint, 2003.

\item{20.} L. Cherkas,  {\it Number of limit cycles of an autonomous
second-order system}~, Differentsial'nye uravneniya  {\bf 12} (1976),
No.5, 944-946

\item{21.} C.J. Christopher, {\it Abel equations: composition conjectures and the
model problem}, Bull. Lond. Math. Soc.
{\bf 32} No. 3, 2000, 332-338. 

\item{22.} J. Devlin, N.G. Lloyd, and J.M. Pearson, {\it Cubic systems and Abel
equations}, J. Diff. Equations {\bf 147} (1998), 435-454.

\item{23.} P. Dolbeault, G. Henkin, {\it Chaines holomorphes de bord donne dans
$CP^n$}, Bull. Soc. Math. France 125 (1997), no. 3, 383-445.

\item{24.} H. Furstenberg, {\it Algebraic functions over finite field}, J. Algebra,
{\bf 7} (1967), 271-277.

\item{25.} H. A. Gasull, J. Llibre, {\it Limit cycles for a class of Abel
equations}, SIAM J. Math. Anal. {\bf 21} (1990), 1235-1244.

\item{26.} L. Gavrilov, {\it Petrov modules and zeros of Abelian integrals}, 
Bull. Sci. Math. {\bf 122} (1998), no. 8, 571-584.

\item{27.} L. Gavrilov, {\it Abelian integrals related to Morse polynomials
and perturbations of plane Hamiltonian vector fields}, Ann. Inst. Fourier {\bf 49}
(1999), no. 2, 611-652.

\item{28.} L. Gavrilov, I.D. Iliev, {\it Second order analysis in polynomially
perturbed reversible quadratic Hamiltonian systems}, Ergodic Theory Dynam. Systems
{\bf 20} (2000), no. 6, 1671-1686.

\item{29.} L. Gavrilov, I.D. Iliev, {\it Two-dimensional Fuchsian systems and the
Chebyshev property}, J. Differential Equations {\bf 191} (2003), no. 1, 105-120.

\item{30.} R. Harvey, B. Lawson, {\it On boundaries of complex analytic
varieties}, Ann. Math. {\bf 102} (1975), 233-290.

\item{31.} Seok Hur, {\it Composition conditions and center problem}, CRAS Paris,
{\bf 333, Ser. I} (2001), 779-784.

\item{32.} Yu. Ilyashenko, {\it Centennial history of Hilbert's 16th problem},
Bull. Amer. Math. Soc. (N.S.) {\bf 39} (2002), no.3, 301-354.

\item {33.} Yu. Il'yashenko, S. Yakovenko, {\it Counting real zeroes of analytic
functions satisfying linear ordinary differential equations}, J. Diff. Equations,
126, no.1 (1996), 87-105.

\item {34.} Yu. Il'yashenko, S. Yakovenko, {\it Double exponential estimate
for the number of zeroes of complete Abelian integrals}, Invent. math., 121,
613-650 (1995).

\item {35.} A.G. Khovanskii, {\it Real analytic varieties with the
finiteness property and complex Abelian integrals}, Funct. Anal. Appl.
{\bf 18} (1984), 119-127.

\item{36.} A. Lins Neto, {\it On the number of solutions of the equation
$x' = P(x,t)$ for which $x(0)=x(1)$}, Inventiones Math., {\bf 59}
(1980), 67-76.

\item{37.} N.G. Lloyd, {\it The number of periodic solutions of the equation
$z'= z^N + p_1(t)z^{N-1} + \dots + p_n(t)$}, Proc. London Math. Soc. {\bf 27}
(1973), 667-700.

\item{38.} N.I. Muskhelishvili, {\it Singular Integral Equations},
P. Noordhoff N.V., Groningen, 1953.

\item{39.} F. Pakovich, {\it A counterexample to the ``Composition Conjecture''}, 
Proc. AMS, 130, no. 12 (2002), 3747-3749.

\item{40.} F. Pakovich, {\it On the polynomial moment problem}, Math. Research Letters 
{\bf 10}, (2003), 401-410.

\item{41.} F. Pakovich, {\it On polynomials orthogonal to all powers 
of a Chebyshev polynomial on a segment}, to appear, Israel J. of Mathematics. 

\item{42.} J.M. Pearson, N.G. Lloyd, and C.J. Christopher, {\it Algorithmic derivation
of centre conditions}, SIAM Review, {\bf 38}, No. 4 (1996), 619-636. 

\item{43.} G.S. Petrov, {\it Complex zeros of an elliptic integral}, Funct. Anal. Appl.
{\bf 23} (1989), no. 2, 88-89. 

\item{44.} J. Ritt, {\it Prime and composite polynomials}, Trans. AMS. {\bf
23} (1922), 51-66.

\item{45.} J. Ritt, {\it Permutable rational functions}, Trans. AMS. {\bf
24} (1922), 399-488.

\item{46.} R. Roussarie, {\it Bifurcation of Planar Vector Fields and Hilbert's sixteenth
Problem}, Progress in Mathematics {\bf 164}, Birkhauser, Basel, 1998.  

\item{47.} N. Roytvarf, {\it Generalized moments, composition of polynomials and Bernstein
classes}, in ``Entire functions in modern analysis. B.Ya. Levin memorial
volume'', Isr. Math. Conf. Proc. {\bf 15}, 339-355 (2001).

\item{48.} K.V. Safonov, {\it On condition for algebraicity and rationality
of the sum of a power series}, Mat. Zametki {\bf 41} (1987), 325-332.

\item{49.} K.V. Safonov, {\it On Power Series of Algebraic and Rational functions},
J. of Math. Anal. and Appl. {\bf 243} (2000), 261-277.

\item{50.} S. Smale, {\it Mathematical problems for the next century},
Math. Intelligencer {\bf 20} (1998), no.2, 7-15.

\item{51.} S. Smale, {\it Dynamic retrospective: great problems, attempts that
failed}, Phys. D {\bf 51} (1991) 267-273.

\item{52.} A. Schinzel, {\it Polynomials with special regard to
reducibility}, Encyclopedia of Mathematics and Its Applications
{\bf 77}, Cambridge University Press, 2000.

\item{53.} D. Schlomiuk, {\it Algebraic particular integrals, integrability
and the problem of the center}~, Trans. AMS {\bf 338}
(1993), No.2, 799-841.

\item{54.} S. Shahshahani, {\it Periodic solutions of polynomial first order
differential equations}, Nonlinear Analysis {\bf 5} (1981), 157-165.

\item{55.} K.S. Sibirsky, {\it Introduction to the algebraic theory of invariants of
differential equations}, Nonlinear \ Science: theory and applications, Manchester
University Press, \ \ Manchester, 1988.

\item{56.} A.N. Varchenko, {\it Estimate of the number of zeroes of Abelian
integrals depending on parameters and limit cycles}, Funct. Anal. Appl. 18 (1984),
98-108.

\item{57.} S. Yakovenko, {\it On functions and curves defined by Ordinary
Differential Equations}, Proceedings of the Arnoldfest, Toronto, 1997,
Fields Institute Communications, Vol. 24, AMS, Providence, RI, (1999),
497-525.

\item{58.} Yang Lijun and Tang Yun, {\it Some new results on Abel Equations},
J. of Math. An. and Appl., {\bf 261} (2001), 100-112.

\item{59.} J. Wermer, {\it Function Rings and Riemann Surfaces}, Annals of
Math. {\bf 67} (1958), 45--71.

\item{60.} J. Wermer, {\it The hull of a curve in $C^n$}, Annals of Math.
{\bf 68} (1958), 550--561.

\item{61.} Y. Yomdin, {\it Center Problem for Abel Equation, Compositions of
Functions and Moment Conditions}, with the Addendum by F. Pakovich,
{\it Polynomial Moment Problem}, to appear in a special volume of the MMJ,
dedicated to the 65 birthday of V. I. Arnol'd.

\item{62.} H. Zoladek, {\it The problem of center for resonant singular points
of polynomial vector fields}, J. Diff. Equations {\bf 137} (1997), no. 1,
94-118.

\par\bye